\documentclass{amsart}
\usepackage{amsmath, amsthm, amsfonts, amssymb}
\usepackage{hyperref}
\hypersetup{
    colorlinks=true,
    linkcolor=black,
    citecolor=black,
    urlcolor=black,
    pdfborder={0 0 0}  
}
\usepackage{mathrsfs,comment}
\usepackage{color}
\usepackage{BOONDOX-cal}
\usepackage{cases}
\usepackage{pifont}
\usepackage{caption}
\usepackage{subcaption} 
\usepackage{tikz}
\usetikzlibrary{patterns} 
\usetikzlibrary{decorations.pathreplacing}  
\usepackage{tikzit}

\newtheorem{thm}{Theorem}[section]
\newtheorem{lem}{Lemma}[section]
\newtheorem{prop}{Proposition}[section]
\newtheorem{rem}{Remark}[section]

\newtheorem{con}{Conjecture}[section]

\allowdisplaybreaks[4]

\tikzstyle{Blue Dot}=[fill=blue, draw=blue, shape=circle]
\tikzstyle{Yellow Dot}=[fill=yellow, draw=yellow, shape=circle]
\tikzstyle{Black}=[fill=black, draw=black, shape=circle]
\tikzstyle{Red}=[fill=red, draw=red, shape=circle]
\tikzstyle{new style 0}=[fill=white, draw=black, shape=circle]

\tikzstyle{new edge style 0}=[-, dashed, dash pattern=on 2mm off 2mm]
\tikzstyle{new edge style 1}=[-, draw=red]
\tikzstyle{new edge style 2}=[-, fill=none, draw=blue]
\tikzstyle{new edge style 3}=[<->, dashed, dash pattern=on 2mm off 2mm]

\makeatletter
\@namedef{subjclassname@2020}{\textup{2020} Mathematics Subject Classification}
\makeatother

\begin{document}
\title[Sharp Poincar\'{e}-Wirtinger Inequalities on Finite Graphs]{Sharp Poincar\'{e}-Wirtinger Inequalities on Finite Graphs}
\author{Chenghao Xiang and Fayou Zhao$^\ast$ }
\address[Chenghao Xiang]{Department of Mathematics, Shanghai University, Shanghai
200444, P.R. China}
\email{chenghaoxiang@shu.edu.cn}
\address[Fayou Zhao]{Newtouch Center for Mathematics of Shanghai University,  Department of Mathematics, Shanghai University,
Shanghai 200444, P.R. China}
\email{fyzhao@shu.edu.cn}
\thanks{The research was supported by National Natural
Science Foundation of China (Grant No. 11971295).\\
$^\ast$ Corresponding author}
\subjclass[2020]{26D15, 05C50, 39A12, 05C12.}

\begin{abstract}
Let $G_n=(V_{G_n},E_{G_n})$ be the finite connected graph with $n\geq3$ vertices.
We study the sharp constant $C_{G_n,p}$ such that the Poincar\'{e}-Wirtinger inequality 
\[
\Vert f-m_f\Vert_p\leq C_{G_n,p} {\rm Var}_{p}^{G_n}(f)
\]
holds for every $f$ on $V$,
where ${\rm Var}_{p}^{G_n}$ stands for the $p$-variation on $G_n$ and $m_f$ stands for the arithmetic mean of $f$.
For complete graphs, the optimal constant is known for $p\in[1,3+\delta_n^1)\cup(3+\delta_n^2,+\infty)$, with $\delta_n^1=\frac{1}{2n^2\ln n}+O(1/n^3)$ and $\delta_n^2=\frac{2}{n}+O(1/n^2)$.
In this paper, we firstly extend the admissible range to all $p\in(0,\infty)$. Secondly, we obtain the sharp constant $C_{S_n,p}$ for the star graph $S_n$ in the range $p\in(0,+\infty)$ and give estimates for the sharp constant $C_{P_n,p}$ on path graphs when $p\in(0,+\infty)$. 
Thirdly, we characterize the extremizers on complete graphs and star graphs, and establish the existence of monotone extremizers on path graphs for $p\in(0,+\infty)$. For path graphs, we obtain explicit sharp constants for $p\in(0,1]\cup\{2\}$,  and interpolation upper bounds for $p\in(1,+\infty)$. We also prove that for any finite connected graph $G_n$, the best constant at $p=+\infty$ is exactly its remoteness. Moreover, for general finite connected graphs $G_n$ and $p\in(0,+\infty)$, we establish optimal bounds for $C_{G_n,p}$;
the complete graph uniquely attains the lower bound for all $p>0$;
the path graph attains the upper bound for $p>1$, but for $0<p\leq1$ there exist nonpath graphs that attain the same supremum.
Finally, we show that our approach yields sharper upper bounds than the path-based estimates obtained via the Cheeger constant.
\end{abstract}

\keywords{Poincar\'{e}-Wirtinger inequality; finite graphs; sharp constants; extremizer; Cheeger constant; $p$-Laplacian}
\maketitle

\tableofcontents

\section{Introduction}\label{s1}
The Poincar\'{e}-Wirtinger inequality is one of the fundamental inequalities in the analysis of partial differential equations and the calculus of variations.
In its classical form on Euclidean domains, it asserts that the deviation of a function from its mean value over a domain can be controlled by the Lebesgue norm of its gradient; see, for instance, \cite{P}, \cite[Chapter 7]{HLP} and \cite[Chapter 5]{E}.

Beyond the Euclidean setting, various extensions of the Poincar\'{e}-Wirtinger inequality have been developed in a broad range of geometric and metric frameworks.
In particular, such inequalities have been studied on noncompact manifolds (see \cite{LW}), Lie groups (see \cite{BPV}), and more generally, on metric measure spaces (see \cite{KZ}), where these inequalities play a significant role in the study of Sobolev spaces, potential theory, and heat kernel estimates.

In view of these developments, and given that many geometric and analytic properties admit discrete counterparts, increasing attention has been devoted to discrete analogues of the Poincar\'{e}-Wirtinger inequality.
Such inequalities have been investigated on meshes, graphs, and other discrete structures from both theoretical and numerical perspectives.
For instance,
Brenner \cite{BRE} established discrete Sobolev and Poincar\'{e}-Wirtinger inequalities for piecewise polynomial functions on two-dimensional triangulations, obtaining estimates with constants of the form $C(1+\vert \ln h\vert)$, where $h$ denotes the mesh size.
These inequalities are applicable to finite element methods that do not require functions to be continuous across element boundaries.
Le and Omnes \cite{LO} proved discrete Poincar\'{e}-Friedrichs inequalities on a two-dimensional polygonal domain covered by arbitrary, possibly nonconforming meshes,
while Carstensen and Hellwig \cite{CH} derived a discrete Poincar\'{e}-Wirtinger inequality on simplices with explicit constants. 
Related discrete inequalities have also been studied on combinatorial structures.
For instance, the optimal constants for Poincar\'{e}-Wirtinger inequalities on Hamming cubes $\{-1,1\}^n$ were obtained in \cite{ELP,ILVHV,IV}, 
and Bobkov and G\"otze \cite{BG} studied discrete isoperimetric and Poincar\'{e}-type inequalities for product probability measures on the discrete cube $\{0,1\}^n$ and on the lattice $\mathbb{Z}^n$.

Among various discrete settings, graphs provide a particularly natural framework in which the Poincar\'{e}-Wirtinger inequality reflects the interaction between geometric structure and analytic properties. Indeed,
the Poincar\'{e}-Wirtinger inequalities on graphs are closely related to geometric quantities such as volume growth and curvature, as well as to analytic and probabilistic phenomena including heat kernel estimates and random walk behavior.
In this direction, 
Delmotte \cite{DE} established that, for reversible Markov chains on weighted graphs, the conjunction of volume regularity and the Poincar\'{e}-Wirtinger inequality characterizes the parabolic Harnack inequality and two-sided Gaussian heat kernel estimates.
Later, Horn et al. \cite{HLLY} extended this criterion to locally finite graphs with non-negative curvature.
Furthermore, 
for finite graphs,
Barlow \cite[Chapter 3.3]{BAR} established the Poincar\'{e}-Wirtinger inequalities on finite weighted graphs through the framework of random walks.

In order to analyze how the constants in the Poincar\'{e}-Wirtinger inequality depend on the size of the graph  and on the exponent $p$, it is natural to seek their optimal values for fundamental graph classes.
Along this direction, Funano and Sakurai \cite[Theorem 8.5]{FS} obtained an explicit upper bound estimates for the sharp constants of high-order Poincar\'{e}-Wirtinger inequalities on weighted finite graphs in terms of the graph distance between disjoint vertex subsets and their weighted measures. 
Levi et al. \cite{LSTV} studied the growth of the sharp constants of local Poincar\'{e}-Wirtinger inequalities on infinite graphs.
For particular finite connected graphs, Hume, Mackay and Tessera \cite[Proposition 10.3]{HMT} obtained non-sharp bounds of the sharp constants of the Poincar\'{e}-Wirtinger inequalities on complete graphs, 
although an explicit expression for the optimal constant was not provided.
More recently, the sharp constant of the Poincar\'{e}-Wirtinger inequalities on complete graphs were obtained by Gonz\'alez-Riquelme and  Madrid in \cite[Theorem 1]{GRM3} for $n\geq 3$ and $p\in[1,3+\delta_n^1)\cup(3+\delta_n^2,+\infty)$, where 
\[
\delta_n^1:=\frac{1}{2n^2\ln n}+O(1/n^3),\quad 
\delta_n^2:=\frac{2}{n}+O(1/n^2).
\] 
The admissible range leaves a gap for $p\in (0,1)\cup(3,4)$, and the lack of sharp constants for path and star graphs motivates this work.

The basic graph notation used herein is as follows. A simple, connected, undirected graph $G=(V_G,E_G)$ is considered, where $V_G$ and $E_G$ denote its vertex and edge sets. Two vertices $v,w\in V_G$ are adjacent (written $v\sim w$) if $\{v,w\}\in E_G$. The graph is connected if, for any distinct $v,w$, there exists a finite path $v=a_0\sim a_1\sim\cdots\sim a_k=w$ with $k\ge 1$. It is finite if $|V_G|_\mu<\infty$, where $\vert\cdot\vert_\mu$ denotes the cardinality of the corresponding vertex set. The degree of vertex $v$ is denoted by $d_G(v)$; since $G$ is connected and has at least two vertices, $d_G(v)\ge 1$, and the maximum degree is defined as $\Delta_G:=\max_{v\in V_G} d_G(v)$.

The following graph classes are used: trees (connected acyclic graphs), complete graphs $K_n$ ($n\ge2$, every pair adjacent), star graphs $S_n$ ($n\ge3$, a tree with one central vertex adjacent to all others), and path graphs $P_n$ ($n\ge3$, vertices ordered linearly with edges only between consecutive vertices).

The object of study is the Poincar\'{e}-Wirtinger inequality on finite graphs. For a finite connected graph $G_n=(V_{G_n},E_{G_n})$ with $n$ vertices $a_1,\dots,a_n$ ($n\ge2$), and for a function $f:V_{G_n}\to\mathbb{R}$ and an exponent $0<p\le\infty$, the sharp constant $C_{G_n,p}$ is defined as the best constant such that
\[
\|f-m_f\|_p \le C_{G_n,p} {\rm Var}_{p}^{G_n}(f)
\]
holds for every $f$ on $V_{G_n}$. Here, $m_f:=\frac{\sum_{a_i\in V_{G_n}}f(a_i)}{\vert V_{G_n}\vert_\mu}$ is the arithmetic mean value of $f$,  and
\begin{equation*}
	\Vert f\Vert_p
	:=
	\begin{cases}
		\Big(\sum\limits_{a_i\in V_{G_n}}\vert f(a_i)\vert^p\Big)^{1/p} & {0<p<+\infty,}\\
		\displaystyle\max_{a_i\in V_{G_n}}\vert f(a_i)\vert & {p=+\infty,}
	\end{cases}
\end{equation*}
is the $p$-norm of $f$ ($p$-quasi norm if $0<p<1$). The notation 
${\rm Var}_{p}^{G_n}(f)$ is the $p$-variation of $f$ on $G_n$ defined by
\[
{\rm Var}_{p}^{G_n}(f)
:=
\begin{cases}
	\Big(\displaystyle \frac{1}{2}\sum_{a_i,a_j\in V_{G_n},a_i\sim a_j} \lvert f(a_i)-f(a_j)\rvert^p \Big)^{1/p} & 0<p<+\infty,\\[1.2ex]
	\displaystyle \max_{a_i,a_j\in V_{G_n},a_i\sim a_j} \lvert f(a_i)-f(a_j)\rvert & p=+\infty.
\end{cases}
\]
Accordingly, the sharp constant is defined by
\begin{align}\label{DBC}
	C_{G_n,p}
	:=\sup_{f:V_{G_n}\to\mathbb{R};{\rm Var}_{p}^{G_n}(f)>0}
	\frac{\Vert f-m_{f}\Vert_{p}}{{\rm Var}_{p}^{G_n}(f)}.
\end{align}
When $p=+\infty$, the sharp constant $C_{G_n,\infty}$ becomes
\begin{align}\label{DBCI}
	C_{G_n,+\infty}
	:=\sup_{f:V_{G_n}\to\mathbb{R};{\rm Var}_{\infty}^{G_n}(f)>0}\frac{\displaystyle\max_{a_j\in V_{G_n}}\vert f(a_j)-m_f\vert}
	{\displaystyle\max_{a_i,a_j\in V_{G_n},a_i\sim a_j}\vert f(a_i)-f(a_j)\vert}.
\end{align}

It is worth emphasizing that the sharp constant $C_{G_n,p}$ is closely related to the smallest positive eigenvalue of the $p$-Laplacian on $G_n$ for $p\in(1,+\infty)$ and the Cheeger constant of the graph $G_n$.
The $p$-Laplacian is a nonlinear generalization of the Laplacian, which corresponds to the case $p=2$, in Euclidean spaces and Riemannian manifolds.
The discrete analogue of $p$-Laplacians on graphs was proposed by Yamasaki (see \cite{Y}).
For a finite connected graph $G_n$ with $n\geq2$ vertices, the $p$-Laplacian on $G_n$ is defined by
\[
\Delta_p f(x)
:=
\sum_{\substack{y\in V_{G_n}\\y\sim x}}\omega_{xy}\vert f(y)-f(x)\vert^{p-2}(f(y)-f(x)),\quad x\in V_{G_n}, 
\]
For $p>1$, every product $|t|^{p-2}t$ is understood by its continuous extension $|0|^{p-2}0:=0$. Here $\omega_{xy}$ is a weight function mapping from $V_{G_n}\times V_{G_n}$ to $\mathbb{R}$.
Throughout the whole paper, we consider the case when the weight function $\omega_{xy}\equiv1$ for every $x\sim y$, and $\omega_{xy}=0$ if two vertices are nonadjacent.
A nonzero function $f$ is called an eigenfunction associated with the eigenvalue $\lambda$ of the $p$-Laplacian on $G_n$ if 
\[
\Delta_p f=-\lambda\vert f\vert^{p-2}f.
\]

The spectral theory of graph $p$-Laplacians has been investigated from various perspectives.
Chung, Grigor'yan and Yau \cite{CGY} developed a unified heat-kernel approach for deriving upper eigenvalue bounds for Laplacians on compact Riemannian manifolds and finite graphs.
Hua and Wang \cite{HW} studied eigenvalues and eigenfunctions of $p$-Laplacian with Dirichlet boundary conditions on symmetric graphs. 
Subsequently, Ge, Hua and Jiang \cite{GHJ} studied the limit of first eigenfunctions of $p$-Laplacian on a finite subset of a graph with Dirichlet boundary condition as $p\to1$. 

A fundamental geometric quantity associated with the spectral theory of graphs is the Cheeger constant, which was initially introduced by Cheeger \cite{C} to estimate the first non-trivial eigenvalue of the Laplace-Beltrami operator on manifolds, and was later generalized to graphs by Alon, Milman \cite{AM} and Dodziuk \cite{DO} respectively to characterize how well-connected a graph is.
The Cheeger constant $h_{\mu,\nu}(G_n)$ of the graph $G_n$ is defined by
\begin{align}\label{DOCC}
	h_{\mu,\nu}(G_n)
	:=
	\min_{S\subset V_{G_n}, 0<\vert S\vert_\mu\leq n/2}
	\frac{\vert\partial S\vert_\nu}{\vert S\vert_\mu},
\end{align}
where $\partial S$ denotes the boundary of $S$, given by
\[
\partial S:=\{\{x,y\}\in E_{G_n}:\ x\in S,\ y\in V_{G_n}\setminus S\},
\] 
and $\vert\cdot\vert_\nu$ denotes the cardinality of the corresponding edge set.

The Cheeger constant provides a geometric characterization of the first positive eigenvalue of the graph $p$-Laplacian, which in turn controls the sharp constant $C_{G_n,p}$. 
On the one hand, for $1<p<+\infty$, by the definition of $C_{G_n,p}$ from \eqref{DBC}, and the Rayleigh quotient characterization (see \cite[Theorem 1]{A1}), we have that
\begin{align}\label{REVPC}
	\lambda_{1,p}(G_n)^{-\frac{1}{p}}\leq C_{G_n,p}
	\leq2\lambda_{1,p}(G_n)^{-\frac{1}{p}},
\end{align}
where $\lambda_{1,p}(G_n)$ is the smallest positive eigenvalue for $\Delta_p$. If $p=2$, then it follows that 
\begin{align*}
	\lambda_{1,2}(G_n)^{-\frac{1}{2}}= C_{G_n,2}.
\end{align*}
On the other hand, adapting the proofs of \cite[Theorem 3]{A1} or \cite[Theorem 3.1]{KM} yields expressions for the lower and upper bounds of $\lambda_{1,p}(G_n)$ in terms of $h_{\mu,\nu}(G_n)$;
in particular, it is shown that the Cheeger constant $h_{\mu,\nu}(G_n)$ on a connected graph $G_n$ is the limit of $\lambda_{1,p}(G_n)$ as $p\to1$ (see \cite[Theorem 5.15]{CHA} and \cite[Proposition 4.1]{BH}).
Such sharp estimates of the Poincar\'{e}-Wirtinger inequality provide an alternative way to study connectivity properties of finite graphs, particularly when the nonlinear eigenvalue problem for the $p$-Laplacian is hard to analyze for $p\neq2$.

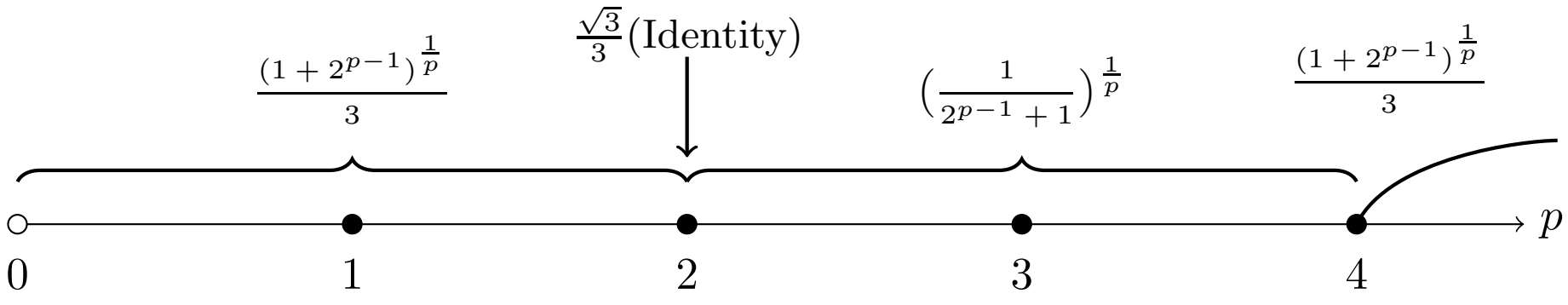
\begin{figure}[htbp]
	\centering
	\begin{tikzpicture}[scale=2.6]
		\draw[->] (0,0) -- (4.5,0) node[right] {$p$};

		\foreach \x in {0,1,2,3,4} {
			\node[anchor=base, below=0.15cm] at (\x,0) {$\x$};
		}

		\filldraw[fill=white, draw=black] (0,0) circle (0.8pt);
		\foreach \x in {1,2,3,4} {
			\filldraw[fill=black, draw=black] (\x,0) circle (0.8pt);
		}
		
		\draw[thick, decorate, decoration={brace, amplitude=5pt, raise=2pt}]
		(0,0.1) -- (2,0.1)
		node[midway, above=20pt, anchor=base, font=\tiny] 
		{$\displaystyle \frac{(1+2^{p-1})^{\frac{1}{p}}}{3}$};    
		
		\draw[thick, decorate, decoration={brace, amplitude=5pt, raise=2pt}]
		(2,0.1) -- (4,0.1)
		node[midway, above=20pt, anchor=base, font=\tiny] 
		{$\displaystyle \Big(\frac{1}{2^{p-1}+1}\Big)^{\frac{1}{p}}$};
		
		\draw[-, thick] (4,0) .. controls (4.1,0.2) and (4.5,0.25) .. (4.6,0.25)
		node[midway, above left=15pt, anchor=base, font=\tiny] 
		{$\displaystyle \frac{(1+2^{p-1})^{\frac{1}{p}}}{3}$};
		
		\draw[->, thick] (2,0.5) -- (2,0.2);
		\node[anchor=north, font=\small] at (2,0.7) {$\frac{\sqrt{3}}{3}$(Identity)};
		
	\end{tikzpicture}
	\caption{The value of $C_{K_3,p}$ for $0<p<+\infty$}
	\label{G1}
\end{figure}

\begin{figure}[htbp]
	\centering
	\begin{tikzpicture}[scale=2.6]
		\draw[->] (0,0) -- (4.5,0) node[right] {$p$};

		\foreach \x in {0,1,2,3,4} {
			\node[anchor=base, below=0.15cm] at (\x,0) {$\x$};
		}
		\node[anchor=base, below=0.15cm] at (3.5,0) {$p_n$};
		
		\filldraw[fill=white, draw=black] (0,0) circle (0.8pt);
		\foreach \x in {1,2,3,4} {
			\filldraw[fill=black, draw=black] (\x,0) circle (0.8pt);
		}
		\filldraw[fill=black, draw=black] (3.5,0) circle (0.8pt);
		
		\draw[thick, decorate, decoration={brace, amplitude=5pt, raise=2pt}]
		(0,0.1) -- (1,0.1)
		node[midway, above=20pt, anchor=base, font=\tiny] 
		{$\displaystyle \left(\frac{1+(n-1)^{p-1}}{n^p}\right)^{\frac{1}{p}}$};
		
		\draw[thick, decorate, decoration={brace, amplitude=5pt, raise=2pt}]
		(1,0.1) -- (2,0.1)
		node[midway, above=20pt, anchor=base, font=\tiny] 
		{$\displaystyle \frac{\left(\lceil \frac{n}{2}\rceil^{p-1} + \lfloor \frac{n}{2}\rfloor^{p-1}\right)^{\frac{1}{p}}}{n}$};
		
		\draw[thick, decorate, decoration={brace, amplitude=5pt, raise=2pt}]
		(2,0.1) -- (3.5,0.1)
		node[midway, above=20pt, anchor=base, font=\tiny] 
		{$\displaystyle \Big(\frac{1}{2^{p-1}+n-2}\Big)^{\frac{1}{p}}$};
		
		\draw[-, thick] (3.5,0) .. controls (3.6,0.2) and (4.5,0.25) .. (4.6,0.25)
		node[midway, above left=15pt, anchor=base, font=\tiny] 
		{$\displaystyle \left(\frac{1+(n-1)^{p-1}}{n^p}\right)^{\frac{1}{p}}$};
		
		\draw[->, thick] (1,-0.4) -- (1,-0.21);
		\node[anchor=north, font=\small] at (1,-0.41) {$\frac{2}{n}$};
		
		\draw[->, thick] (2,-0.4) -- (2,-0.21);
		\node[anchor=north, font=\small] at (2,-0.41) {$\frac{1}{\sqrt{n}}$(Identity)};
	\end{tikzpicture}
	\caption{The value of $C_{K_n,p}$ for $n\geq4$ and $0<p<+\infty$, where
		$p_n=3+\frac{2}{(4\ln2-1)n}+O\Big(\frac{\ln n}{n^2}\Big)$, $\lfloor x\rfloor:=\max\{k\in\mathbb{Z}:k\leq x\}$ and $\lceil x\rceil:=\min\{k\in\mathbb{Z}:k\geq x\}$}
	\label{G2}
\end{figure}
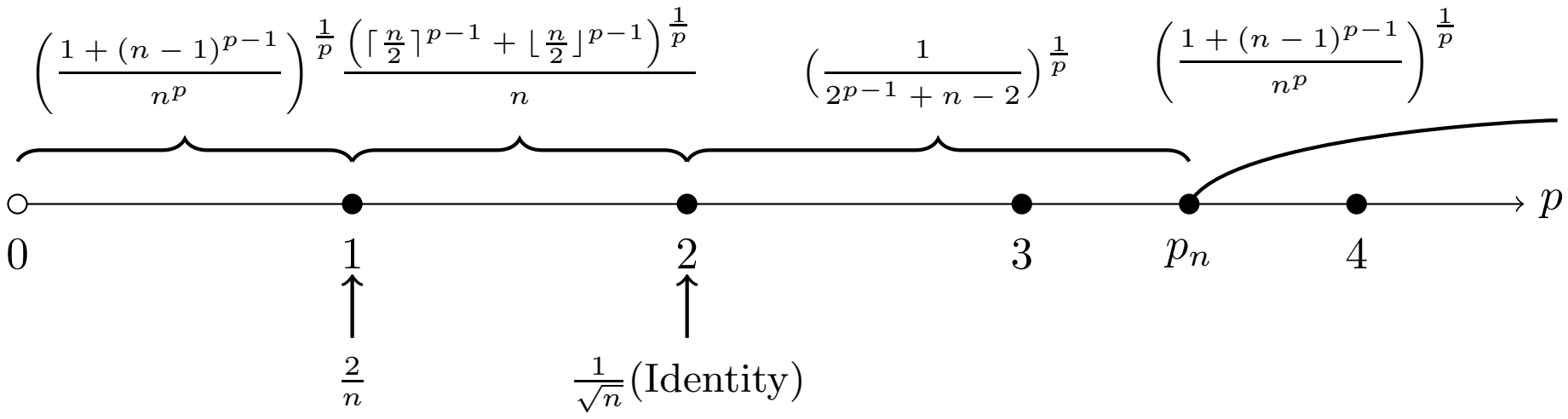

Following the above line of investigation, the paper is organized as follows.
In Section \ref{s2}, firstly,
to complete the remaining gap of the admissible range of $p$ for the sharp constants for the Poincar\'{e}-Wirtinger inequality on $K_n$ in \cite[Theorem 1]{GRM3},  we introduce the threshold value $p_n$ in Theorem \ref{KNSC} determined by 
\[
((n-1)^{p_n-1}+1)(n-2+2^{p_n-1}) = n^{p_n},
\]
and obtain the sharp constants for $p$ in $[3,p_n)$ and $[p_n,4]$ respectively.
Moreover, we establish the sharp constant for $p\in(0,1)\cup[3,4]$ and characterize the corresponding extremal functions, extending the range beyond the case $p\ge4$ treated in \cite{GRM3}. Combining the results in \cite[Theorem 1]{GRM3} and Theorem \ref{KNSC},
we summarize the complete description of $C_{K_n,p}$ for all $n\geq3$ and $p\in(0,+\infty)$ in Figures \ref{G1} and \ref{G2}.

Sections \ref{s3} and \ref{s4} are devoted to the best constants in the Poincar\'{e}-Wirtinger inequality for two other fundamental classes of finite connected graphs: star graphs and path graphs. For star graphs, Theorem \ref{SPCS} provides an explicit formula for $C_{S_n,p}$ for all $p\in(0,+\infty)$, together with a characterization of all extremizers attaining the best constant.
A key difficulty lies in determining the structure of these extremizers. In Propositions \ref{EOSN2} and \ref{EOSN3}, we prove that, for $p\in(1,2)$ and $p\in(2,+\infty)$, respectively, every extremizer takes exactly two distinct values on the noncentral vertices, both different from its value at the center.
Theorem \ref{EOSN4} then determines the numbers of noncentral vertices whose values lie above and below the central value, thereby yielding the characterization of the extremizers.

For path graphs $P_n$ with $n\geq4$ vertices, Theorem \ref{EOPN} establishes the existence of monotone extremizers for $C_{P_n,p}$ for all $p\in(0,+\infty)$, which are strictly monotone when $p>1$.
Furthermore, an explicit expression of $C_{P_n,p}$ is obtained for $p\in(0,1]\cup\{2\}$, while lower bounds and interpolation upper bounds are established for general $p\in(1,+\infty)$ in Theorem \ref{SPCP}. Reflection antisymmetry holds for $p=2$, but need not hold for general $p>1$, $p\neq2$. 
Moreover, for complete graphs $K_n$ and star graphs $S_n$ with $n\geq3$ vertices, and for path graphs $P_n$ with $n\geq4$ vertices, we prove that
\[
\lim_{p\to+\infty}C_{K_n,p}=1-\frac{1}{n},\quad \lim_{p\to+\infty}C_{S_n,p}=2-\frac{3}{n},\quad
\lim_{p\to+\infty}C_{P_n,p}=\frac{n-1}{2}.
\]
These limits are in fact special cases of a more general result established in Section \ref{s5}, where we determine the sharp constant for $p=+\infty$: for every finite connected graph $G_n$, it is given by the maximum average distance from a vertex to all vertices.

The limits obtained above are special manifestations of a more general extremal phenomenon. A natural question in the study of graph inequalities is whether certain analytic quantities can characterize the underlying geometric structures of graphs. For instance, Soria and Tradacete \cite[Theorem 3.1]{ST} showed that, for $p\in(0,1]$, the complete graph and the star graph are the extremal graphs attaining the infimum and supremum, respectively, of the $p$-quasi norm of the Hardy--Littlewood maximal operator on finite graphs (and also of the corresponding fractional quasi-norm; see \cite{XZ}). Motivated by these extremal results, in Section \ref{s5} we study the sharp constants for finite connected graphs. We prove that, for any $0<p<+\infty$ and any finite connected graph $G_n$,
\[
C_{K_n,p}\leq C_{G_n,p}\leq C_{P_n,p}.
\]
Furthermore, for every $0<p<+\infty$, the complete graph $K_n$ is precisely the graph attaining the infimum of the sharp constant.
For the upper bound, the path graph attains it when $1<p<+\infty$; however, for $0<p\leq1$, there exist graphs not isomorphic to path graphs attaining the upper bound as well (see Remark \ref{CTE}).

In Section \ref{s6}, we establish upper bounds for the sharp constants of Poincar\'{e}-Wirtinger inequalities on finite connected graphs in terms of the Cheeger constant (see \cite[Theorem 3]{A1} and Proposition \ref{p41}). 
It is concluded that the path-based Cheeger-type bounds are classical but weaker than those obtained from Theorems \ref{SPCP} and~\ref{BFG}, which give sharper upper estimates in terms of the number of vertices.

\section{Sharp Poincar\'{e}-Wirtinger inequalities on complete graphs for $p\in(0,1)\cup[3,4]$}\label{s2}

We start by proving the sharp Poincar\'{e}-Wirtinger inequalities constants on graphs. Before the proof, we introduce some basic lemmas. 

\begin{lem}[Karamata's inequality \cite{Ka}]
Let $f$ be a convex function on some interval $I$ of the real line. If $a_1,\ldots,a_n$, and $b_1,\ldots,b_n$ are numbers in $I$ such that the sequence $\{a_j\}_{j=1}^n$ majorizes $\{b_j\}_{j=1}^n$, which means, $a_1\geq\cdots\geq a_n$, $b_1\geq\cdots\geq b_n$, $\sum_{j=1}^{n}a_j=\sum_{j=1}^{n}b_j$ and $\sum_{j=1}^{k}a_j\geq\sum_{j=1}^{k}b_j$ for every $1\leq k\leq n-1$, then we have
\[
\sum_{j=1}^{n}f(a_j)\geq\sum_{j=1}^{n}f(b_j).
\]
\end{lem}

\begin{lem}\label{REU}
For any integer $n\geq4$, and real number $p\in[3,4]$, the equation of $p$
\begin{align}\label{EQKN01}
((n-1)^{p-1}+1)(n-2+2^{p-1})=n^p
\end{align}
admits unique solution $p_n\in\big(3+\frac{1}{n},3+\frac{3}{n}\big)$ and
\[
p_n
=3+\frac{2}{(4\ln2-1)n}+O\Big(\frac{\ln n}{n^2}\Big)
\] as $n\to \infty$.
\end{lem}

\begin{proof}
Define $H(p):=\ln((n-1)^{p-1}+1)+\ln(n-2+2^{p-1})-p\ln n, \ p\in[3,4]$. Then we have that 
\begin{align*}
H'(p)
&=    \frac{(n-1)^{p-1}\ln(n-1)}{(n-1)^{p-1}+1}+\frac{2^{p-1}\ln2}{n-2+2^{p-1}}-\ln n\\
&\geq \frac{\ln(n-1)}{1+(n-1)^{-2}}+\frac{4\ln2}{n+2}-\ln n\\
&=    \ln\Big(1-\frac{1}{n}\Big)+\frac{4\ln2}{n+2}-\frac{\ln(n-1)}{1+(n-1)^2}.
\end{align*}
Since $\ln(1-\frac{1}{n})>-\frac{1}{n-1}$ for $n\geq4$, it follows that
\begin{align*}
H'(p)
&> -\frac{1}{n-1}+\frac{4\ln2}{n+2}-\frac{\ln(n-1)}{(n-1)^2}\\
&=    \frac{(4\ln2-1)(n-1)^2-3(n-1)-(n+2)\ln(n-1)}{(n-1)^2 (n+2)}.
\end{align*}

Define $s(t):=(4\ln2-1)t^2-3t-(t+3)\ln t$, $t\in [3,+\infty)$. We have that 
\[
s'(t)=2(4\ln2-1)t-4-\frac{3}{t}-\ln t,
\]
and 
\[
s''(t)
=2(4\ln2-1)+\frac{3}{t^2}-\frac{1}{t}
\geq s''(6)
=8\ln2-\frac{25}{12}
\approx3.4618>0.
\]
Thus,  
\[
s'(t)
\geq s'(3)
=24\ln2-\ln3-11
\approx4.5369>0,
\]
and consequently, 
\begin{align*}
(4\ln2-1)(n-1)^2-3(n-1)-(n+2)\ln(n-1)
\geq s(3)
\approx0.3616>0,
\end{align*}
implying that $H'(p)>0$ on $p\in[3,4]$ for all $n\geq4$.


We claim that $H(3+\frac{1}{n})<0$ and $H(3+\frac{3}{n})>0$ for all $n\geq4$. Indeed, $H(3+\frac{1}{n})<0$ and $H(3+\frac{3}{n})>0$ are equivalent to 
\begin{align}\label{H31N}
\mathcal{H}_1(n):=\Big(\Big(1-\frac{1}{n}\Big)^{2+\frac{1}{n}}+\Big(\frac{1}{n}\Big)^{2+\frac{1}{n}}\Big)
\Big(1+\frac{2}{n}(2^{\frac{1}{n}+1}-1)\Big)<1
\end{align}
and
\begin{align}\label{H33N}
\mathcal{H}_2(n):=\Big(\Big(1-\frac{1}{n}\Big)^{2+\frac{3}{n}}+\Big(\frac{1}{n}\Big)^{2+\frac{3}{n}}\Big)
\Big(1+\frac{2}{n}(2^{\frac{3}{n}+1}-1)\Big)>1
\end{align}
respectively for any $n\geq4$.

We first prove \eqref{H31N}. Since $(1-\frac{1}{n})^{\frac{1}{n}}<1-\frac{1}{n^2}$ from Bernoulli's inequality and $2^{\frac{1}{n}}<1+\frac{\ln2}{n}+\frac{\ln2}{2n^2}$ for $n\geq4$, it then follows that
\begin{align*}
\mathcal{H}_1(n)
&<
\Big(\Big(1-\frac{1}{n}\Big)^{2}\Big(1-\frac{1}{n^2}\Big)+\frac{1}{n^2}\Big)
\Big(1+\frac{2}{n}\Big(1+\frac{2\ln2}{n}+\frac{\ln2}{n^2}\Big)\Big).
\end{align*}
Let
\[G_1(t):=
((1-t)^2(1-t^2)+t^2)(1+2t(1+2t\ln2 +t^2\ln2)),\ t\in (0,1/4].\]
Then we have
\begin{align*}
G_1(t)
&=1
+t^2\Big(
(4\ln2-3)
+2(2-3\ln2)t
+3t^2  
+2(5\ln2-1)t^3 
-2t^5\ln2\Big) \\
&< 1
+t^2\Big(
(4\ln2-3)
+2(2-3\ln2)t
+\frac{5}{2}(1+\ln2)t^2\Big).
\end{align*}

Let 
\[G_{1,1}(t):=(4\ln2-3)
+2(2-3\ln2)t
+\frac{5}{2}(1+\ln2)t^2,\ t\in [0,1/4].\]
It is readily verified that
\[G_{1,1}(t)\leq \max\{G_{1,1}(0),\  G_{1,1}(1/4)\}=\max\Big\{4\ln2-3,\ \frac{85\ln2-59}{32}\Big\}<0\] for $t\in [0,1/4]$, which implies  \eqref{H31N} holds for any $n\geq4$.

We next prove \eqref{H33N}. 
For $n=4,5,6,7$, calculations yield that 
\[
\mathcal{H}_2(4)
=
\frac{(3^{\frac{11}{4}}+1)(2+2^{\frac{11}{4}})}{4^{\frac{15}{4}}}
\approx 1.0373,
\]
\[
\mathcal{H}_2(5)
=
\frac{(4^{\frac{13}{5}}+1)(3+2^{\frac{13}{5}})}{5^{\frac{18}{5}}}
\approx 1.0423,
\] 
\[
\mathcal{H}_2(6)
=
\frac{(5^{\frac{5}{2}}+1)(4+2^{\frac{5}{2}})}{6^{\frac{7}{2}}}
\approx 1.0386, 
\]
\[
\mathcal{H}_2(7)
=\frac{(6^{\frac{17}{7}}+1)(5+2^{\frac{17}{7}})}{7^{\frac{24}{7}}}
\approx 1.0334.
\]

For $n\geq8$, since $2+\frac{3}{n}<3$, we have 
\begin{align*}
\mathcal{H}_2(n)>
\Big(\Big(1-\frac{1}{n}\Big)^{2+\frac{3}{n}}+\Big(\frac{1}{n}\Big)^{3}\Big)
\Big(1+\frac{2}{n}(2^{\frac{3}{n}+1}-1)\Big).
\end{align*}
Let 
\[G_2(t):=((1-t)^2(1-t)^{3t}+t^3)(1+2t(2^{1+3t}-1)),\ t\in (0, 1/8].\]
Note that 
\[
(1-t)^{3t}
=e^{3t\ln(1-t)}>e^{\frac{-3t^2}{1-t}}>1-\frac{3t^2}{1-t},
\] and $2^{1+3t}-1>1+(6\ln2)t$ for $t\in (0, 1/8]$. It follows that
\begin{align*}
G_2(t)&>
\Big((1-t)^2\Big(1-\frac{3t^2}{1-t}\Big)+t^3\Big)(1+2t(1+(6\ln2)t))\\
&=
1+t^2\Big((12\ln2-6)-(24\ln2)t-(24\ln2-8)t^2+(48\ln2)t^3\Big)\\
&>
1+t^2\Big((12\ln2-6)-(24\ln2)t-(24\ln2-8)t^2\Big).
\end{align*}

Let 
\[G_{2,1}(t):=(12\ln2-6)-(24\ln2)t-(24\ln2-8)t^2, \ t\in [0,1/8].\]
A direct computation gives
\[
\min_{t\in[0,1/8]}G_{2,1}(t)=G_{2,1}(1/8)
= \frac{69\ln2-47}{8} \approx 0.1034>0,
\]
which implies that \eqref{H33N} holds for $n\ge 8$. Hence, together with the cases $4\le n\le 7$ already established, \eqref{H33N} holds for all $n\ge4$.

The strict monotonicity of $H(p)$ on $[3,4]$, together with \eqref{H31N} and \eqref{H33N}, implies that \eqref{EQKN01} has a unique solution $p=p_n$, which satisfies
\[
p_n\in\Big(3+\frac{1}{n},\,3+\frac{3}{n}\Big)
\quad\text{and}\quad
\lim_{n\to\infty} p_n = 3.
\]

We now turn to the asymptotics  of $p_n$ as $n\to \infty$. To this end, for $x\in[3,4]$ and $n\geq4$, define
\[
F(x):=((n-1)^{x-1}+1)(n-2+2^{x-1})-n^x,
\]
and set $p_n=3+\varepsilon_n$ with $1/n<\varepsilon_n<3/n$ for $n\geq 4$. An application of the sandwich theorem yields
\[
\lim_{n\to\infty} \varepsilon_n \ln n = 0.
\]
To justify the expansion with a uniform remainder, write $F(x)=n^x(F(x)/n^x)$. For $3\leq x\leq3+3/n$, direct differentiation gives
\[
\frac{F(x)}{n^x}=O(n^{-2}),\qquad
\left(\frac{F(x)}{n^x}\right)'=O(n^{-1}),\qquad
\left(\frac{F(x)}{n^x}\right)''=O(n^{-1}).
\]
Hence $F''(x)=O(n^2\ln n)$ uniformly on this interval. Since $\varepsilon_n=O(1/n)$, Taylor's formula yields
\begin{align*}
F(p_n)
&=F(3)+F'(3)\varepsilon_n+O(\ln n)\\
&=4-2n+\Big((2-3n)\ln n+(n+2)(n-1)^2\ln\Big(1-\frac1n\Big)\\
&\qquad+4(n^2-2n+2)\ln2\Big)\varepsilon_n+O(\ln n),
\end{align*}
as $n\to\infty$.

The condition $F(p_n)=0$ gives
\begin{align*}
	&\Big((2-3n)\ln n+(n+2)(n-1)^2\ln\Big(1-\frac{1}{n}\Big)+4(n^2-2n+2)\ln2\Big)\varepsilon_n\\
	=& 2n-4 +O(\ln n),
\end{align*} as $n\to \infty$. By Taylor's formula,  $\ln(1-\frac{1}{n})
=-\frac{1}{n}-\frac{1}{2n^2}+O(\frac{1}{n^3})$ as $n\to \infty$. It then follows that
\begin{align*}
	&(2-3n)\ln n+(n+2)(n-1)^2\ln\Big(1-\frac{1}{n}\Big)+4(n^2-2n+2)\ln2\\
	=&(2-3n)\ln n +(n^3-3n+2)\Big(-\frac{1}{n}-\frac{1}{2n^2}+O\Big(\frac{1}{n^3}\Big)\Big) +4(n^2-2n+2)\ln2 \\
	=& (4\ln2-1)n^2+O(n\ln n),
\end{align*} as $n\to \infty$.
Hence using the expansion $\frac{1}{1+x}=1+O(x)$ as $x\to 0$, we have
\begin{align*}
	\varepsilon_n
	&=
	\frac{2n-4}{(4\ln2-1)n^2+O(n\ln n)}+O\Big(\frac{\ln n}{n^2}\Big)\\
	&=
	\frac{2}{(4\ln2-1)n}\cdot\frac{1-\frac{2}{n}}{1+O(\frac{\ln n}{n})}
	+O\Big(\frac{\ln n}{n^2}\Big)\\
	&=
	\frac{2}{(4\ln2-1)n}\Big(1-\frac{2}{n}\Big)\Big(1+O\Big(\frac{\ln n}{n}\Big)\Big)
	+O\Big(\frac{\ln n}{n^2}\Big)\\
	&=
	\frac{2}{(4\ln2-1)n}+O\Big(\frac{\ln n}{n^2}\Big),
\end{align*}
implying that
\[
p_n=3+\frac{2}{(4\ln2-1)n}+O\Big(\frac{\ln n}{n^2}\Big),
\]
as $n\to \infty$.
This completes the proof.
\end{proof}

\begin{lem}\label{N3CASE}
Let $p\in[3,4]$, then $3^p\geq(2^{p-1}+1)^2$.
\end{lem}
\begin{proof}
Let $f(x)=3^x-(2^{x-1}+1)^2$, $x\in[3,4]$. Then we have that 
\[
f'(x)=3^x \ln3  -2(4^{x-1}+2^{x-1})\ln2
\]
and
\begin{align*}
f''(x)
=3^x \Big(\ln^2 3-\Big(\Big(\frac{4}{3}\Big)^x+\Big(\frac{2}{3}\Big)^x\Big)\ln^2 2\Big).
\end{align*}
From $3^5=243<2^8=256$, we can compute that $5\ln3<8\ln2$ and $\Big(\frac{\ln3}{\ln2}\Big)^2<\frac{64}{25}$. Thus, it follows that 
\[
\Big(\frac{4}{3}\Big)^x+\Big(\frac{2}{3}\Big)^x
>\Big(\frac{4}{3}\Big)^3+\Big(\frac{2}{3}\Big)^4
=\frac{208}{81}>\frac{64}{25}>\Big(\frac{\ln3}{\ln2}\Big)^2,
\]
implying $f''(x)<0$ and $f$ is concave on the interval $[3,4]$. Thus we obtain that 
\[
\min\limits_{x\in[3,4]}f(x)=\min\{f(3),\ f(4)\}=0.
\]
Therefore,  $3^p\geq(2^{p-1}+1)^2$ holds for every $p\in[3,4]$. This completes the proof.
\end{proof}

\begin{lem}\label{NOMAX}
For any integer $n\geq3$, and real number $p\in(3,4)$, the function $\psi_n:[0,1]\to\mathbb{R}$ is defined as
\[
\psi_n(x):=\frac{((n-1)-(n-2)x)^{p-1}+(1+(n-2)x)^{p-1}}{(n-2)((1-x)^{p-1}+x^{p-1})+2},
\quad x\in[0,1], 
\] 
then for $n\geq4$,
\begin{align*}
\max_{x\in[0,1]}\psi_n(x)
=
\begin{cases}
\frac{n^{p-1}}{2^{p-1}+n-2} & {p\in(3,p_n],}\\[0.2cm]
\frac{(n-1)^{p-1}+1}{n} & {p\in(p_n,4),}
\end{cases}
\end{align*}
where $p_n$ is the solution to the equation 
\[
f(p)=((n-1)^{p-1}+1)(n-2+2^{p-1})-n^p=0,\ p\in (3,4).
\]
If $n=3$, then 
\begin{align}\label{eq6333}
\max_{x\in[0,1]}\psi_3(x)=\frac{3^{p-1}}{2^{p-1}+1}.
\end{align}
\end{lem}

\begin{proof}
For $n\geq3$, one can check that $\psi_n(x)=\psi_n(1-x)$, $\psi_n\in C^2([0,1])\cap C^\infty((0,1))$ and 
\[\psi_n(0)=\frac{(n-1)^{p-1}+1}{n}, \  \psi_n(1/2)=\frac{n^{p-1}}{2^{p-1}+n-2}.\] 
We consider two cases for $\psi_n$ on the interval $[0,1/2]$, depending on whether $\psi_n$ is monotonic there.

\textbf{Case 1.} If $\psi_n$ is monotonically decreasing or increasing on $[0,1/2]$, then
\[
\max_{x\in[0,1]}\psi_n(x)=\max\Big\{\psi_n(0),\ \psi_n(1/2)\Big\}.
\]

\textbf{Case 2.} Assume that there exists $x_0\in(0,1/2)$ satisfying $\psi_n'(x_0)=0$. We claim that $\psi_n''(x_0)>0$. Once this is established, $x_0$ is a local minimum, and consequently the maximum of $\psi_n$ is attained at $x=0$ or $x=1/2$.

We now prove the claim.  For $x\in[0,1/2]$, set
\[
g_n(x):=((n-1)-(n-2)x)^{p-1}+(1+(n-2)x)^{p-1},
\]
and 
\[
h_n(x):=(n-2)((1-x)^{p-1}+x^{p-1})+2.
\]
Then $\psi_n(x)=\frac{g_n(x)}{h_n(x)}$ and we have that 
\begin{align}\label{eqx21}
\frac{\psi'_n(x)}{\psi_n(x)}=\frac{g'_n(x)}{g_n(x)}-\frac{h'_n(x)}{h_n(x)},
\end{align}
and
\begin{align*}
\Big(\frac{\psi'_n(x)}{\psi_n(x)}\Big)'
=\frac{\psi''_n(x)}{\psi_n(x)}-\Big(\frac{\psi'_n(x)}{\psi_n(x)}\Big)^2
=\frac{g''_n(x)}{g_n(x)}-\Big(\frac{g'_n(x)}{g_n(x)}\Big)^2
-\Big(\frac{h''_n(x)}{h_n(x)}-\Big(\frac{h'_n(x)}{h_n(x)}\Big)^2\Big).
\end{align*}

The conditions $\psi'_n(x_0)=0$ and $\psi_n(x_0)>0$ yield that
\[
\frac{\psi''_n(x_0)}{\psi_n(x_0)}
=\frac{g''_n(x_0)}{g_n(x_0)}-\frac{h''_n(x_0)}{h_n(x_0)},
\]
which implies that $\psi''_n(x_0)>0$ is equivalent to \begin{align}\label{eqx22}
\frac{g''_n(x_0)}{g_n(x_0)}-\frac{h''_n(x_0)}{h_n(x_0)}>0.
\end{align}
A straightforward calculation yields
\[
g'_n(x_0)=(p-1)(n-2)\Big((1+(n-2)x_0)^{p-2}-(n-1-(n-2)x_0)^{p-2}\Big),
\]
\[
g''_n(x_0)=(p-1)(p-2)(n-2)^2\Big((1+(n-2)x_0)^{p-3}+(n-1-(n-2)x_0)^{p-3}\Big),
\]
and
\[
h'_n(x_0)=(p-1)(n-2)\Big(x_0^{p-2}-(1-x_0)^{p-2}\Big),
\]
\[
h''_n(x_0)=(p-1)(p-2)(n-2)\Big(x_0^{p-3}+(1-x_0)^{p-3}\Big).
\]

Noting that $x_0\in(0,1/2)$, and using $\psi'_n(x_0)=0$ together with \eqref{eqx21}, we have
\begin{align*}
&(n-2)\Big((1-x_0)^{p-1}+x_0^{p-1}\Big)+2\\
=&
\frac{\Big(((n-1)-(n-2)x_0)^{p-1}+(1+(n-2)x_0)^{p-1}\Big)(x_0^{p-2}-(1-x_0)^{p-2})}
{(1+(n-2)x_0)^{p-2}-(n-1-(n-2)x_0)^{p-2}}.
\end{align*}
Therefore, \eqref{eqx22} is equivalent to
\begin{align}\label{PPEQ1}
&(n-2)\Big((1+(n-2)x_0)^{p-3}+(n-1-(n-2)x_0)^{p-3}\Big)\\ \nonumber
>&
\frac{(1-x_0)^{p-3}+x_0^{p-3}}{(1-x_0)^{p-2}-x_0^{p-2}}
\Big((n-1-(n-2)x_0)^{p-2}-(1+(n-2)x_0)^{p-2}\Big).
\end{align} 
Define $a:=1+(n-2)x_0$, $b:=n-1-(n-2)x_0$, then $x_0=\frac{a-1}{n-2}$, $1-x_0=\frac{b-1}{n-2}$ and \eqref{PPEQ1} can be written as
\begin{align*}
(n-2)(a^{p-3}+b^{p-3})
&>
\frac{\Big(\frac{a-1}{n-2}\Big)^{p-3}+\Big(\frac{b-1}{n-2}\Big)^{p-3}}
{\Big(\frac{b-1}{n-2}\Big)^{p-2}-\Big(\frac{a-1}{n-2}\Big)^{p-2}}
(b^{p-2}-a^{p-2})\\
&=
(n-2)\frac{b^{p-2}-a^{p-2}}{(b-1)^{p-2}-(a-1)^{p-2}}((a-1)^{p-3}+(b-1)^{p-3}),
\end{align*}
which is equivalent to 
\begin{align}\label{PPEQ2}
\frac{a^{p-3}+b^{p-3}}{b^{p-2}-a^{p-2}}
>
\frac{(a-1)^{p-3}+(b-1)^{p-3}}{(b-1)^{p-2}-(a-1)^{p-2}}.
\end{align}

Now define $\alpha:=p-3\in(0,1)$, $t:=a-1$, $l:=b-a=(n-2)(1-2x_0)\in(0,n-2)$, and define the function $\mathcal{F}_{l}:(0,+\infty)\to\mathbb{R}$ as
\[
\mathcal{F}_{l}(t):=\frac{t^\alpha+(t+l)^\alpha}{(t+l)^{\alpha+1}-t^{\alpha+1}}.
\]
To prove that \eqref{PPEQ2} holds, it suffices to prove that $\mathcal{F}_{l}$ is monotonically increasing. Indeed, letting $u:=1+{l}/{t}={1}/{x_0}-1\in(1,+\infty)$, we have that
\begin{align*}
\mathcal{F}'_{l}(t)
&=
\frac{t^{2\alpha}}{\big((t+l)^{\alpha+1}-t^{\alpha+1}\big)^2}
\Bigg(
\alpha\Big(1+\Big(1+\frac{l}{t}\Big)^{\alpha-1}\Big)
\Big(\Big(1+\frac{l}{t}\Big)^{\alpha+1}-1\Big)
\\
&\qquad\qquad
-(\alpha+1)
\Big(\Big(1+\frac{l}{t}\Big)^\alpha-1\Big)
\Big(\Big(1+\frac{l}{t}\Big)^\alpha+1\Big)
\Bigg)
\\
&=
\frac{t^{2\alpha}}{\big((t+l)^{\alpha+1}-t^{\alpha+1}\big)^2}
\Big(
1-u^{2\alpha}+\alpha(u^{\alpha+1}-u^{\alpha-1})
\Big).
\end{align*}

Define $\mathcal{P}(u):=1-u^{2\alpha}+\alpha(u^{\alpha+1}-u^{\alpha-1}),\ u\in(1,+\infty)$. Then we have that
\[
\mathcal{P}'(u)
=\alpha u^{\alpha-2}\Big((\alpha+1)u^2-2u^{\alpha+1}+1-\alpha\Big)
=:\alpha u^{\alpha-2}\mathcal{S}(u),
\]
and
\[
\mathcal{S}'(u)=2(\alpha+1)(u-u^\alpha)>0
\]
because of $\alpha\in(0,1)$ and $u\in(1,+\infty)$. Hence we have that
\[
\mathcal{S}(u)>\lim_{u\to 1^+}\mathcal{S}(u)=0,
\]
which leads to 
\[
\mathcal{P}(u)>\lim_{u\to 1^+}\mathcal{P}(u)=0.
\]
Since $\mathcal{F}'_l(t)>0$ on $(0,+\infty)$, $\mathcal{F}_l$ is strictly increasing there.

Combining Cases 1 and 2, for $n\geq3$, we obtain that 
\[
\max_{x\in[0,1]}\psi_n(x)=\max\Big\{\psi_n(0),\ \psi_n(1/2)\Big\}.
\]

If $n=3$, then from Lemma \ref{N3CASE}, we obtain that
\[
\max_{x\in[0,1]}\psi_n(x)
=
\max\Big\{\frac{2^{p-1}+1}{3},\ \frac{3^{p-1}}{2^{p-1}+1}\Big\}
=
\frac{3^{p-1}}{2^{p-1}+1}.
\]
If $n\geq4$, then by the proof of  Lemma \ref{REU}, we know that $\psi_n(1/2)\geq\psi_n(0)$ whenever $p\in (3,p_n]$; $\psi_n(1/2)<\psi_n(0)$ whenever $p\in(p_n,4)$. Therefore,
\begin{align*}
\max_{x\in[0,1]}\psi_n(x)
=
\begin{cases}
\frac{n^{p-1}}{2^{p-1}+n-2} & {p\in (3,p_n],}\\[0.2cm]
 \frac{(n-1)^{p-1}+1}{n} & {p\in(p_n,4).}
\end{cases}
\end{align*}
This completes the proof.
\end{proof}

As in \cite[Lemma 2]{GRM3}, we may prove the existence of extremizers of Poincar\'e-Wirtinger inequalities on general finite connected graphs.

\begin{lem}\label{EX}
Let $G_n$ be the finite connected graph with $n\geq3$ vertices, and $0<p<+\infty$. Then there exists $f:V_{G_n}\to\mathbb{R}^+$ for which ${\rm Var}_{p}^{G_n}(f)>0$ and
\begin{align}\label{CS1}
C_{G_n,p}
=\frac{\Vert f-m_f\Vert_p}{{\rm Var}_{p}^{G_n}(f)}.
\end{align}
\end{lem}
\begin{proof}
Let $A=\big\{x=(x_1,x_2,\ldots,x_n):x_i\in[0,1],\ \displaystyle\min_{i=1,2,\ldots,n}x_i=0,\ \displaystyle\max_{i=1,2,\ldots,n}x_i=1\big\}$. We claim that there exists a function $f_x:V_{G_n}\to\mathbb{R}$ with $f_x(i)=x_i$ for any $i\in V_{G_n}$ satisfying
\begin{align*}
C_{G_n,p}=\frac{\Vert f_x-m_{f_x}\Vert_p}{{\rm Var}_{p}^{G_n}(f_x)}.
\end{align*}
In fact, for every nonconstant $g: V_{G_n}\to\mathbb{R}$ we have that the quantity
\begin{align*}
\frac{\Vert g-m_g\Vert_p}{{\rm Var}_{p}^{G_n}(g)}
\end{align*}
remains unchanged by normalization $f\mapsto\frac{f}{\displaystyle\max_{i=1,2,\ldots,n}f(i)}$
and translation $f\mapsto f+c$, $c\in\mathbb{R}$. Every nonconstant function can therefore be translated and normalized to a point of $A$. Since $G_n$ is connected, the denominator is strictly positive on $A$. The quotient is continuous on the compact set $A$, so there exists $x_0\in A$ such that
\begin{align*}
\frac{\Vert f_{x_0}-m_{f_{x_0}}\Vert_p}{{\rm Var}_{p}^{G_n}(f_{x_0})}
=\sup_{x\in A}\frac{\Vert f_x-m_{f_x}\Vert_p}{{\rm Var}_{p}^{G_n}(f_x)}
=C_{G_n,p}.
\end{align*}
Taking $f=f_{x_0}+1$ gives strictly positive values and leaves the quotient unchanged. This completes the proof.
\end{proof}

We now present the explicit expression for $C_{K_n,p}$ for $p\in(0,1)\cup(3,4)$, thus filling the gap left open in \cite[Theorem 1]{GRM3}.
\begin{thm}\label{KNSC}
	Let $K_n$ be the complete graph with $n\geq3$ vertices $a_1$, $a_2,\ldots, a_n$. Let $f:V_{K_n}\to\mathbb{R}$. The following statements hold:
	
	(i) If $0<p<1$ and $n\geq3$, then
	\begin{align}\label{SPCKN}
		\sum_{j=1}^{n}\vert f(a_j)-m_f\vert^p
		\leq
		\frac{1+(n-1)^{p-1}}{n^p}
		{\rm Var}_{p}^{K_n}(f)^p.
	\end{align}
	The inequality is sharp and the equality attained for the Dirac delta functions $\delta_{a_i}$ ($i=1,\dots,n$), defined by
	\[
	\delta_{a_i}(a_j)=
	\begin{cases}
		1 & j=i,\\
		0 & \text{otherwise}.
	\end{cases}
	\]

	(ii) If $3<p<4$ and $n=3$, then 
	\[
	\sum_{j=1}^{3}\vert f(a_j)-m_f\vert^p
	\leq
	\frac{1}{2^{p-1}+1}
	{\rm Var}_{p}^{K_3}(f)^p.
	\]
	The inequality is sharp and the equality is attained for any function $f:V_{K_3}\to \mathbb{R}$ assuming values
	\[
	\big(f(a_1),f(a_2),f(a_{3})\big)=\big(a+b,a,a-b\big)
	\]
	for $a,b\in\mathbb{R}$ with $b\neq0$.
	
	(iii) If $3<p\leq p_n$ and $n\geq4$, then
	\[
	\sum_{j=1}^{n}\vert f(a_j)-m_f\vert^p
	\leq
	\frac{1}{2^{p-1}+n-2}
	{\rm Var}_{p}^{K_n}(f)^p,
	\]
	where $p_n$ is the root of the equation 
	\[
	((n-1)^{p-1}+1)(n-2+2^{p-1})=n^p
	\]
	in the interval $(3,4)$. The inequality is sharp and the equality is attained for any function $f:V_{K_n}\to \mathbb{R}$ assuming values
	\[
	\Big(f(a_1),f(a_2),\ldots,f(a_{n-1}),f(a_n)\Big)=\Big(a+b,a,\ldots,a,a-b\Big)
	\]
	for $a,b\in\mathbb{R}$ with $b\neq0$.
	
	(iv) If $p_n<p\leq 4$ and $n\geq4$, then
	\[
	\sum_{j=1}^{n}\vert f(a_j)-m_f\vert^p
	\leq
	\frac{1+(n-1)^{p-1}}{n^p}
	{\rm Var}_{p}^{K_n}(f)^p.
	\]
	The inequality is sharp and the equality is attained for Dirac delta functions.
\end{thm}

\begin{proof} For $p\in(0,1)\cup[3,4]$, by Lemma \ref{EX}, there exists an extremizer function $f$ on $K_n$ such that equality 
	\[
	\sum_{i=1}^{n}\vert f(a_i)-m_f\vert^p=C_{K_n,p}^p\sum_{1\leq i\leq j\leq n}\vert f(a_i)-f(a_j)\vert^p
	\]
	holds. Without loss of generality, assume the vertices are ordered so that
	\[
	f(a_1)\ge f(a_2)\ge\cdots \ge f(a_{n-1})\ge f(a_n).
	\]
	Then there exists an integer $k\in\{1,\dots,n-1\}$ such that
	\[
	f(a_1)\ge\cdots\ge f(a_k)\ge m_f\ge f(a_{k+1})\ge\cdots\ge f(a_n).
	\]

	(i) Assume $0<p<1$ and $n\geq3$. We first prove that
\[
C_{K_n,p}^p=\frac{1+(n-1)^{p-1}}{n^p}.
\]
By translation and dilation invariance, normalize the extremizer so that
\[
\min_{1\leq j\leq n}f(a_j)=0,\qquad \max_{1\leq j\leq n}f(a_j)=1.
\]
Set
\[
I_+:=\{j:f(a_j)=1\},\quad I_-:=\{j:f(a_j)=0\},\quad
I_0:=\{j:0<f(a_j)<1\},
\]
and let $r:=|I_+|$, $s:=|I_-|$ and $k:=|I_0|$, so that $r,s\geq1$ and $r+s+k=n$. For sufficiently small $|\varepsilon|$, define
\[
f_\varepsilon(a_j):=
\begin{cases}
1+\frac{\varepsilon}{r} & {j\in I_+,}\\[0.1cm]
-\frac{\varepsilon}{s} & {j\in I_-,}\\[0.1cm]
f(a_j) & {j\in I_0.}
\end{cases}
\]
Then $m_{f_\varepsilon}=m_f$. Define
\[
R(f_\varepsilon)(\varepsilon):=\Vert f_\varepsilon-m_f\Vert_p^p
-C_{K_n,p}^p{\rm Var}_p^{K_n}(f_\varepsilon)^p.
\]
The extremality of $f$ gives $R(f_\varepsilon)(0)=0$ and $R(f_\varepsilon)(\varepsilon)\leq0$. Under this perturbation, all zero differences inside $I_+$ or $I_-$ remain identically zero, differences inside $I_0$ remain unchanged, and every varying difference is nonzero near $0$. The values in $I_0$ and the mean are fixed. Thus $R$ is differentiable near $0$, and Fermat's theorem gives $R'(f_\varepsilon)(0)=0$. Consequently,
\[
C_{K_n,p}^p=
\frac{(1-m_f)^{p-1}+m_f^{p-1}}
{r+s+\displaystyle\sum_{j\in I_0}\big((1-f(a_j))^{p-1}+f(a_j)^{p-1}\big)}.
\]
Since $1/n\leq m_f\leq1-1/n$, convexity of $\Phi(t):=(1-t)^{p-1}+t^{p-1}$ gives
\[
(1-m_f)^{p-1}+m_f^{p-1}
\leq n^{1-p}\big(1+(n-1)^{p-1}\big).
\]
For every $j\in I_0$, both $(1-f(a_j))^{p-1}$ and $f(a_j)^{p-1}$ exceed $1$. Hence the denominator is at least $r+s+2k=n+k\geq n$, and
\begin{align}\label{eq6688}
C_{K_n,p}^p\leq\frac{1+(n-1)^{p-1}}{n^p}.
\end{align}
To establish the optimality of the constant, we consider the test function defined by
\[
\widetilde{f}(a_1)=1,\quad \widetilde{f}(a_2)=\cdots=\widetilde{f}(a_n)=0.
\]
This gives the lower estimate
\[
C_{K_n,p}^p\geq\frac{\Vert\widetilde{f}-m_{\widetilde{f}}\Vert_p^p}{{\rm Var}_p^{K_n}(\widetilde{f})^p}
=\frac{1+(n-1)^{p-1}}{n^p}.
\]
Together with \eqref{eq6688}, this proves the claimed formula.

We now consider $p\in[3,4]$. Choose $\varepsilon>0$ sufficiently small such that
\[
f(a_1)+\varepsilon>f(a_2),\qquad f(a_n)-\varepsilon<f(a_{n-1}),
	\]
	and define the auxiliary function $f_{\varepsilon,K_n}^{(s,t)}:V_{K_n}\to\mathbb{R}$ by
	\begin{equation*}
		f_{\varepsilon,K_n}^{(s,t)}(a_j):=
		\begin{cases}
			f(a_s)+\varepsilon & {j=s,}\\
			f(a_t)-\varepsilon & {j=t,}\\
			f(a_j) & {j\neq s,t.}
		\end{cases}
	\end{equation*}
	Then $m_{f}=m_{f_{\varepsilon,K_n}^{(s,t)}}$.

	For $\varepsilon>0$ given above, define the remainder function
	\begin{align*}
		R(f_{\varepsilon,K_n}^{(1,n)})(\varepsilon)
		&:=
		\sum_{j=1}^{k}\Big(f_{\varepsilon,K_n}^{(1,n)}(a_j)-m_{f_{\varepsilon,K_n}^{(1,n)}}\Big)^p
		+\sum_{j=k+1}^{n}\Big(m_{f_{\varepsilon,K_n}^{(1,n)}}-f_{\varepsilon,K_n}^{(1,n)}(a_j)\Big)^p\\
		& \quad -C_{K_n,p}^p
		\sum_{1\leq i<j\leq n}
		\Big(f_{\varepsilon,K_n}^{(1,n)}(a_i)-f_{\varepsilon,K_n}^{(1,n)}(a_j)\Big)^p.
	\end{align*}
For sufficiently small $|\varepsilon|$, the remainder is defined by the same difference of the two energies, with absolute values where necessary. Since $p\geq3$, it is differentiable, is nonpositive and vanishes at $0$. Fermat's theorem therefore gives $R(f_{\varepsilon,K_n}^{(1,n)})'(0)=0$, that is,
	\begin{align}\label{R0}
		C_{K_n,p}^p=\frac{(f(a_1)-m_f)^{p-1}+(m_f-f(a_n))^{p-1}}
		{\sum_{j=1}^{n}\Big((f(a_1)-f(a_j))^{p-1}+(f(a_j)-f(a_n))^{p-1}\Big)}.
	\end{align}

Define the function $g:V_{K_n}\to\mathbb{R}$ as 
\begin{equation*}
		g(a_j)=
		\begin{cases}
			f(a_1) & {j=1,}\\
			\frac{1}{n-2}\sum_{j=2}^{n-1}f(a_j) & {2\leq j\leq n-1,}\\
			f(a_n) & {j=n,}
		\end{cases}
\end{equation*}
	then $m_g=m_f$. For $p\in[3,4]$, by Karamata's inequality, we have 
\begin{align*}
		&\sum_{j=1}^{n}\Big((f(a_1)-f(a_j))^{p-1}+(f(a_j)-f(a_n))^{p-1}\Big)\\
		&\geq
		\sum_{j=1}^{n}\Big((g(a_1)-g(a_j))^{p-1}+(g(a_j)-g(a_n))^{p-1}\Big).
\end{align*}
	Thus, it follows that 
\begin{align*}
        &\frac{(f(a_1)-m_f)^{p-1}+(m_f-f(a_n))^{p-1}}
		{\sum_{j=1}^{n}\Big((f(a_1)-f(a_j))^{p-1}+(f(a_j)-f(a_n))^{p-1}\Big)}\\
		&\leq
		\frac{(g(a_1)-m_g)^{p-1}+(m_g-g(a_n))^{p-1}}
		{\sum_{j=1}^{n}\Big((g(a_1)-g(a_j))^{p-1}+(g(a_j)-g(a_n))^{p-1}\Big)},
\end{align*}
implying that 
\begin{align*}\label{eq666-1}
C_{K_n,p}^p	
&\leq
\sup_{h\in \Lambda}\frac{(h(a_1)-m_h)^{p-1}+(m_h-h(a_n))^{p-1}}
{\sum_{j=1}^{n}\Big((h(a_1)-h(a_j))^{p-1}+(h(a_j)-h(a_n))^{p-1}\Big)} \nonumber\\
&=
\sup_{h\in \Lambda_1 \cup\Lambda_2 \cup \Lambda_3}
\frac{(h(a_1)-m_h)^{p-1}+(m_h-h(a_n))^{p-1}}
{\sum_{j=1}^{n}\Big((h(a_1)-h(a_j))^{p-1}+(h(a_j)-h(a_n))^{p-1}\Big)}
\end{align*}
where
\begin{equation*}
		\begin{cases}
			\Lambda:=\{h: h(a_1)>h(a_n),\ h(a_1)\geq h(a_2)=\cdots=h(a_{n-1})\geq h(a_n)\}, \\
			\Lambda_1:=\{h: h(a_1)> h(a_2)=\cdots=h(a_{n-1})= h(a_n)\}, \\
			\Lambda_2:=\{h: h(a_1)> h(a_2)=\cdots=h(a_{n-1})> h(a_n)\}, \\
            \Lambda_3:=\{h: h(a_1)= h(a_2)=\cdots=h(a_{n-1})>h(a_n)\}.
		\end{cases}
\end{equation*}

By the translation and dilation invariance of the quotient in \eqref{R0}, we may take $h(a_1)=1$, $h(a_2)=x$ and $h(a_n)=0$, where $x\in [0,1]$. Therefore we set
\begin{equation*}
\begin{cases}
\Lambda_1:=\{h: h(a_1)=1, h(a_2)=\cdots=h(a_{n-1})= h(a_n)=0\}, \\
\Lambda_2:=\{h: h(a_1)=1, h(a_2)=\cdots=h(a_{n-1})=x\in (0,1), h(a_n)=0\}, \\
\Lambda_3:=\{h: h(a_1)= h(a_2)=\cdots=h(a_{n-1})=1, h(a_n)=0\}.
\end{cases}
\end{equation*}

(ii) Assume $p\in(3,4)$ and  $n=3$. For $h\in\Lambda$, by a simple calculation and \eqref{eq6333} from Lemma \ref{NOMAX}, it follows that
\begin{align*}
	&\sup_{h\in \Lambda}
	\frac{(h(a_1)-m_h)^{p-1}+(m_h-h(a_3))^{p-1}}
	{\sum_{j=1}^{3}\Big((h(a_1)-h(a_j))^{p-1}+(h(a_j)-h(a_3))^{p-1}\Big)} \nonumber\\
	&=
	\max_{x\in[0,1]}
	\frac{\big(1-\frac{1+x}{3}\big)^{p-1}+\big(\frac{1+x}{3}\big)^{p-1}}
	{(1-x)^{p-1}+x^{p-1}+2}\\
	&=\frac{1}{2^{p-1}+1}.	
\end{align*}
Hence $C_{K_3,p}^p\leq\frac{1}{2^{p-1}+1}$. 

The lower bound for $C_{K_3,p}^p$ is obtained from the choice 	$\widetilde{f}(a_1)=1$, $\widetilde{f}(a_2)=1/2$ and $\widetilde{f}(a_3)=0$, yielding
\[
C_{K_3,p}^p\geq\frac{\Vert \widetilde{f}-m_{\widetilde{f}}\Vert_p^p}{{\rm Var}_{p}^{K_3}(\widetilde{f})^{p}}=\frac{1}{2^{p-1}+1}.
\]
Hence $C_{K_3,p}^p=\frac{1}{2^{p-1}+1}$ is sharp.

The invariance of the quotient in \eqref{R0} under translations and dilations allows us to normalize $f:V_{K_3}\to\mathbb{R}$ so that
\[
(f(a_1),f(a_2),f(a_3))=(a+b,a,a-b)
\]
for some $a,b\in\mathbb{R}$ with $b\neq0$.
	
We now prove (iii) and (iv). Assume $p\in(3,4)$ and $n\geq4$. For $h\in\Lambda$, a direct calculation yields
\begin{align*}
	&\sup_{h\in \Lambda}
	\frac{(h(a_1)-m_h)^{p-1}+(m_h-h(a_n))^{p-1}}
	{\sum_{j=1}^{n}\Big((h(a_1)-h(a_j))^{p-1}+(h(a_j)-h(a_n))^{p-1}\Big)} \nonumber\\
	&=\frac{1}{n^{p-1}}\max_{x\in[0,1]}	\psi_n(x).
\end{align*}

From Lemma \ref{NOMAX}, we obtain the upper bound of $C_{K_n,p}^p$. The maximum is attained at $x=1/2$ for $p\in(3,p_n]$, and at $x=0$ or $x=1$ for $p\in[p_n,4)$. At $p=p_n$, all three values $0$, $1/2$ and $1$ attain the maximum.

The lower bound can be obtained by taking $\widetilde{f}(a_1)=1$, $\widetilde{f}(a_2)=\cdots=\widetilde{f}(a_{n-1})=1/2$, $\widetilde{f}(a_n)=0$, or $\bar{f}(a_1)=1$, $\bar{f}(a_2)=\cdots=\bar{f}(a_{n-1})=\bar{f}(a_n)=0$. Then we have
\[
C_{K_n,p}^p
\geq
\max\Big\{\frac{1}{2^{p-1}+n-2},\frac{1+(n-1)^{p-1}}{n^p}\Big\}.
\]
The case $p=4$ in (iv) follows by continuity: apply the inequality for $p<4$ to each fixed function and let $p\to4^-$. The Dirac delta functions give equality at the limiting constant as well. This completes the proof.
\end{proof}

\section{Sharp Poincar\'{e}-Wirtinger inequalities on star graphs}\label{s3}

In this section, we study the sharp constants in the Poincar\'{e}-Wirtinger inequalities on star graphs. Before obtaining the sharp constants on star graphs, we first characterize the extremizers on star graphs.

\begin{lem}\label{QGEQ1}
Let $S_n$ be the star graph with $n\geq4$ vertices $a_1$, $a_2,\ldots, a_n$ where $a_1$ is the center. Then 
\[
C_{S_n,p}>1
\]
for any $p\in(1,2)\cup(2,+\infty)$.
\end{lem}

\begin{proof}
Set $t\in(0,\frac{\vert p-2\vert}{p-1}]$. For $1<p<2$ and $n\geq5$, take
\begin{equation*}
		f_1(a_j)=
		\begin{cases}
			0 & {j=1,}\\
            2+t & {j=2,}\\
			0 & {3\leq j\leq n-2,}\\
			-1 & {n-1\leq j\leq n.}
		\end{cases}
	\end{equation*}
From Bernoulli's inequality, we derive that
\begin{align*}
&\Vert f_1-m_{f_1}\Vert_p^p-{\rm Var}_{p}^{S_n}(f_1)^p\\
=&
\frac{(n-3)t^p}{n^p}+2\Big(1+\frac{t}{n}\Big)^p+\Big(2+t-\frac{t}{n}\Big)^p-(2+t)^p-2\\
\geq&
\frac{(n-3)t^p}{n^p}+2+\frac{2pt}{n}+(2+t)^p-\frac{pt(2+t)^{p-1}}{n}-(2+t)^p-2\\
=&
\frac{(n-3)t^p}{n^p}+\frac{pt(2-(2+t)^{p-1})}{n}.
\end{align*}
Since $0<t\leq\frac{2-p}{p-1}$, it follows that
\begin{align*}
(2+t)^{p-1}
\leq
1+(p-1)(1+t)
\leq
1+(p-1)\Big(1+\frac{2-p}{p-1}\Big)=2,
\end{align*}
implying that 
\[
C_{S_n,p}\geq
\frac{\Vert f_1-m_{f_1}\Vert_p}{{\rm Var}_{p}^{S_n}(f_1)}>1.
\]

For $p>2$ and $n\geq5$, take
\begin{equation*}
		f_2(a_j)=
		\begin{cases}
			0 & {j=1,}\\
            1+t & {j=2,}\\
            1 & {j=3,}\\
			0 & {4\leq j\leq n-1,}\\
			-2 & {j=n.}
		\end{cases}
	\end{equation*}
It follows from Bernoulli's inequality that
\begin{align*}
&\Vert f_2-m_{f_2}\Vert_p^p-{\rm Var}_{p}^{S_n}(f_2)^p\\
=&
\frac{(n-3)t^p}{n^p}+\Big(2+\frac{t}{n}\Big)^p+\Big(1-\frac{t}{n}\Big)^p+\Big(1+t-\frac{t}{n}\Big)^p
-(1+t)^p-1-2^p\\
\geq&
\frac{(n-3)t^p}{n^p}+2^p+\frac{2^{p-1}pt}{n}+1-\frac{pt}{n}+(1+t)^p-\frac{pt(1+t)^{p-1}}{n}
-(1+t)^p-1-2^p\\
=&
\frac{(n-3)t^p}{n^p}+\frac{pt}{n}(2^{p-1}-1-(1+t)^{p-1}).
\end{align*}
From the convexity of $x\mapsto x^{p-1}$, we have
\begin{align*}
(1+t)^{p-1}\leq(1-t)+t2^{p-1}=1+t(2^{p-1}-1)
\leq
1+\frac{p-2}{p-1}(2^{p-1}-1).
\end{align*}
It then follows that
\begin{align*}
2^{p-1}-1-(1+t)^{p-1}
&\geq
2^{p-1}-1-\Big(1+\frac{p-2}{p-1}(2^{p-1}-1)\Big)\\
&=
\frac{2^{p-1}-p}{p-1}>0,
\end{align*}
implying that 
\[
C_{S_n,p}\geq
\frac{\Vert f_2-m_{f_2}\Vert_p}{{\rm Var}_{p}^{S_n}(f_2)}>1.
\]

If $n=4$, then by taking the value 	
\[
	\big(f_1^\ast(a_1),f_1^\ast(a_2),f_1^\ast(a_{3}), f_1^\ast(a_4)\big)=\big(0,2+t,-1,-1\big)
	\]
when $1<p<2$ and  
\[
	\big(f_2^\ast(a_1),f_2^\ast(a_2),f_2^\ast(a_{3}), f_2^\ast(a_4)\big)=\big(0,1+t,1,-2\big)
	\]
when $2<p<+\infty$, we can similarly derive that 
$
C_{S_4,p}>1
$
for any $p\in(1,2)\cup(2,+\infty)$. This completes the proof.
\end{proof}

\begin{prop}\label{EOSN1}
Let $S_n$ be the star graph with $n\geq3$ vertices $a_1$, $a_2,\ldots, a_n$ where $a_1$ is the center, and let $0<p<+\infty$. There exists $f:V_{S_n}\to\mathbb{R}$ with
\begin{align*}
C_{S_n,p}
=\frac{\Vert f-m_f\Vert_p}{{\rm Var}_{p}^{S_n}(f)},
\end{align*}
such that $\min\limits_{2\leq j\leq n}f(a_j)\leq f(a_1)\leq \max\limits_{2\leq j\leq n}f(a_j)$.
Moreover for $n\geq4$, if $1<p<+\infty$ and $p\neq2$, then $f_{\vert V_{S_n}\setminus\{a_1\}}$ takes at least 2 different values, and at most 3 different values.
\end{prop}

\begin{proof}
By Lemma \ref{EX}, there exists $f$ such that
\begin{align}\label{EXTRESN}
C_{S_n,p}
=\frac{\Vert f-m_f\Vert_p}{{\rm Var}_{p}^{S_n}(f)}.
\end{align}
Since the quotient on the right hand side of \eqref{EXTRESN} is translation and dilation invariant with respect to $f$, we can assume that $f(a_1)=0$, and $m_f\geq0$. 
We now proceed in several steps.

\textit{\textbf{Step 1:} We prove that $f_{\vert V_{S_n}\setminus\{a_1\}}$  takes at least 2 different values.} Suppose that $f$ satisfies \eqref{EXTRESN}, and $f(a_i)=f(a_j)$ for any $2\leq i\leq j\leq n$, then for $1<p<+\infty$,
\[
C_{S_n,p}=\frac{(1+(n-1)^{p-1})^{1/p}}{n}.
\]
However, if we set $\widetilde{g}(a_1)=0$, $\widetilde{g}(a_2)=1$, $\widetilde{g}(a_3)=\cdots=\widetilde{g}(a_n)=0$, then it follows that 
\[
\frac{(1+(n-1)^{p-1})^{1/p}}{n}
=
C_{S_n,p}
\geq
\frac{\Vert \widetilde{g}-m_{\widetilde{g}}\Vert_p}{{\rm Var}_{p}^{S_n}(\widetilde{g})}
=
\frac{(n-1+(n-1)^{p})^{1/p}}{n},
\]
which is a contradiction.

\textit{\textbf{Step 2:} We prove that $\min\limits_{2\leq j\leq n}f(a_j)\leq f(a_1)\leq \max\limits_{2\leq j\leq n}f(a_j)$.} We consider the following cases. 	
	
\textbf{Case 1.} Suppose that the maximum of $f$ is attained at $a_1$ and
\[
	f(a_1)> f(a_2)\ge f(a_3)\ge\cdots\ge f(a_{n-1})\ge f(a_n).
\]
We now consider the auxiliary function $g_1$ defined by
	\begin{equation*}
		g_1(a_j)=
		\begin{cases}
			f(a_2) & {j=1,}\\
			f(a_1) & {j=2,}\\
			f(a_j) & {3\leq j\leq n}.
		\end{cases}
	\end{equation*}
Then we have that $m_{g_1}=m_{f}$, $g_1(a_2)\geq g_1(a_j)$ for any $1< j\leq n$, and 
	\begin{align*}
		C_{S_n,p}^p
		\geq\frac{\Vert g_1-m_{g_1}\Vert_p^p}{{\rm Var}_{p}^{S_n}(g_1)^p}.
	\end{align*}

	For $1< p<+\infty$, since $\Vert g_1-m_{g_1}\Vert_p^p=\Vert f-m_f\Vert_p^p$ and
	\begin{align*}
		{\rm Var}_{p}^{S_n}(g_1)^p&=\sum_{j=2}^{n}\vert g_1(a_j)-g_1(a_1)\vert^p\\
		&<(f(a_1)-f(a_2))^p+\sum_{j=3}^{n}(f(a_1)-f(a_j))^p\\
		&={\rm Var}_{p}^{S_n}(f)^p,
	\end{align*} 
	it follows that
	\begin{align*}
		\frac{\Vert g_1-m_{g_1}\Vert_p^p}{{\rm Var}_{p}^{S_n}(g_1)^p}>\frac{\Vert f-m_f\Vert_p^p}{{\rm Var}_{p}^{S_n}(f)^p},
	\end{align*} 
	which contradicts \eqref{EXTRESN}.
	
\textbf{Case 2.} Suppose that $f$ attains its minimum at $a_1$ and that
	\[
	f(a_2)\geq f(a_3)\geq\cdots\geq f(a_{n-1})\geq f(a_n)>f(a_1). 
	\]
Following the same argument as in Case 1, we take
	\begin{equation*}
		g_2(a_j)=
		\begin{cases}
			f(a_n) & {j=1,}\\
			f(a_j) & {2\leq j\leq n-1,}\\
			f(a_1) & {j=n.}
		\end{cases}
	\end{equation*}
Then we have
\begin{align*}
\frac{\Vert g_2-m_{g_2}\Vert_p^p}{{\rm Var}_{p}^{S_n}(g_2)^p}
>
\frac{\Vert f-m_f\Vert_p^p}{{\rm Var}_{p}^{S_n}(f)^p},
\end{align*} 
which contradicts \eqref{EXTRESN}. Hence, we can assume that 
\[
\min\limits_{2\leq j\leq n}f(a_j)\leq f(a_1)\leq \max\limits_{2\leq j\leq n}f(a_j).
\]

\textit{\textbf{Step 3:} We prove that $f_{\vert V_{S_n}\setminus\{a_1\}}$  takes at most 3 different values.} If $n=3$, $4$, the conclusion holds trivially; for this reason, we consider the case $n\geq5$ below. For $1<p<+\infty$ and $p\neq2$, we consider the auxiliary function
\begin{align*}
L(x_2,\ldots,x_n)
:=\left\vert\frac{\sum_{j=2}^{n}x_j}{n}\right\vert^p+\sum_{j=2}^{n}\left\vert x_j-\frac{\sum_{j=2}^{n}x_j}{n}\right\vert^p-C_{S_n,p}^p\sum_{j=2}^{n}\vert x_j\vert^p.
\end{align*}
Note that $f$ is the extremizer satisfying \eqref{EXTRESN}. 
Let $x_k=f(a_k)$ for any $2\leq k\leq n$, then
\[
L(x_2,\ldots,x_n)=0, \quad  \frac{\partial}{\partial x_k}L(x_2,\ldots,x_n)=0, 
\] 
which is equivalent to 
\begin{align}\label{LAE}
&\vert m_f - f(a_k) \vert ^ { p - 2 } ( m_f - f(a_k) ) 
+C_{S_n,p}^p \vert f(a_k) \vert ^ { p - 2 } f(a_k) \nonumber\\
&\qquad= \frac { \vert m_f \vert ^ { p - 2 } m_f + \sum _ { j = 2 } ^ { n } \vert m_f - f(a_j) \vert ^ { p - 2 } (m_f - f(a_j)) } { n }\\
&\qquad=:K_f \nonumber.
\end{align}

We claim that the equation \eqref{LAE} in terms of $f(a_k)$ has at most 3 different roots for $1<p<\infty$, $p\neq2$ and any $m_f$, $C_{S_n,p}$, $K_f\in\mathbb{R}$. In fact, we define the function 
\begin{equation}\label{phi}
\phi(x):=
\begin{cases}
			\vert m_f - x \vert ^ { p - 2 } ( m_f - x ) 
+C_{S_n,p}^p \vert x \vert ^ { p - 2 } x & {x\in\mathbb{R}\setminus\{0, m_f\},}\\
			C_{S_n,p}^p \vert m_f \vert ^ { p - 2 } m_f & {x=m_f,}\\
			\vert m_f \vert ^ { p - 2 }m_f   & {x=0.}
\end{cases}
\end{equation}
Notice that $\phi$ is continuous on $\mathbb{R}$, and without loss of generality, we can assume that $m_f\geq0$, then we have that 
\[
\phi'(x)
=
(p-1)\Big(C_{S_n,p}^p\vert x\vert^{p-2}-\vert m_f-x\vert^{p-2}\Big),
\quad x\in\mathbb{R}\setminus\{0, m_f\}.
\]
Then it follows that $\phi'(x)=0$ is equivalent to 
\[
\vert m_f-x\vert=(C_{S_n,p})^{\frac{p}{p-2}}\vert x\vert.
\]
Thus we have that 
\[
\Big(1-(C_{S_n,p})^{\frac{2p}{p-2}}\Big)x^2-2m_f x+m_f^2=0.
\]

The roots of the quadratic equation are $x=m_f(1-(C_{S_n,p})^{\frac{p}{p-2}})^{-1}$ and $x=m_f(1+(C_{S_n,p})^{\frac{p}{p-2}})^{-1}$. Hence, the equation \eqref{LAE} has at most 3 different roots for $1<p<\infty$, $p\neq2$ and any $m_f$, $K_f\in\mathbb{R}$.
This completes the proof.
\end{proof}

\begin{prop}\label{EOSN2}
Let $S_n$ be the star graph with $n\geq4$ vertices $a_1$, $a_2,\ldots, a_n$ where $a_1$ is the center. For $1<p<2$, if
\begin{align*}
C_{S_n,p}
=\frac{\Vert f-m_f\Vert_p}{{\rm Var}_{p}^{S_n}(f)},
\end{align*}
then $f_{\vert V_{S_n}\setminus\{a_1\}}$ only takes 2 values not equal to $f(a_1)$.
\end{prop}
\begin{proof}
By Proposition \ref{EOSN1} we can assume that $f$ is such that
$C_{S_n,p}=\frac{\Vert f-m_f\Vert_p}{{\rm Var}_{p}^{S_n}(f)}$.
Moreover, we may assume that $f(a_1)=0$, and assume that there exists $r\in[2,n-1]$ such that
\[
f(a_2)\geq\cdots\geq f(a_r)\geq 0\geq f(a_{r+1})\geq\cdots\geq f(a_n).
\]

If $m_f=0$, then we would get $\frac{\Vert f-m_f\Vert_p}{{\rm Var}_{p}^{S_n}(f)}=1$, contradicting Lemma \ref{QGEQ1}. Hence, we may assume that $m_f>0$ without loss of generality.
We now proceed in several steps.

\textit{\textbf{Step 1:} We prove that $f(a_{r+1})=\cdots=f(a_n)<0$.} Assume that for some $r+1\leq k<l\leq n$, $f(a_k)\neq f(a_l)$.
For $\varepsilon>0$ sufficiently small and $1\leq j\leq n$, choose the function $f_{\varepsilon}^{(s,t)}:V_{S_n}\to\mathbb{R}$ given by
\begin{equation}\label{RDF}
		f_{\varepsilon,S_n}^{(s,t)}(a_j):=
		\begin{cases}
			f(a_s)+\varepsilon & {j=s,}\\
			f(a_t)-\varepsilon & {j=t,}\\
			f(a_j) & {j\neq s,t,}
		\end{cases}
\end{equation}
and define the remainder functions
\begin{align*}
		R(f_{\varepsilon,S_n}^{(2,j)})(\varepsilon)
        &:=
        \Vert f_{\varepsilon,S_n}^{(2,j)}-m_{f_{\varepsilon,S_n}^{(2,j)}} \Vert_p^p
        -C_{S_n,p}^p{\rm Var}_{p}^{S_n}(f_{\varepsilon,S_n}^{(2,j)})^p,
\end{align*} 
for $j=k,l$ respectively.
For sufficiently small $|\varepsilon|$, define both remainders by the same formulas. Since $p>1$, they are differentiable, are nonpositive and vanish at $0$. Thus $R(f_{\varepsilon,S_n}^{(2,k)})'(0)=R(f_{\varepsilon,S_n}^{(2,l)})'(0)=0$ and 
\begin{align*} 
C_{S_{n},p}^p
&=
\frac{\left(f(a_{2})-m_{f}\right)^{p-1}+\left(m_{f}-f(a_{k})\right)^{p-1}}
{f(a_{2})^{p-1}+\left(-f(a_{k})\right)^{p-1}}\\
&=
\frac{\left(f(a_{2})-m_{f}\right)^{p-1}+\left(m_{f}-f(a_{l})\right)^{p-1}}
{f(a_{2})^{p-1}+\left(-f(a_{l})\right)^{p-1}}.
\end{align*}

We define the function $\mathfrak{G}(x)$ for $x\geq 0$ by
\[
\mathfrak{G}(x):=
\frac{\bigl(f(a_{2})-m_{f}\bigr)^{p-1}+\bigl(m_{f}+x\bigr)^{p-1}}
{f(a_{2})^{p-1}+x^{p-1}}.
\]
Then $x=-f(a_k)$ and $x=-f(a_l)$ are two roots of $\mathfrak{G}(x)-C_{S_{n},p}^p=0$. Since $C_{S_{n},p}>1$ for any $p\in(1,2)$ and $n\geq4$, at either of these roots it follows that 
\begin{align*}
&\frac{\bigl(f(a_{2})-m_{f}\bigr)^{p-1}+\bigl(m_{f}+x\bigr)^{p-1}}
{f(a_{2})^{p-1}+x^{p-1}}-1\\
=&
\frac{(m_f+x)^{p-1}-x^{p-1}-(f(a_2)^{p-1}-(f(a_2)-m_f)^{p-1})}{f(a_{2})^{p-1}+x^{p-1}}
>0,
\end{align*}
implying that 
$(m_f+x)^{p-1}-x^{p-1}>f(a_2)^{p-1}-(f(a_2)-m_f)^{p-1}$. From the monotonicity of $x\mapsto (m_f+x)^{p-1}-x^{p-1}$ for $p\in(1,2)$ and $x\geq0$, we know that $x\in[0,f(a_2)-m_f)\subset[0,f(a_2))$. Thus, we have that 
\begin{align*}
&\frac{(f(a_2)^{p-1}+x^{p-1})^2}{p-1}\mathfrak{G}'(x)\\
=&
(m_f+x)^{p-2}(f(a_2)^{p-1}-m_fx^{p-2})-x^{p-2}(f(a_2)-m_f)^{p-1}\\
\leq&
(m_f+x)^{p-2}(f(a_2)^{p-1}-m_f f(a_2)^{p-2})-x^{p-2}(f(a_2)-m_f)^{p-1}\\
=&
f(a_2)(f(a_2)-m_f)\Big((f(a_2)(m_f+x))^{p-2}-(x(f(a_2)-m_f))^{p-2}\Big)<0,
\end{align*}
which implies, by continuity at $0$ and the derivative estimate for $0<x<f(a_2)-m_f$, that $\mathfrak{G}$ is monotonically decreasing on $x\in[0,f(a_2)-m_f)$ and the equation $\mathfrak{G}(x)-C_{S_{n},p}^p=0$ has a unique root on $[0,f(a_2)-m_f)$; hence $f(a_k)=f(a_l)$. Since $k$ and $l$ are arbitrary, we obtain 
$f(a_{r+1})=\cdots=f(a_n)=-c\leq0$.

We claim $c>0$. In fact, Suppose that $c=0$, then for $\varepsilon\geq0$ sufficiently small, choose the function $\widehat{f}_{\varepsilon}:V_{S_n}\to\mathbb{R}$ given by
\begin{equation*}
	\widehat{f}_{\varepsilon}(a_j)=
	\begin{cases}
		0  & {j=1,}\\
		a  & {2\leq j\leq u+1,}\\
		(1-\varepsilon)b  & {u+2\leq j\leq u+v+1,}\\
		-\frac{va}{w}\varepsilon & {u+v+2\leq j\leq n,}
	\end{cases}
\end{equation*}
and define the quotient function
\begin{align*}
	\widehat{Q}(\varepsilon)
	:=
	\frac{\sum_{j=1}^{n}\vert \widehat{f}_{\varepsilon}(a_j)-m_{\widehat{f}_{\varepsilon}}\vert^p}
	{\sum_{j=1}^{n}\vert\widehat{f}_{\varepsilon}(a_j)\vert^p}
	=:\frac{\widehat{N}(\varepsilon)}{\widehat{D}(\varepsilon)}.
\end{align*}
From the definition of $C_{S_n,p}$, $\widehat{Q}(0)=C_{S_n,p}^p$ and $\widehat{Q}(\varepsilon)\leq C_{S_n,p}^p$ for any sufficiently small $\varepsilon\geq0$. 
However, the following computation yields 
\[
\widehat{Q}'_+(0)>0,
\]
which contradicts $\widehat{Q}(\varepsilon)\leq \widehat{Q}(0)$. In fact, since $
\widehat{D}'_+(0)=-pvb^{p-1}\leq0$,
\begin{align*}
	\frac{\widehat{N}'_+(0)}{p}
	&=
	um_f(a-m_f)^{p-1}-v(m_f-b)^p-m_f^p+w\Big(\frac{ua}{w}-m_f\Big)m_f^{p-1}\\
	&=
	um_f(a-m_f)^{p-1}-v(m_f-b)^p+m_f^{p-1}(ua-m_f(w+1))\\
	&=
	um_f(a-m_f)^{p-1}-v(m_f-b)^p+m_f^{p-1}(nm_f-vb-m_f(w+1))\\
	&=
	um_f(a-m_f)^{p-1}-v(m_f-b)^p+m_f^{p-1}(um_f+v(m_f-b))\\
	&=
	um_f((a-m_f)^{p-1}+m_f^{p-1})+v(m_f-b)(m_f^{p-1}-(m_f-b)^{p-1})>0,
\end{align*}
leading to 
\begin{align*}
	\widehat{Q}'_+(0)
	=\frac{\widehat{N}'_+(0)\widehat{D}(0)-\widehat{N}(0)\widehat{D}'_+(0)}
	{\widehat{D}(0)^2}>0,
\end{align*}
which contradicts. Hence $c>0$.

Now we can assume that $f_{\vert V_{P_n}\setminus\{a_1\}}$ takes three values, let us say $a$ $u$-times, $b$ $v$-times and $-c$ $w$-times, with $a>m_f\geq b\geq0>-c$, that is 
\begin{equation*}
	f(a_j)=
	\begin{cases}
		0  & {j=1,}\\
		a  & {2\leq j\leq u+1,}\\
		b  & {u+2\leq j\leq u+v+1,}\\
		-c & {u+v+2\leq j\leq n.}
	\end{cases}
\end{equation*}

\textit{\textbf{Step 2:} Under the assumption $a>b$, we claim $b=0$.} In fact, suppose that $0<b\leq m_f$, then for $\varepsilon\geq0$ sufficiently small, choose the function $\breve{f}_{\varepsilon}:V_{S_n}\to\mathbb{R}$ given by
\begin{equation*}
\breve{f}_{\varepsilon}(a_j)=
\begin{cases}
0  & {j=1,}\\
a  & {2\leq j\leq u+1,}\\
b-\varepsilon  & {u+2\leq j\leq u+v+1,}\\
-c+\frac{v}{w}\varepsilon & {u+v+2\leq j\leq n,}
\end{cases}
\end{equation*}
and define the remainder function
\begin{align*}
R(\breve{f}_{\varepsilon})(\varepsilon)
:=
\sum_{j=1}^{n}\vert \breve{f}_{\varepsilon}(a_j)-m_{\breve{f}_{\varepsilon}}\vert^p
-
C_{S_n,p}^p\sum_{j=1}^{n}\vert\breve{f}_{\varepsilon}(a_j)\vert^p.
\end{align*}
Since $R(\breve{f}_{\varepsilon})(0)=0$ and $R(\breve{f}_{\varepsilon})(\varepsilon)\leq 0$ for any sufficiently small $\varepsilon\geq0$, it follows that $R(\breve{f}_{\varepsilon})'_+(0)\leq0$. 
However, since $C_{S_n,p}^p>1$, then
\begin{align*}
R(\breve{f}_{\varepsilon})'_+(0)
&=
pv(m_f-b)^{p-1}-pv(m_f+c)^{p-1}+pvC_{S_n,p}^p(b^{p-1}+c^{p-1})\\
&>
pv\Big(b^{p-1}+(m_f-b)^{p-1}+c^{p-1}-(m_f+c)^{p-1}\Big)\\
&>
0,
\end{align*}
where $b^{p-1}+(m_f-b)^{p-1}\geq m_f^{p-1}$ and $(m_f+c)^{p-1}<m_f^{p-1}+c^{p-1}$. This is a contradiction. Thus $b$ cannot be positive. 

Therefore, we may assume that   
\[
\Big(f(a_1),f(a_2),\ldots,f(a_{n-1}),f(a_n)\Big)
=
\left(
0,
\underbrace{a,\ldots,a}_{u\text{-times}},
\underbrace{0,\ldots,0}_{v\text{-times}},
\underbrace{-c,\ldots,-c}_{w\text{-times}}
\right).
\]

\textit{\textbf{Step 3:} We prove that no leaf has value $0$.} If $v=0$, the conclusion is already established. Suppose $v>0$, and consider
\begin{equation*}
\widetilde{f}(a_j)=
\begin{cases}
0  & {j=1,}\\
a  & {2\leq j\leq u+1,}\\
-\frac{wc}{v+w} & {u+2\leq j\leq n.}
\end{cases}
\end{equation*}
Observe that
\begin{align*}
\frac{\Vert \widetilde{f}-m_{\widetilde{f}}\Vert_p^p}{{\rm Var}_{p}^{S_n}(\widetilde{f})^p}
&=
\frac{m_f^p+u(a-m_f)^p+(v+w)\Big(m_f+\frac{wc}{v+w}\Big)^p}
{ua^p+(v+w)\Big(\frac{wc}{v+w}\Big)^p}\\
&=
\frac{\Vert f-m_f\Vert_p^p-\Big(vm_f^p+w(c+m_f)^p-(v+w)\Big(m_f+\frac{wc}{v+w}\Big)^p\Big)}
{{\rm Var}_{p}^{S_n}(f)^p-\Big(wc^p-(v+w)\Big(\frac{wc}{v+w}\Big)^p\Big)}\\
&=:
\frac{\Vert f-m_f\Vert_p^p-\Delta_N}
{{\rm Var}_{p}^{S_n}(f)^p-\Delta_D}.
\end{align*}
Then we have that
\begin{align*}
\frac{\Vert \widetilde{f}-m_{\widetilde{f}}\Vert_p^p}{{\rm Var}_{p}^{S_n}(\widetilde{f})^p}
-
\frac{\Vert f-m_f\Vert_p^p}{{\rm Var}_{p}^{S_n}(f)^p}
&=
\frac{\Vert f-m_f\Vert_p^p-\Delta_N}
{{\rm Var}_{p}^{S_n}(f)^p-\Delta_D}
-
\frac{\Vert f-m_f\Vert_p^p}{{\rm Var}_{p}^{S_n}(f)^p}\\
&=
\frac{-{\rm Var}_{p}^{S_n}(f)^p\Delta_N+\Vert f-m_f\Vert_p^p\Delta_D }
{{\rm Var}_{p}^{S_n}(f)^p({\rm Var}_{p}^{S_n}(f)^p-\Delta_D)}\\
&=
\frac{\Vert f-m_f\Vert_p^p}{{\rm Var}_{p}^{S_n}(f)^p}
\frac{\Delta_D}{{\rm Var}_{p}^{S_n}(f)^p-\Delta_D}
-
\frac{\Delta_N}{{\rm Var}_{p}^{S_n}(f)^p-\Delta_D}\\
&>
\frac{\Delta_D-\Delta_N}{{\rm Var}_{p}^{S_n}(f)^p-\Delta_D}.
\end{align*}

We claim that $\Delta_D\geq\Delta_N$. In fact, since $x\mapsto x^{p-1}$ is concave on $[0,+\infty)$, it follows that 
\begin{align*}
\Big(x+\frac{wc}{v+w}\Big)^{p-1}\geq\frac{v}{v+w}x^{p-1}+\frac{w}{v+w}(x+c)^{p-1}
\end{align*}
holds for any $x\in[0,m_f]$. Then
\begin{align*}
(v+w)\int_{0}^{m_f}\Big(x+\frac{wc}{v+w}\Big)^{p-1}dx
\geq
v\int_{0}^{m_f}x^{p-1}dx+w\int_{0}^{m_f}(x+c)^{p-1}dx,
\end{align*}
which leads to 
\begin{align*}
(v+w)\left(\Big(m_f+\frac{wc}{v+w}\Big)^{p}-\Big(\frac{wc}{v+w}\Big)^{p} \right)
\geq
vm_f^p+w\left((c+m_f)^p-c^p\right).
\end{align*}
Thus
\begin{align*}
wc^p-(v+w)\Big(\frac{wc}{v+w}\Big)^p
\geq
vm_f^p+w(c+m_f)^p-(v+w)\Big(m_f+\frac{wc}{v+w}\Big)^p.
\end{align*}
Hence, we derive that 
\[
\frac{\Vert \widetilde{f}-m_{\widetilde{f}}\Vert_p^p}{{\rm Var}_{p}^{S_n}(\widetilde{f})^p}
>
\frac{\Vert f-m_f\Vert_p^p}{{\rm Var}_{p}^{S_n}(f)^p}=C_{S_n,p}^p,
\]
which contradicts. This completes the proof.
\end{proof}

For $2<p<\infty$, we aim to establish a result analogous to that in the case $1<p<2$. However, since the convexity properties of $x\mapsto x^{p-2}$ change, the proof requires a different approach. To this end, we first prove the following lemma.
\begin{lem}\label{TAIN}
Let $b,c\in\mathbb{R}$ satisfying $m_f>0>-b>-c$, and $c>m_f>0$. Then
\begin{align}\label{TAIN1}
\frac{c^{p-1}(m_f+b)^{p-1}-b^{p-1}(m_f+c)^{p-1}}{c^{p-1}-b^{p-1}}
>
\frac{m_f(b+c)(m_f+b)^{p-2}(m_f+c)^{p-2}}{b(m_f+b)^{p-2}+c(m_f+c)^{p-2}}  
\end{align}
holds for any $2<p<+\infty$.
\end{lem}

\begin{proof}
Let $x:=\frac{b}{m_f+b}$ and $y:=\frac{c}{m_f+c}$. Then $b=\frac{x m_f }{1-x}$, $c=\frac{y m_f}{1-y}$, $m_f+b=\frac{m_f}{1-x}$ and $m_f+c=\frac{m_f}{1-y}$. Thus the left hand side of \eqref{TAIN1} becomes 
\begin{align*}
&\frac{c^{p-1}(m_f+b)^{p-1}-b^{p-1}(m_f+c)^{p-1}}{c^{p-1}-b^{p-1}}\\
=&
\frac{\Big(\frac{y m_f}{1-y}\Big)^{p-1}\Big(\frac{m_f}{1-x}\Big)^{p-1}-\Big(\frac{x m_f}{1-x}\Big)^{p-1}\Big(\frac{m_f}{1-y}\Big)^{p-1}}
{\Big(\frac{y m_f}{1-y}\Big)^{p-1}-\Big(\frac{x m_f}{1-x}\Big)^{p-1}}\\
=&
m_f^{p-1}\frac{y^{p-1}-x^{p-1}}{(y(1-x))^{p-1}-(x(1-y))^{p-1}},
\end{align*}
and the right hand side of \eqref{TAIN1} becomes 
\begin{align*}
&\frac{m_f(b+c)(m_f+b)^{p-2}(m_f+c)^{p-2}}{b(m_f+b)^{p-2}+c(m_f+c)^{p-2}}\\
=&
\frac{m_f^2\Big(\frac{x+y-2xy}{(1-x)(1-y)}\Big)\Big(\frac{m_f}{1-x}\Big)^{p-2}\Big(\frac{m_f}{1-y}\Big)^{p-2}}
{\frac{x m_f}{1-x}\Big(\frac{m_f}{1-x}\Big)^{p-2}+\frac{y m_f}{1-y}\Big(\frac{m_f}{1-y}\Big)^{p-2}}\\
=&
m_f^{p-1}\frac{x+y-2xy}{x(1-y)^{p-1}+y(1-x)^{p-1}}.
\end{align*}

Since $m_f>0$, \eqref{TAIN1} is equivalent to 
\begin{align}\label{TAIN2}
\frac{y^{p-1}-x^{p-1}}{(y(1-x))^{p-1}-(x(1-y))^{p-1}}
>
\frac{x+y-2xy}{x(1-y)^{p-1}+y(1-x)^{p-1}},
\end{align}
for $0<x<y<1$, $1/2<y<1$, and $2<p<+\infty$. Let $A:=\frac{y}{x}$ and $B:=\frac{1-x}{1-y}$. 
One can check that $B>1$, $A>1$, and 
\[
x=\frac{B-1}{AB-1},\quad
y=\frac{A(B-1)}{AB-1},\quad
1-x=\frac{B(A-1)}{AB-1},\quad
1-y=\frac{A-1}{AB-1}.
\]

Then the left hand side of \eqref{TAIN2} becomes
\begin{align*}
\frac{y^{p-1}-x^{p-1}}{(y(1-x))^{p-1}-(x(1-y))^{p-1}}
&=
\frac{(AB-1)^{p-1}(A^{p-1}-1)}{((AB)^{p-1}-1)(A-1)^{p-1}},
\end{align*}
and the right hand side of \eqref{TAIN2} becomes
\begin{align*}
\frac{x+y-2xy}{x(1-y)^{p-1}+y(1-x)^{p-1}}
=
\frac{(AB-1)^{p-2}(AB+1)}{(A-1)^{p-2}(AB^{p-1}+1)}.
\end{align*}
Hence, \eqref{TAIN2} is equivalent to 
\begin{align}\label{TAIN3}
\frac{(A^{p-1}-1)(AB-1)(AB^{p-1}+1)}{(A-1)((AB)^{p-1}-1)(AB+1)}>1.
\end{align}

If $A=B$, then
\begin{align*}
\frac{(A^{p-1}-1)(AB-1)(AB^{p-1}+1)}{(A-1)((AB)^{p-1}-1)(AB+1)}
&=\frac{A^{p-1}-1}{A-1}\frac{(A^2-1)(A^p+1)}{(A^{2p-2}-1)(A^2+1)}\\
&=1+\frac{A(A-1)(A^{p-2}-1)}{(A^2+1)(A^{p-1}+1)}>1.
\end{align*}

If $A\ne B$, first assume that $B>A>1$; the case $A>B>1$ will be considered below. Rewrite \eqref{TAIN3} as
\[
\frac{A^{p-1}-1}{A-1}(AB^{p-1}+1)-\frac{(AB)^{p-1}-1}{AB-1}(AB+1)>0.
\]

Since $t\mapsto \frac{t^p-1}{t-1}$ is convex on $(1,+\infty)$ for $2<p<+\infty$, and 
\begin{align*}
B=\frac{AB-B}{AB-A}A+\frac{B-A}{AB-A}AB,
\end{align*}
it follows that
\begin{align*}
\frac{B^p-1}{B-1}
<\frac{AB-B}{AB-A}\Big(\frac{A^p-1}{A-1}\Big)+\frac{B-A}{AB-A}\Big(\frac{(AB)^p-1}{AB-1}\Big),
\end{align*}
which implies that 
\begin{align}\label{TAIN4}
\frac{B^{p-1}-A^{p-1}}{B-A}<\frac{(AB)^{p-1}-1}{AB-1}.
\end{align}

Since $t\mapsto \frac{t^{p-1}-1}{t-1}$ is monotonically increasing on $(1,+\infty)$ for $2<p<+\infty$, it follows that
\begin{align}\label{TAIN5}
\frac{B^{p-1}-1}{B-1}>\frac{A^{p-1}-1}{A-1}.
\end{align}
One can check that 
\begin{align*}
&\frac{A^{p-1}-1}{A-1}(AB^{p-1}+1)-\frac{(AB)^{p-1}-1}{AB-1}(AB+1)\\
=&
\frac{AB-2A+1}{A-1}\Big(\frac{(AB)^{p-1}-1}{AB-1}-\frac{B^{p-1}-A^{p-1}}{B-A}\Big)\\
&\qquad+\frac{(A-1)(B-1)}{B-A}
\Big(\frac{B^{p-1}-1}{B-1}-\frac{A^{p-1}-1}{A-1}\Big)>0.
\end{align*}
Since $y>{1}/{2}$, it follows that $AB>2A-1$, so the first coefficient is positive. Therefore, \eqref{TAIN4} and \eqref{TAIN5} imply \eqref{TAIN3} when $B>A$. If $A>B$, interchange $A$ and $B$ in the argument for \eqref{TAIN4}; that inequality remains unchanged. Both the inequality in \eqref{TAIN5} and the sign of $1/(B-A)$ reverse, so the second term in the same identity is still positive. Thus \eqref{TAIN3} holds in this case as well. This completes the proof.
\end{proof}

\begin{prop}\label{EOSN3}
Let $S_n$ be the star graph with $n\geq4$ vertices $a_1$, $a_2,\ldots, a_n$ where $a_1$ is the center. For $2<p<+\infty$, If
\begin{align*}
C_{S_n,p}
=\frac{\Vert f-m_f\Vert_p}{{\rm Var}_{p}^{S_n}(f)},
\end{align*}
then $f_{\vert V_{S_n}\setminus\{a_1\}}$ only takes 2 values not equal to $f(a_1)$.
\end{prop}

\begin{proof}
Analogously to the analysis of the case $p\in(1,2)$ in Proposition \ref{EOSN2}, we may assume that $f$ satisfies
\[
C_{S_n,p}=\frac{\lVert f-m_f\rVert_p}{\operatorname{Var}_{p}^{S_n}(f)},
\]
and that $f(a_1)=0$, $m_f>0$, and $f(a_2)\geq\cdots\geq f(a_n)$. Let $r\in[2,n]$ be the largest integer such that $f(a_r)\geq0$.

We proceed the proof into three steps.

\textit{\textbf{Step 1:} We can assume that $f(a_{2})=\cdots=f(a_r)>m_f$.} Since $f(a_1)=0$ and $m_f>0$, there must exist some $f(a_k)$, $2\leq k\leq r$ such that $f(a_k)>m_f$. It follows from the proof of Proposition \ref{EOSN1} that $x=f(a_k)$ is a root of the equation
\begin{align}\label{2PS1}
C_{S_n,p}^p x^{ p - 1 }-(x-m_f)^ { p - 1 }
=& K_f.
\end{align}
From the monotonicity of $\phi$ defined in \eqref{phi} on $[m_f,+\infty)$ for $p\in(2,+\infty)$, it follows that if there exists another $f(a_l)$, $2\leq l\leq r$ such that $f(a_l)>m_f$ and is a root of \eqref{2PS1} as well, then $f(a_k)=f(a_l)$.
Moreover, we have that
\begin{align*}
K_f
=\phi(f(a_k))
>
\phi(m_f)
=
C_{S_n,p}^pm_f^{p-1}.
\end{align*}

Suppose there exists $f(a_i)\in[0,m_f]$ for some $2\leq i\leq r$, then it follows that 
\begin{align}\label{LEQKF}
\phi(f(a_i))&=(m_f - f(a_i))^ { p - 1 }+C_{S_n,p}^p f(a_i)^{ p - 1 } \nonumber\\
&\leq
C_{S_n,p}^p\Big((m_f - f(a_i))^ { p - 1 }+f(a_i)^{ p - 1 }\Big)\nonumber\\
&\leq 
C_{S_n,p}^p m_f^{p-1}
<
K_f,
\end{align}
which implies that $f(a_i)$ cannot satisfy \eqref{LAE}. This indicates that all nonnegative leaf values of $f$ are equal and are greater than $m_f$.

Since $f|_{V_{S_n}\setminus\{a_1\}}$ takes at most three values by Proposition \ref{EOSN1}, we may now assume that it takes three values, say $a$, $-b$, and $-c$, with multiplicities $u$, $v$, and $w$, respectively, and with $a>m_f>0>-b>-c$; that is,
\[
\Big(f(a_1),f(a_2),\ldots,f(a_{n-1}),f(a_n)\Big)
=
\left(
0,
\underbrace{a,\ldots,a}_{u\text{-times}},
\underbrace{-b,\ldots,-b}_{v\text{-times}},
\underbrace{-c,\ldots,-c}_{w\text{-times}}
\right).
\]

\textit{\textbf{Step 2:} We can assume that $0\leq v<2$.} Suppose that $v\geq2$, and without loss of generality, we assume that $f(a_{r+1})=f(a_{r+2})=-b$, $f(a_l)=-c$, for some $u+v+2\leq l\leq n$.
For $\varepsilon>0$ sufficiently small such that $0<\varepsilon<\min\{b,c-b\}$, choose the function $f_{\varepsilon, S_n}^{(r+1,l)}$ defined in \eqref{RDF} and the remainder function:
\begin{align*}
R(f_{\varepsilon,S_n}^{(r+1,l)})(\varepsilon)
&:=
\Vert f_{\varepsilon,S_n}^{(r+1,l)}-m_{f_{\varepsilon,S_n}^{(r+1,l)}} \Vert_p^p
-C_{S_n,p}^p{\rm Var}_{p}^{S_n}(f_{\varepsilon,S_n}^{(r+1,l)})^p. 
\end{align*}
For sufficiently small $|\varepsilon|$, define the remainder by the same formula. It is differentiable and nonpositive, and $R(f_{\varepsilon,S_n}^{(r+1,l)})(0)=0$. Thus $R(f_{\varepsilon,S_n}^{(r+1,l)})'(0)=0$ and
\begin{align*}
C_{S_n,p}^p
=
\frac{(m_f+c)^{p-1}-(m_f+b)^{p-1}}{c^{p-1}-b^{p-1}}.
\end{align*}
Then from Cauchy's mean value theorem, there exist $\xi\in(b,c)$ such that
\begin{align}\label{2PS2}
\Big(\frac{m_f+c}{c}\Big)^{p-2}<C_{S_n,p}^p=\Big(\frac{m_f+\xi}{\xi}\Big)^{p-2}<\Big(\frac{m_f+b}{b}\Big)^{p-2}.
\end{align}

Define $\theta(x):=(m_f+x)^p-C_{S_n,p}^px^p$, $x\geq0$. Then from \eqref{2PS2}, we obtain that $\theta''(b)>0$, and
\begin{align*}
R(f_{\varepsilon,S_n}^{(r+1,r+2)})(\varepsilon)
&=
\Vert f_{\varepsilon}^{(r+1,r+2)}-m_{f_{\varepsilon}^{(r+1,r+2)}} \Vert_p^p
-C_{S_n,p}^p{\rm Var}_{p}^{S_n}(f_{\varepsilon}^{(r+1,r+2)})^p \\
&=
\Vert f_{\varepsilon}^{(r+1,r+2)}-m_{f_{\varepsilon}^{(r+1,r+2)}} \Vert_p^p
-\Vert f-m_{f} \Vert_p^p\\
&\qquad
-C_{S_n,p}^p\Big({\rm Var}_{p}^{S_n}(f_{\varepsilon}^{(r+1,r+2)})^p-{\rm Var}_{p}^{S_n}(f)^p\Big)\\
&=\theta(b+\varepsilon)+\theta(b-\varepsilon)-2\theta(b)>0,
\end{align*}
for all sufficiently small $0<\varepsilon<\min\{b,c-b\}$. This indicates that 
\[
\frac{\Vert f_{\varepsilon}^{(r+1,r+2)}-m_{f_{\varepsilon}^{(r+1,r+2)}} \Vert_p}
{{\rm Var}_{p}^{S_n}(f_{\varepsilon}^{(r+1,r+2)})} 
>C_{S_n,p},
\]
which contradicts to the definition of $C_{S_n,p}$. Hence $0\leq v<2$. 
If $v=0$, then the proof is finished. Thus we only to consider the case $v=1$. 

\textit{\textbf{Step 3:} We claim that $f$ cannot be an extremizer if $v=1$.} In fact, suppose that the extremizer $f$ has the form 
\begin{align}\label{ASSMP2P}
\Big(f(a_1),f(a_2),\ldots,f(a_{n-1}),f(a_n)\Big)
=
\left(
0,
\underbrace{a,\ldots,a}_{u\text{-times}},
-b,
\underbrace{-c,\ldots,-c}_{w\text{-times}}
\right),
\end{align}
with $a>m_f>0>-b>-c$, $1\leq u\leq n-2$, and $1\leq w\leq n-2$. Let 
\begin{equation*}
\begin{cases}
&\alpha:=C_{S_n,p}^p a^{p-2}-(a-m_f)^{p-2};\\
&\beta:=(m_f+b)^{p-2}-C_{S_n,p}^p b^{p-2};\\
&\gamma:=C_{S_n,p}^p c^{p-2}-(m_f+c)^{p-2}.\\
\end{cases}
\end{equation*}
From the fact $C_{S_n,p}>1$ and \eqref{2PS2}, we know that $\alpha$, $\beta$, $\gamma>0$.
Given $\varepsilon$ vary in a sufficiently small neighborhood of $0$ such that
\[
a+\frac{\gamma\varepsilon}{u\gamma+w\alpha}
>
-b-\varepsilon
>
-c+\frac{\alpha\varepsilon}{u\gamma+w\alpha},
\]
we can construct an auxiliary function $g_{\varepsilon}$ satisfying 
\begin{equation*}
g_{\varepsilon}(a_j)=
\begin{cases}
0  & {j=1,}\\
a+\frac{\gamma\varepsilon}{u\gamma+w\alpha}  & {2\leq j\leq u+1,}\\
-b-\varepsilon  & {j=u+2,}\\
-c+\frac{\alpha\varepsilon}{u\gamma+w\alpha} & {u+3\leq j\leq n,}
\end{cases}
\end{equation*}
$g_0=f$, $m_{g_{\varepsilon}}=m_f$, such that the remainder function 
\begin{align*}
R(g_{\varepsilon})(\varepsilon)
:=
\Vert g_{\varepsilon}-m_{g_{\varepsilon}} \Vert_p^p
-C_{S_n,p}^p{\rm Var}_{p}^{S_n}(g_{\varepsilon})^p
\end{align*}
satisfies $R(g_{\varepsilon})(0)=0$, $R(g_{\varepsilon})'(0)=0$, $R(g_{\varepsilon})''(0)>0$, contradicting the extremality of $f$.

Specifically, we have that
\begin{align*}
\frac{R(g_{\varepsilon})''(0)}{p(p-1)}
&=
\beta
-u\alpha\Big(\frac{\gamma}{u\gamma+w\alpha}\Big)^2
-w\gamma\Big(\frac{\alpha}{u\gamma+w\alpha}\Big)^2.
\end{align*}
Since 
\[
\Big(u\alpha\Big(\frac{\gamma}{u\gamma+w\alpha}\Big)^2
+w\gamma\Big(\frac{\alpha}{u\gamma+w\alpha}\Big)^2\Big)
\Big(\frac{u}{\alpha}+\frac{w}{\gamma}\Big)
=\Big(\frac{u\gamma}{u\gamma+w\alpha}+\frac{w\alpha}{u\gamma+w\alpha}\Big)^2=1,
\]
it follows that 
\begin{align*}
\frac{R(g_{\varepsilon})''(0)}{p(p-1)}
=\beta-\Big(\frac{u}{\alpha}+\frac{w}{\gamma}\Big)^{-1}
\geq
\beta-\Big(\frac{u}{\alpha}+\frac{1}{\gamma}\Big)^{-1}.
\end{align*}
From \eqref{LAE} and the assumption \eqref{ASSMP2P}, we know that 
\begin{align}
K_f
&=C_{S_n,p}^p a^{p-1}-(a-m_f)^{p-1} \label{KFA}\\
&=(m_f+b)^{p-1}-C_{S_n,p}^p b^{p-1} \label{KFB}\\
&=(m_f+c)^{p-1}-C_{S_n,p}^p c^{p-1} \label{KFC}.
\end{align}
Thus, by \eqref{KFA} and the definition of $\alpha$, we have that
\[
(a-m_f)\alpha=C_{S_n,p}^p a^{p-2}(a-m_f)-(a-m_f)^{p-1}<K_f,
\]
implying
\begin{align*}
\frac{R(g_{\varepsilon})''(0)}{p(p-1)}
&>
\beta-\Big(\frac{u(a-m_f)}{K_f}+\frac{1}{\gamma}\Big)^{-1}\\
&=
\beta-\Big(\frac{nm_f+b+wc-um_f}{K_f}+\frac{1}{\gamma}\Big)^{-1}\\
&>
\beta-\Big(\frac{b+c}{K_f}+\frac{1}{\gamma}\Big)^{-1}.
\end{align*}

We show that $\beta-\Big(\frac{b+c}{K_f}+\frac{1}{\gamma}\Big)^{-1}>0$, which is equivalent to  
\[
(b+c)\beta\gamma+K_f(\beta-\gamma)>0.
\]
By the definition of $\beta$ and $\gamma$, and by \eqref{KFB}, \eqref{KFC}, we have
that 
\[
b\beta=K_f-m_f(m_f+b)^{p-2},\quad c\gamma=m_f(m_f+c)^{p-2}-K_f. 
\]   
Hence
\begin{align*}
bc\Big((b+c)\beta\gamma&+K_f(\beta-\gamma)\Big)
=
(b+c)(b\beta)(c\gamma)+K_f\big(c(b\beta)-b(c\gamma)\big)\\
&=
(b+c)\Big(K_f-m_f(m_f+b)^{p-2}\Big)\Big(m_f(m_f+c)^{p-2}-K_f\Big)\\
&\qquad+K_f\left(c\Big(K_f-m_f(m_f+b)^{p-2}\Big)-b\Big(m_f(m_f+c)^{p-2}-K_f\Big)\right)\\
&=
m_f\Big(K_f\Big( b(m_f+b)^{p-2}+c(m_f+c)^{p-2}\Big)-\\
&\qquad\qquad\qquad
    m_f(b+c)(m_f+b)^{p-2}(m_f+c)^{p-2}\Big).
\end{align*}

From \eqref{KFB} and \eqref{KFC}, it follows that
\[
c^{p-1}(m_f+b)^{p-1}-c^{p-1}K_f=b^{p-1}(m_f+c)^{p-1}-b^{p-1}K_f=C_{S_n,p}^p(bc)^{p-1}.
\]
We thus obtain that
\[
K_f
=
\frac{c^{p-1}(m_f+b)^{p-1}-b^{p-1}(m_f+c)^{p-1}}{c^{p-1}-b^{p-1}},
\]
and that $(b+c)\beta\gamma+K_f(\beta-\gamma)>0$ is equivalent to 
\[
\frac{c^{p-1}(m_f+b)^{p-1}-b^{p-1}(m_f+c)^{p-1}}{c^{p-1}-b^{p-1}}
>
\frac{m_f(b+c)(m_f+b)^{p-2}(m_f+c)^{p-2}}{b(m_f+b)^{p-2}+c(m_f+c)^{p-2}}.
\]

It now remains to verify that $c>m_f>0$; then we can apply Lemma \ref{TAIN} to complete the proof. Since $\gamma=C_{S_n,p}^p c^{p-2}-(m_f+c)^{p-2}>0$, we have that
\begin{align*}
K_f&=(m_f+c)^{p-1}-C_{S_n,p}^p c^{p-1}\\
&<
(m_f+c)^{p-1}-c(m_f+c)^{p-2}\\
&=
m_f(m_f+c)^{p-2}\\
&<
m_f C_{S_n,p}^p c^{p-2}.
\end{align*}
Combining with \eqref{LEQKF}, it follows that 
\[
C_{S_n,p}^p m_f^{p-1}
<
K_f<m_f C_{S_n,p}^p c^{p-2},
\]
indicating $c>m_f>0$. This completes the proof.
\end{proof}

Combining Lemma \ref{QGEQ1}, and Propositions \ref{EOSN1}, \ref{EOSN2} and \ref{EOSN3}, we obtain the following characterization of the extremizers for $C_{S_n,p}$ when $p\in(1,2)\cup(2,+\infty)$.

\begin{thm}\label{EOSN4}
Let $S_n$ be the star graph with $n\geq4$ vertices $a_1$, $a_2,\ldots, a_n$, where $a_1$ is the center, and let $p\in(1,2)\cup(2,+\infty)$. There exists $f:V_{S_n}\to\mathbb{R}$ with
\begin{align*}
C_{S_n,p}
=\frac{\Vert f-m_f\Vert_p}{{\rm Var}_{p}^{S_n}(f)},\quad
m_f\neq f(a_1),
\end{align*}
and 
\begin{align*}
f\left(V_{S_n}\setminus\{a_1\}\right)=\{a,-b\}
\quad
\text{for some}\ a>f(a_1)>-b,\ a,b\in\mathbb{R}.
\end{align*}
Moreover, let
\[
u:=\vert \{a_j\in V_{S_n}: f(a_j)=a\}\vert_\mu,
\quad
v:=\vert \{a_j\in V_{S_n}: f(a_j)=-b\}\vert_\mu.
\]
If $m_f>f(a_1)$, then $u\leq v$ if $p\in(1,2)$, $u>v$ if $p\in(2,+\infty)$.
\end{thm}  

\begin{proof}
By Lemma \ref{QGEQ1}, Propositions \ref{EOSN1}, \ref{EOSN2}, and \ref{EOSN3}, for every $p\in(1,2)\cup(2,+\infty)$, and $n\geq4$, there exists $f$ such that
\[
\frac{\Vert f-m_f\Vert_p}{{\rm Var}_{p}^{S_n}(f)}=C_{S_n,p}>1
\]
and we can assume that $f$ has the form
\[
\bigl(f(a_1),f(a_2),\ldots,f(a_{n-1}),f(a_n)\bigr)
=
\bigl(0,\underbrace{a,\ldots,a}_{u\text{ times}},\underbrace{-b,\ldots,-b}_{v\text{ times}}\bigr),
\]
where $1\leq u,v\leq n-2$, $a,b>0$, and $ua-vb>0$. We claim that $u\leq v$ when $1<p<2$, whereas $u>v$ when $2<p<+\infty$.

In fact, on the one hand, suppose $u>v$ when $1<p<2$, then we have that
\begin{align*}
\frac{\Vert f-m_f\Vert_p^p}{{\rm Var}_{p}^{S_n}(f)^p}-1
&=
\frac{m_f^p+v((m_f+b)^p-b^p)-u(a^p-(a-m_f)^p)}{ua^p+vb^p}.
\end{align*}
Since $x\mapsto (m_f+x)^p-x^p$ is concave and monotonically increasing on $x\in[0,+\infty)$ for $m_f>0$, $1<p<2$, and
\begin{align*}
a-m_f&=\frac{u(a-m_f)+vb-vb}{u}\\
&=
\frac{v}{u}b+\frac{u-v}{u}\frac{u(a-m_f)-vb}{u-v}\\
&=
\frac{v}{u}b+\frac{u-v}{u}\frac{(v+1)m_f}{u-v}
\end{align*}
it follows that 
\begin{align*}
a^p-(a-m_f)^p
&\geq
\frac{v}{u}\Big((m_f+b)^p-b^p\Big)\\
&\qquad
+\frac{u-v}{u}\left(\Big(m_f+\frac{(v+1)m_f}{u-v}\Big)^p-\Big(\frac{(v+1)m_f}{u-v}\Big)^p\right)\\
&>
\frac{v}{u}\Big((m_f+b)^p-b^p\Big)
+\frac{(u-v)m_f^p}{u}.
\end{align*}
Thus
\begin{align*}
\frac{\Vert f-m_f\Vert_p^p}{{\rm Var}_{p}^{S_n}(f)^p}-1
&<
\frac{m_f^p+v((m_f+b)^p-b^p)-v((m_f+b)^p-b^p)-(u-v)m_f^p}{ua^p+vb^p}\\
&=
\frac{(1-(u-v))m_f^p}{ua^p+vb^p}\leq0,
\end{align*}
which contradicts to $C_{S_n,p}>1$. Hence $1\leq u\leq v\leq n-2$. 

On the other hand, suppose $u\leq v$ when $2<p<+\infty$, by Lagrange's mean value theorem, there exist $\xi\in(a-m_f,a)$ and $\eta\in(b,m_f+b)$ such that
\begin{align*}
&m_f^p+v((m_f+b)^p-b^p)-u(a^p-(a-m_f)^p)\\
&\qquad\qquad=
m_f^p+pm_f(v\eta^{p-1}-u\xi^{p-1})\\
&\qquad\qquad<
m_f^p+pm_f(v(m_f+b)^{p-1}-u(a-m_f)^{p-1}).
\end{align*}
Since 
\begin{align*}
u(a-m_f-(m_f+b))
&=
ua-vb-(u+v)m_f+(v-u)(m_f+b)\\
&=
nm_f-(u+v)m_f+(v-u)(m_f+b)\\
&=m_f+(v-u)(m_f+b)>0,
\end{align*}
it follows that
\begin{align*}
&m_f^p+v((m_f+b)^p-b^p)-u(a^p-(a-m_f)^p)\\
&\qquad\qquad<
m_f^p+pm\Big(v(m_f+b)^{p-1}-u(a-m_f)(m_f+b)^{p-2}\Big)\\
&\qquad\qquad=
m_f^p+pm_f(m_f+b)^{p-2}(v(m_f+b)-u(a-m_f))\\
&\qquad\qquad=
m_f^p-pm_f^2(m_f+b)^{p-2}\\
&\qquad\qquad<
(1-p)m_f^p<0,
\end{align*}
which implies that 
\[
\frac{\Vert f-m_f\Vert_p^p}{{\rm Var}_{p}^{S_n}(f)^p}<1,
\]
contradicting to $C_{S_n,p}>1$. Hence $1\leq v< u\leq n-2$. This completes the proof.
\end{proof}

\begin{thm}\label{SPCS}
	Let $S_n$ be the star graph with $n\geq3$ vertices $a_1$, $a_2,\ldots, a_n$, where $a_1$ is the center. Let $f:V_{S_n}\to\mathbb{R}$. We have that

    (i) If $0<p<1$ and $n\geq3$, then 	
\[
	\sum_{j=1}^{n}\vert f(a_j)-m_f\vert^p
	\leq
	\frac{n-1+(n-1)^p}{n^p}{\rm Var}_{p}^{S_n}(f)^p.
\]
The inequality is sharp and, for nonconstant $f$, the equality is attained if and only if there exists a unique $k\in\{2,\ldots,n\}$ such that $f(a_k)\neq f(a_1)$ and $f(a_j)=f(a_1)$ for all $j\neq k$.

	(ii) If $p=1$ and $n\geq3$, then 
	\[
	\sum_{j=1}^{n}\vert f(a_j)-m_f\vert
	\leq
	\Big(2-\frac{2}{n}\Big)
	{\rm Var}_{1}^{S_n}(f).
	\]
	The inequality is sharp and the equality is attained for Dirac delta functions supported at a leaf.
	
(iii) If $1<p<2$ and $n\geq3$, then 
	\begin{align}\label{SPCS12}
	\sum_{j=1}^{n}\vert f(a_j)-m_f\vert^p
	\leq
	M_{n,p}^{(1)}
	{\rm Var}_{p}^{S_n}(f)^p,
	\end{align}
where
\[
M_{n,p}^{(1)}
:=
\max_{\substack{
u\in\mathbb Z,1\leq u\leq\lfloor\frac{n-1}{2}\rfloor\\
0<x<1
}}
\frac{
(n-1-u+x)^p+u^{p-1}(n-1-u+x^p)
}{
\bigl(n-1-u+(u+1)x\bigr)^p
+(n-1-u)u^{p-1}(1-x)^p
}.
\]
The inequality is sharp and the equality is attained for any $f:V_{S_n}\to \mathbb{R}$ assuming values
\[
\bigl(f(a_1),f(a_2),\ldots,f(a_{n-1}),f(a_n)\bigr)
=
\bigl(a,\underbrace{a+b,\ldots,a+b}_{u\text{ times}},\underbrace{a-c,\ldots,a-c}_{n-u-1\text{ times}}\bigr),
\]
where $u\in\mathbb Z$, $1\leq u\leq\lfloor\frac{n-1}{2}\rfloor$, $a\in\mathbb{R}$, $b,c>0$, $ub>(n-u-1)c$. In addition, with
\[
x:=\frac{ub-(n-u-1)c}{ub+(u+1)c},
\]
the parameters must satisfy
\[
\frac{(n-1-u+x)^p+u^{p-1}(n-1-u+x^p)}
{\bigl(n-1-u+(u+1)x\bigr)^p+(n-1-u)u^{p-1}(1-x)^p}
=M_{n,p}^{(1)}.
\]

(iv) If $p=2$, and $n\geq3$, then 
	\[
	\sum_{j=1}^{n}\vert f(a_j)-m_f\vert^p
	\leq
	{\rm Var}_{p}^{S_n}(f)^p.
	\]
	The inequality is sharp and the equality is attained if and only if $f(a_1)=m_f$.

(v) If $2<p<+\infty$ and $n\geq4$, then 
	\[
	\sum_{j=1}^{n}\vert f(a_j)-m_f\vert^p
	\leq
	M_{n,p}^{(2)}
	{\rm Var}_{p}^{S_n}(f)^p,
	\]
where
\[
M_{n,p}^{(2)}
:=
\max_{\substack{
u\in\mathbb Z,\lfloor\frac{n-1}{2}\rfloor+1\leq u\leq n-2\\
0<x<1
}}
\frac{
(n-1-u+x)^p+u^{p-1}(n-1-u+x^p)
}{
\bigl(n-1-u+(u+1)x\bigr)^p
+(n-1-u)u^{p-1}(1-x)^p
}.
\]
The inequality is sharp and the equality is attained for any $f:V_{S_n}\to \mathbb{R}$ assuming values
\[
\bigl(f(a_1),f(a_2),\ldots,f(a_{n-1}),f(a_n)\bigr)
=
\bigl(a,\underbrace{a+b,\ldots,a+b}_{u\text{ times}},\underbrace{a-c,\ldots,a-c}_{n-u-1\text{ times}}\bigr),
\]
where $u\in\mathbb Z$, $\lfloor\frac{n-1}{2}\rfloor+1\leq u\leq n-2$, $a\in\mathbb{R}$, $b,c>0$, $ub>(n-u-1)c$. In addition, with
\[
x:=\frac{ub-(n-u-1)c}{ub+(u+1)c},
\]
the parameters must satisfy
\[
\frac{(n-1-u+x)^p+u^{p-1}(n-1-u+x^p)}
{\bigl(n-1-u+(u+1)x\bigr)^p+(n-1-u)u^{p-1}(1-x)^p}
=M_{n,p}^{(2)}.
\]

If $n=3$, then
	\[
	\sum_{j=1}^{n}\vert f(a_j)-m_f\vert^p
	\leq
	{\rm Var}_{p}^{S_n}(f)^p.
	\]
	The inequality is sharp and the equality is attained for any function $f:V_{S_n}\to \mathbb{R}$ assuming values
	\[
	\Big(f(a_1),f(a_2),f(a_3)\Big)=\Big(\frac{a+b}{2},a,b\Big)
	\]
	for distinct $a,b\in\mathbb{R}$.

\end{thm}

\begin{proof}
(i) If $0<p<1$ and $n\geq3$, then for any $f:V_{S_n}\to\mathbb{R}$, we have that
\begin{align*}
\sum_{j=1}^{n}\vert f(a_j)-m_f\vert^p
&=
\sum_{j=1}^{n}\Big\vert f(a_j)-f(a_1)-(m_f-f(a_1))\Big\vert^p\\
&=
\sum_{j=1}^{n}\Big\vert \frac{n-1}{n}(f(a_j)-f(a_1))-\frac{1}{n}\sum_{i\neq j}(f(a_i)-f(a_1))\Big\vert^p\\
&\leq
\frac{(n-1)^p}{n^p}\sum_{j=1}^{n}\vert f(a_j)-f(a_1)\vert^p
+\frac{1}{n^p}\sum_{j=1}^{n}\sum_{i\neq j}\vert f(a_i)-f(a_1)\vert^p\\
&=
\frac{n-1+(n-1)^p}{n^p}{\rm Var}_{p}^{S_n}(f)^p,
\end{align*}
and, for nonconstant $f$, the equality is attained if and only if there exists a unique $k\in\{2,\ldots,n\}$ such that $f(a_k)\neq f(a_1)$ and $f(a_j)=f(a_1)$ for all $j\neq k$.

(ii) If $p=1$ and $n\geq3$. The choice $\widetilde{f}(a_2)=1$ and $\widetilde{f}(a_i)=0$ for $i\ne2$ gives that
	\[C_{S_n,1}\geq \frac{\Vert \widetilde{f}-m_{\widetilde{f}}\Vert_{1}}{{\rm Var}_{1}^{S_n}(\widetilde{f})}=2-\frac{2}{n}.\] 
	It remains to prove the upper bound $C_{S_n,1}\leq 2-2/n$.
	Since the quotient in \eqref{DBC} is invariant under $f\mapsto C(f-f(a_1))$ ($C\ne0$), we may normalize $f$ so that
	\[
	f(a_1)=0,\qquad \sum_{j=2}^n f(a_j)\ge0,
	\]
	with not all $f(a_j)$ zero, and with the values ordered as
	\[
	f(a_2)\ge f(a_3)\ge\cdots\ge f(a_n).
	\]

	If $\sum_{j=2}^{n}f(a_j)=0$, then $m_f=0$, and
	\[
	\frac{\Vert f-m_{f}\Vert_{1}}{{\rm Var}_{1}^{S_n}(f)}=1.
	\]

	If $\sum_{j=2}^{n}f(a_j)>0$, then
	\begin{align*}
		\frac{\Vert f-m_{f}\Vert_{1}}{{\rm Var}_{1}(f)}
		&=\frac{\sum_{i=1}^n \left\vert f(a_i)
			-\frac{1}{n}\sum_{j=2}^{n}f(a_j)\right\vert}{\sum_{i=2}^n\vert f(a_i)\vert}\\
		&\leq\frac{1}{n}
		+\frac{\sum_{i=2}^{n}\left\vert nf(a_i)-\sum_{j=2}^{n}f(a_j)\right\vert}{n\sum_{i=2}^{n}f(a_i)}\\
		&=
		\frac{1}{n}
		+\sum_{i=2}^{n}\left\vert\frac{f(a_i)}{\sum_{j=2}^{n}f(a_j)}
		-\frac{1}{n}\right\vert.
	\end{align*}
	The equality above is attained if and only if $f(a_2),\cdots,f(a_n)$ are all nonnegative.

	Let  $x_i:=\frac{f(a_i)}{\sum_{j=2}^{n}f(a_j)}$. If $x_i\leq 1/n$ for all $2\leq i\leq n$, then 
	\[
	\frac{\sum_{i=2}^{n}f(a_i)}{\sum_{j=2}^{n}f(a_j)}\leq\frac{n-1}{n},
	\]
	which contradicts. We now consider the following two cases.
	
	\textbf{Case A.} Assume $x_i>1/n$ for all $2\leq i\leq n$.
	Then 
	\[
	\frac{1}{n}
	+\sum_{i=2}^{n}\left\vert\frac{f(a_i)}{\sum_{j=2}^{n}f(a_j)}-\frac{1}{n}\right\vert=\frac{2}{n}.
	\]

	\textbf{Case B.} Assume that there exists $k\in\mathbb{Z}^{+}$, $2\leq k\leq n-1$ such that $x_i\geq1/n$ whenever $2\leq i\leq k$ and $x_i<1/n$ whenever $k+1\leq i\leq n$. Then it follows that
	\begin{align*}
\frac{\Vert f-m_f\Vert_1}{{\rm Var}_1^{S_n}(f)}
&=\frac{\frac1n+\sum_{i=2}^k(x_i-\frac1n)+\sum_{i=k+1}^n(\frac1n-x_i)}
{\sum_{i=2}^n|x_i|}\\
&=\frac{2\sum_{i=2}^k x_i-\frac{2(k-1)}n}{\sum_{i=2}^n|x_i|}\\
&\leq\frac{\sum_{i=2}^n|x_i|+1-\frac{2(k-1)}n}{\sum_{i=2}^n|x_i|}\\
&\leq1+\frac{1-\frac2n}{\sum_{i=2}^n|x_i|}\leq2-\frac2n,
\end{align*}

	Combining the three estimates above yields
	\[
	2-\frac{2}{n}\ge \max\Big\{1,\frac{2}{n}\Big\}\quad (n\ge3).
	\]
	Therefore, the sharp constant is
	\[
	C_{S_n,1}=2-\frac{2}{n}.
	\]

(iii) If $1<p<2$ and $n\geq4$, then it follows from Theorem \ref{EOSN4} that there exists $f:V_{S_n}\to\mathbb{R}$ such that
\begin{align*}
C_{S_n,p}
=\frac{\Vert f-m_f\Vert_p}{{\rm Var}_{p}^{S_n}(f)},\quad
m_f\neq f(a_1),
\end{align*} 
and $f$ is given by
\begin{equation*}
		f(a_j):=
		\begin{cases}
			0 & {j=1,}\\
			t & {2\leq j\leq u+1,}\\
			-1 & {u+2\leq j\leq n.}
		\end{cases}
\end{equation*}
Here $u$ and $t$ are to be determined and satisfy $1\leq u\leq\lfloor\frac{n-1}{2}\rfloor$ and $t>\frac{n-u-1}{u}$. Hence we can write
\begin{align*}
C_{S_n,p}
=
\max_{\substack{
u\in\mathbb Z,\;1\leq u\leq\lfloor\frac{n-1}{2}\rfloor\\
t>\frac{n-u-1}{u}
}}
\Bigg(\frac{m_f^p+u(t-m_f)^p+(n-u-1)(m_f+1)^p}{ut^p+n-u-1}\Bigg)^{\frac{1}{p}},
\end{align*} 
where $m_f=\frac{ut-n+u+1}{n}>0$.

Set $x:=\frac{m_f}{m_f+1}\in(0,1)$. Then $m_f=\frac{x}{1-x}$, $t=\frac{(u+1)x+n-u-1}{u(1-x)}$, $t-m_f=\frac{x+n-u-1}{u(1-x)}$. Substituting these into the above expression, we obtain 
\begin{align*}
C_{S_n,p}
&=
\max_{\substack{
u\in\mathbb Z,\;1\leq u\leq\lfloor\frac{n-1}{2}\rfloor\\
0<x<1
}}
\Bigg(\frac{\Big(\frac{x}{1-x}\Big)^p+u\Big(\frac{n-u-1+x}{u(1-x)}\Big)^p+\frac{n-u-1}{(1-x)^p}}
{u\Big(\frac{(u+1)x+n-u-1}{u(1-x)}\Big)^p+n-u-1}\Bigg)^{\frac{1}{p}}\\
&=
\max_{\substack{
u\in\mathbb Z,\;1\leq u\leq\lfloor\frac{n-1}{2}\rfloor\\
0<x<1
}}
\Bigg(\frac{
(n-1-u+x)^p+u^{p-1}(n-1-u+x^p)
}{
\bigl(n-1-u+(u+1)x\bigr)^p
+(n-1-u)u^{p-1}(1-x)^p
}\Bigg)^{\frac{1}{p}}\\
&=
(M_{n,p}^{(1)})^{\frac{1}{p}}.
\end{align*}

For $n=3$, from Proposition \ref{EOSN1}, we know that 
\begin{align*}
C_{S_3,p}^p
&=\max_{t\in[0,1]}\frac{1}{3^p}\frac{(1+t)^p+(2-t)^p+\vert 1-2t\vert^p}{t^p+(1-t)^p}\\
&=\max_{t\in[0,1/2]}\frac{1}{3^p}\frac{(1+t)^p+(2-t)^p+(1-2t)^p}{t^p+(1-t)^p}\\
&=
\max
\Big\{
1,\frac{2^p+2}{3^p},
\max_{t\in(0,1/2)}\frac{1}{3^p}\frac{(1+t)^p+(2-t)^p+(1-2t)^p}{t^p+(1-t)^p}
\Big\}.
\end{align*}
By Bernoulli's inequality, it follows that
\begin{align*}
&(1+t)^p+(2-t)^p+(1-2t)^p-(2^p+2)(t^p+(1-t)^p)\\
&\qquad\geq
1+pt+2^p-p2^{p-1}t+1-2p t-(2^p+2)(t^p+(1-t)^p)\\
&\qquad>
2^p+2-p(1+2^{p-1})t-(2^p+2)(t^p+1-t)\\
&\qquad=
(2-p)(1+2^{p-1})t-(2^p+2)t^p.
\end{align*}
Thus, it suffices to take $t>0$ sufficiently small such that 
\[
\frac{(1+t)^p+(2-t)^p+(1-2t)^p}{t^p+(1-t)^p}>2^p+2.
\]
Consequently, taking $x=(1-2t)/(1+t)$, we have
\begin{align*}
C_{S_3,p}^p
&=
\max
\Big\{1,\max_{t\in(0,1/2)}\frac{1}{3^p}\frac{(t+1)^p+(2-t)^p+(1-2t)^p}{t^p+(1-t)^p}\Big\}\\
&=\max\Big\{1, \max_{0<x<1}\frac{(1+x)^p+(1+x^p)}{(1+2x)^p+(1-x)^p}\Big\}.
\end{align*}

Moreover, by Bernoulli's inequality, it follows that
\[
(1+x)^p\geq 1+px,\quad (1-x)^p\leq(1-x)(1-(p-1)x)=1-px+(p-1)x^2,
\]
and
\[
(1+2x)^p\leq(1+2x)(1+2(p-1)x)=1+2px+4(p-1)x^2
\]
for $0<x<1$ and $1<p<2$. Thus there exists $\delta\in(0,1)$ sufficiently small such that $x^p>5(p-1)x^2$ for any $x\in(0,\delta)$, and
\[
\frac{(1+x)^p+(1+x^p)}{(1+2x)^p+(1-x)^p}
\geq
\frac{2+px+x^p}{2+px+5(p-1)x^2}>1,\quad x\in(0,\delta),
\]
which implies that
\begin{align*}
C_{S_3,p}^p
&=
\max\Big\{1, \max_{0<x<1}\frac{(1+x)^p+(1+x^p)}{(1+2x)^p+(1-x)^p}\Big\}\\
&=
\max_{0<x<1}\frac{(1+x)^p+(1+x^p)}{(1+2x)^p+(1-x)^p}\\
&=M_{3,p}^{(1)}.
\end{align*}

(iv) If $p=2$ and $n\ge 3$, then for any $f:V_{S_n}\to\mathbb{R}$, by setting $g(a_j)=f(a_j)-f(a_1)$ for $1\le j\le n$, we have
\begin{align*}
&\frac{\Vert f-m_{f}\Vert_{2}^2}{{\rm Var}_{2}^{S_n}(f)^2}
=
\frac{\Vert g-m_{g}\Vert_{2}^2}{{\rm Var}_{2}^{S_n}(g)^2}\\
&\quad=
\frac{\Big(\displaystyle\sum_{j=2}^{n}g(a_j)\Big)^2}{n^2\displaystyle\sum_{j=2}^{n} g(a_j)^2}
+
\frac{\displaystyle\sum_{j=2}^{n}\Big( ng(a_j)-\displaystyle\sum_{i=2}^{n}g(a_i)\Big)^2}{n^2\displaystyle\sum_{j=2}^{n}g(a_j)^2}\\
&\quad=
\frac{\Big(\displaystyle\sum_{j=2}^{n}g(a_j)\Big)^2+\displaystyle\sum_{j=2}^{n}\Big(n^2g(a_j)^2
+
\Big(\displaystyle\sum_{i=2}^{n}g(a_i)\Big)^2-2ng(a_j)\displaystyle\sum_{i=2}^{n}g(a_i)\Big)}
{n^2\displaystyle\sum_{j=2}^{n}g(a_j)^2}\\
&\quad=
\frac{n\Big(\displaystyle\sum_{j=2}^{n}g(a_j)\Big)^2
+n^2\displaystyle\sum_{j=2}^{n}g(a_j)^2
-2n\Big(\displaystyle\sum_{j=2}^{n}g(a_j)\Big)^2}
{n^2\displaystyle\sum_{j=2}^{n}g(a_j)^2}\\
&\quad=
1-\frac{1}{n}\frac{\Big(\displaystyle\sum_{j=2}^{n}g(a_j)\Big)^2}
{\displaystyle\sum_{j=2}^{n} g(a_j)^2}\leq1.
\end{align*}
The equality is attained if and only if $\sum_{j=2}^{n}g(a_j)=0$, which is equivalent to $f(a_1)=m_f$.

By taking $\widetilde{f}(a_2)=1$, $\widetilde{f}(a_1)=\widetilde{f}(a_3)=\cdots=\widetilde{f}(a_{n-1})=0$, $\widetilde{f}(a_n)=-1$ then we have $C_{S_n,2}\geq1$. 

(v) If $2<p<+\infty$ and $n\geq4$, then by Theorem \ref{EOSN4} and an argument analogous to that in (iii), we obtain 
\begin{align*}
C_{S_n,p}^p
=
\max_{\substack{
u\in\mathbb Z,\lfloor\frac{n-1}{2}\rfloor+1\leq u\leq n-2\\
0<x<1
}}
\frac{
(n-1-u+x)^p+u^{p-1}(n-1-u+x^p)
}{
\bigl(n-1-u+(u+1)x\bigr)^p
+(n-1-u)u^{p-1}(1-x)^p
}.
\end{align*}

For $n=3$, since 
\begin{align*}
C_{S_3,p}^p
&=\max_{t\in[0,1/2]}\frac{1}{3^p}\frac{(1+t)^p+(2-t)^p+(1-2t)^p}{t^p+(1-t)^p},
\end{align*}
the function $x\mapsto x^{p/2}$ is convex on $[0,+\infty)$ for $2<p<+\infty$,
and the sequence $\{(3t)^2, (3-3t)^2,0\}$ majorizes $\{(2-t)^2+(12t^2-12t+3), (1+t)^2, (1-2t)^2\}$ for $t\in[0,1/2]$, we obtain from Karamata's inequality:
\begin{align*}
(3t)^p+(3-3t)^p
&\geq \Big((2-t)^2+(12t^2-12t+3)\Big)^{\frac{p}{2}}+(1+t)^p+(1-2t)^p\\
&\geq (2-t)^p+(1+t)^p+(1-2t)^p,
\end{align*}
implying $C_{S_3,p}^p\leq1$. By taking $\widetilde{f}(a_1)=1/2$, $\widetilde{f}(a_2)=1$, and $\widetilde{f}(a_3)=0$, we have that 
\[
C_{S_{3},p}^p
\geq \frac{\Vert \widetilde{f}-m_{\widetilde{f}}\Vert_p^p}{{\rm Var}_{p}^{S_n}(\widetilde{f})^p}=
1.
\] 
Therefore $C_{S_3,p}^p=1$. This completes the proof.
\end{proof}

For $p\in(1,2)\cup(2,+\infty)$, we obtained the expression for $C_{S_n,p}$ in which the number of variables is reduced from $n$ to $2$.
Moreover, we found that for relatively small $n$, the number of variables can be reduced to $1$ (since for $n=4$, $u$ can take only one value for $p\in(1,2)$ and only one value for $p\in(2,+\infty)$).
In this case, we have that 
\begin{equation*}
C_{S_4,p}=
\begin{cases}
\displaystyle\max_{x\in(0,1)}\Big(\frac{(2+x)^p+2+x^p}{(2+2x)^p+2(1-x)^p}\Big)^\frac{1}{p} & {1<p<2,}\\[0.4cm]
\displaystyle\max_{x\in(0,1)}\Big(\frac{(1+x)^p+2^{p-1}(1+x^p)}{(1+3x)^p+2^{p-1}(1-x)^p}\Big)^\frac{1}{p} & {2<p<+\infty.}
\end{cases}
\end{equation*}
For general $n$, we propose the following conjecture:
\begin{con}
Let $S_n$ be the star graph with $n\geq 5$ vertices. Then there exists a constant $N(p)$ depending only on $p$ such that for every $4\leq n\leq N(p)$, we have that
\begin{equation*}
C_{S_n,p}=
\begin{cases}
\displaystyle\max_{x\in(0,1)}\Big(\frac{(n-2+x)^p+n-2+x^p}{(n-2+2x)^p+(n-2)(1-x)^p}\Big)^\frac{1}{p} & {1<p<2,}\\[0.4cm]
\displaystyle\max_{x\in(0,1)}\Big(\frac{(1+x)^p+(n-2)^{p-1}(1+x^p)}{(1+(n-1)x)^p+(n-2)^{p-1}(1-x)^p}\Big)^\frac{1}{p} & {2<p<+\infty.}
\end{cases}
\end{equation*}  
\end{con}

For $p\in(1,2)\cup(2,+\infty)$, $C_{S_n,p}$ can be expressed as the maximum of a bivariate function over its domain. However, we can prove that $C_{S_n,p}$ is uniformly bounded for all $p\in[1,+\infty)$. We have the following result.

\begin{prop}
Let $S_n$ be the star graph with $n\geq3$ vertices $a_1$, $a_2,\ldots, a_n$ where $a_1$ is the center, then we have 
\[
1\leq C_{S_n,p}\leq 2^{\vert2/p-1\vert},\quad
\sup_{n\in\mathbb{Z}^+}C_{K_n,p}<+\infty.
\]
for all $1\leq p<+\infty$.
\end{prop}

\begin{proof}
For every $p\in(1,2)\cup(2,+\infty)$, and $n\geq3$, assume $f$ is the function such that
\[
\frac{\Vert f-m_f\Vert_p}{{\rm Var}_{p}^{S_n}(f)}=C_{S_n,p},
\]
and without loss of generality, we can assume that $f(a_2)\geq\cdots\geq f(a_n)$.
For $\varepsilon>0$ sufficiently small, choose the function $f_{\varepsilon}^{(2,n)}:V_{S_n}\to\mathbb{R}$ and consider the remainder function
\begin{align*}
R(f_{\varepsilon,S_n}^{(2,n)})(\varepsilon)
&:=
\Vert f_{\varepsilon,S_n}^{(2,n)}-m_{f_{\varepsilon,S_n}^{(2,n)}} \Vert_p^p
-C_{S_n,p}^p{\rm Var}_{p}^{S_n}(f_{\varepsilon,S_n}^{(2,n)})^p. 
\end{align*}
For sufficiently small $|\varepsilon|$, define the same remainder by the displayed formula. It is differentiable, is nonpositive and vanishes at $0$. Hence $R(f_{\varepsilon,S_n}^{(2,n)})'(0)=0$ and 
	\begin{align*}		C_{S_{n},p}^p=\frac{\left(f(a_{2})-m_{f}\right)^{p-1}+\left(m_{f}-f(a_{n})\right)^{p-1}}
		{\left(f(a_{2})-f(a_{1})\right)^{p-1}+\left(f(a_{1})-f(a_{n})\right)^{p-1}}.
	\end{align*}

If $p=1$, then $C_{S_n,1}=2-2/n<2$; if $p=2$, then $C_{S_n,2}=1$. If $1<p<2$, then the function $z\mapsto z^{p-1}$ is concave on $[0,+\infty)$. By Jensen's inequality and $(a+b)^{p-1}\leq a^{p-1}+b^{p-1}$ for any $a,b\geq0$, we have 
	\begin{align*}
\frac{\left(f(a_{2})-m_{f}\right)^{p-1}+\left(m_{f}-f(a_{n})\right)^{p-1}}
{\left(f(a_{2})-f(a_{1})\right)^{p-1}+\left(f(a_{1})-f(a_{n})\right)^{p-1}}
<2^{2-p}.
	\end{align*}

If $2<p<+\infty$, then 
\begin{align*}
C_{S_{n},p}^p
<\frac{\left(f(a_2)-f(a_n)\right)^{p-1}}
{\left(f(a_{2})-f(a_{1})\right)^{p-1}+\left(f(a_{1})-f(a_{n})\right)^{p-1}}
\leq2^{p-2}.
\end{align*}

The lower bound $C_{S_n,p}\geq1$ follows by taking the center value $0$ and two leaf values $1$ and $-1$, with all remaining values $0$. For complete graphs, by \cite[Theorem 1]{GRM3} and Theorem \ref{KNSC}, it is easy to see that $\sup_{n\in\mathbb{Z}^+}C_{K_n,p}<+\infty$ holds for all $1\leq p<+\infty$.
This completes the proof.
\end{proof}

\section{Sharp Poincar\'{e}-Wirtinger inequalities on path graphs}\label{s4}

We devote this section to studying the sharp constants in the Poincar\'{e}-Wirtinger inequalities on path graphs. We first characterize the extremizers on the path graphs.

\begin{thm}\label{EOPN}
Let $P_n$ be the path graph with $n\geq3$ vertices $a_1$, $a_2,\ldots, a_n$ such that $a_1\sim a_2\sim\cdots\sim a_n$. For $0<p<+\infty$, there exists $f:V_{P_n}\to\mathbb{R}$ with
\begin{align}\label{EXTREPN}
C_{P_n,p}
=\frac{\Vert f-m_f\Vert_p}{{\rm Var}_{p}^{P_n}(f)},
\end{align}
such that $f$ is monotonic with respect to the natural order of vertices. Moreover, if $1<p<+\infty$, then $f$ given above is strictly monotonic.
\end{thm}

\begin{proof}
For $0<p<+\infty$, from Lemma \ref{EX}, there exists $f$ such that \eqref{EXTREPN} holds.

If $f$ is not monotone, then there is some $j\in\{2,\dots,n-2\}$ satisfying
\begin{align}\label{eqp662}
f(a_{j-1})\ge f(a_{j+1})\ge f(a_j)\ge f(a_{j+2}).
\end{align}
Taking the function
\begin{equation*}
g(a_i)=
\begin{cases}
f(a_{j+1}) & {i=j,}\\
f(a_j) & {i=j+1,}\\
f(a_i) & {\text{otherwise}},
\end{cases}
\end{equation*}
we have $\Vert g-m_g\Vert_p=\Vert f-m_f\Vert_p$. Using \eqref{eqp662}, it follows that 
\begin{align*}
&\vert g(a_{j-1})-g(a_{j})\vert^p
+\vert g(a_{j})-g(a_{j+1})\vert^p
+\vert g(a_{j+1})-g(a_{j+2})\vert^p\\
&=(f(a_{j-1})-f(a_{j+1}))^p
+(f(a_{j+1})-f(a_{j}))^p
+(f(a_{j})-f(a_{j+2}))^p\\
&\leq(f(a_{j-1})-f(a_{j}))^p
+(f(a_{j+1})-f(a_{j}))^p
+(f(a_{j+1})-f(a_{j+2}))^p.
\end{align*}
Consequently, we have 
\[
C_{P_n,p}\geq \frac{\Vert g-m_g\Vert_p^p}{{\rm Var}_{p}^{P_n}(g)^p}
\geq\frac{\Vert f-m_f\Vert_p^p}{{\rm Var}_{p}^{P_n}(f)^p}=C_{P_n,p},
\]
which implies that $g$ is also an extremal function and is monotone. Without loss of generality, we may assume that $f$ is monotonically decreasing.

Now fix $p\in(1,+\infty)$, and suppose there exists an index $k\in\{1,\dots,n-1\}$ such that $f(a_k)=f(a_{k+1})$.
Considering the three cases $k=1$, $k=n-1$, and $2\le k\le n-2$ separately yields a contradiction in each. Therefore, the assumption is false, and $f$ must be strictly monotone.

To this end, for $\varepsilon>0$ sufficiently small, define the function $f_{\varepsilon,P_n}^{(s,t)}$ on $P_n$ by
\begin{equation*}
	f_{\varepsilon,P_n}^{(s,t)}(a_j):=
	\begin{cases}
		f(a_s)+\varepsilon & {j=s,}\\
		f(a_t)-\varepsilon & {j=t,}\\
		f(a_j) & {j\neq s,t.}
	\end{cases}
\end{equation*} 
Then $m_{f_{\varepsilon,P_n}^{(s,t)}}=m_f$.

\textbf{Case 1.} Assume $k=1$. Then $f(a_1)=f(a_{2})>f(a_3)$ and $f(a_2)>m_f$. For $0<\varepsilon<\min\Big\{f(a_2)-m_f,\frac{f(a_2)-f(a_3)}{1+2^{\frac{p}{p-1}}}\Big\}$, we have that
\begin{align*}
	\frac{\Big\Vert f_{\varepsilon,P_n}^{(1,2)}-m_{f_{\varepsilon,P_n}^{(1,2)}}\Big\Vert_p^p}
	{{\rm Var}_p^{P_n}(f_{\varepsilon,P_n}^{(1,2)})^p}
	=
	\frac{
		\Vert f-m_f\Vert_p^p
		+N_1(\varepsilon)
	}
	{
		{\rm Var}_{p}^{P_n}(f)^p
		+D_1(\varepsilon)
	},
\end{align*} 
where
\begin{align*}
	N_1(\varepsilon)&:=
	(f(a_1)-m_f+\varepsilon)^p
	+( f(a_{2})-m_f-\varepsilon)^p\\
	&\quad -( f(a_1)-m_f)^p
	-(f(a_{2})-m_f)^p,
\end{align*}
and
\begin{align*}
	D_1(\varepsilon)
	&:=(2\varepsilon)^p
	+(f(a_{2})-f(a_{3})-\varepsilon)^p
	-(f(a_{2})-f(a_{3}))^p.
\end{align*}

Since 
\[
\frac{\partial}{\partial\varepsilon}D_1(\varepsilon)
=p\Big(2\cdot(2\varepsilon)^{p-1}-(f(a_2)-f(a_3)-\varepsilon)^{p-1}\Big)<0,
\]
it follows that 
\[
D_1(\varepsilon)<\lim_{\varepsilon\to0^+}D_1(\varepsilon)=0.
\]

Since $x\mapsto x^p$ is strictly convex for $p>1$, applying Karamata's inequality gives $N_1(\varepsilon)>0$. Thus
\begin{align*}
	\frac{\Big\Vert f_{\varepsilon,P_n}^{(1,2)}-m_{f_{\varepsilon,P_n}^{(1,2)}}\Big\Vert_p^p}
	{{\rm Var}_p^{P_n}(f_{\varepsilon,P_n}^{(1,2)})^p}
	>
	\frac{\Vert f-m_f\Vert_p^p}{{\rm Var}_{p}^{P_n}(f)^p},
\end{align*}
which contradicts to that $f$ is the extremizer.

\textbf{Case 2.} Assume $2\le k\le n-2$. Then $f(a_{k-1})>f(a_k)=f(a_{k+1})>f(a_{k+2})$. Define
\[
\eta:=\min\Big\{(f(a_{k-1})-f(a_k))/3,\ (f(a_{k+1})-f(a_{k+2}))/3,\ \beta \Big\},
\]
where
\[
\beta:=
\begin{cases}
	|f(a_k)-m_f|, & f(a_k)\neq m_f,\\
	+\infty, & f(a_k)=m_f.
\end{cases}
\]
Then choose $\varepsilon$ such that $0<\varepsilon<\eta$.

It follows that
\begin{align*}
\frac{\Big\Vert f_{\varepsilon,P_n}^{(k,k+1)}-m_{f_{\varepsilon,P_n}^{(k,k+1)}}\Big\Vert_p^p}
{{\rm Var}_p^{P_n}(f_{\varepsilon,P_n}^{(k,k+1)})^p}
=
\frac{
\Vert f-m_f\Vert_p^p
+N_2(\varepsilon)
}
{
{\rm Var}_{p}^{P_n}(f)^p
+D_2(\varepsilon)
},
\end{align*}
where
\begin{align*}
N_2(\varepsilon)&:=
\vert f(a_k)-m_f+\varepsilon\vert^p
+\vert f(a_{k+1})-m_f-\varepsilon\vert^p\\
&\quad -\vert f(a_k)-m_f\vert^p
-\vert f(a_{k+1})-m_f\vert^p,
\end{align*}
and
\begin{align*}
D_2(\varepsilon)&:=
( f(a_{k-1})-f(a_k)-\varepsilon)^p
+(2\varepsilon)^p
+( f(a_{k+1})-f(a_{k+2})-\varepsilon)^p\\
&\quad-(f(a_{k-1})-f(a_k))^p
-( f(a_{k+1})-f(a_{k+2}))^p.
\end{align*}

For $1<p<+\infty$, Karamata's inequality gives $N_2(\varepsilon)>0$ for the three cases $m_f>f(a_k)$, $m_f=f(a_k)$, $m_f<f(a_k)$, and $D_2(\varepsilon)<0$ for the two cases $f(a_{k-1})-f(a_k)\le f(a_{k+1})-f(a_{k+2})$ and $f(a_{k-1})-f(a_k)\ge f(a_{k+1})-f(a_{k+2})$. Hence we have
\begin{align*}
\frac{\Big\Vert f_{\varepsilon,P_n}^{(k,k+1)}-m_{f_{\varepsilon,P_n}^{(k,k+1)}}\Big\Vert_p^p}
{{\rm Var}_p^{P_n}(f_{\varepsilon,P_n}^{(k,k+1)})^p}
>
\frac{\Vert f-m_f\Vert_p^p}{{\rm Var}_{p}^{P_n}(f)^p},
\end{align*}
which contradicts to that $f$ is the extremizer.

\textbf{Case 3.}  Assume $k=n-1$. By the same argument as in Case 1, we obtain a contradiction to the extremality of $f$; hence we omit the details.

This proves the theorem.
\end{proof}

\begin{thm}\label{SPCP}
Let $P_n$ be the path graph with $n\geq4$ vertices $a_1$, $a_2,\ldots, a_n$ such that $a_1\sim a_2\sim\cdots\sim a_n$. Let $f:V_{P_n}\to\mathbb{R}$. The following estimates hold:

(i) If $0<p\leq1$, then 
\[
\sum_{j=1}^{n}\vert f(a_j)-m_f\vert^p
\leq
\frac{\lfloor\frac{n}{2}\rfloor\lceil\frac{n}{2}\rceil^p
+\lfloor\frac{n}{2}\rfloor^p\lceil\frac{n}{2}\rceil}{n^p}
{\rm Var}_{p}^{P_n}(f)^p,
\]
The inequality is sharp and the equality is attained if and only if
$f(a_1)=f(a_2)=\cdots=f(a_r)$, and $f(a_{r+1})=\cdots=f(a_n)$ hold for nonconstant $f$, where $r\in\{\lfloor n/2\rfloor,\lceil n/2\rceil\}$.

(ii) If $p=2$, then 
\[
\sum_{j=1}^{n}\vert f(a_j)-m_f\vert^2
\leq
\frac{1}{4\sin^2(\frac{\pi}{2n})}
{\rm Var}_{2}^{P_n}(f)^2.
\]
The inequality is sharp and, for nonconstant $f$, the equality is attained if and only if $f(a_j)=k\cos\Big(\frac{(2j-1)\pi}{2n}\Big)+l$ for all $1\leq j\leq n$ with $k,l\in\mathbb{R}$, $k\neq0$.

(iii) If $1<p<+\infty$, then
\begin{align*}
\max\bigg\{\frac{n}{\pi},\frac{(n-1)}{2(p+1)^{\frac{1}{p}}}\bigg\}
\leq C_{P_n,p}\leq U_{n,p}
\end{align*}
when $n$ is odd, and
\begin{align*}
\max\bigg\{\frac{n}{2^{\frac{p-1}{p}}\pi},\frac{n-1}{2(p+1)^{\frac{1}{p}}}\bigg\}
\leq C_{P_n,p}\leq U_{n,p}
\end{align*}
when $n$ is even, where
\begin{equation}\label{PATHINTERP}
U_{n,p}:=
\begin{cases}
\displaystyle
\left(\frac{2\lfloor n/2\rfloor\lceil n/2\rceil}{n}\right)^{\frac{2}{p}-1}
\left(\frac{1}{2\sin(\frac{\pi}{2n})}\right)^{2-\frac{2}{p}}
& {1<p<2,}\\[0.4cm]
\displaystyle
\left(\frac{1}{2\sin(\frac{\pi}{2n})}\right)^{\frac{2}{p}}
\left(\frac{n-1}{2}\right)^{1-\frac{2}{p}}
& {2\leq p<+\infty.}
\end{cases}
\end{equation}
Moreover, for each fixed $1\leq p<+\infty$, $C_{P_n,p}=O(n)$ as $n\to+\infty$. For each fixed $0<p<1$, $C_{P_n,p}$ has order $n^{1/p}$.
\end{thm}

\begin{proof}
(i) First assume $0<p<1$. By Theorem \ref{EOPN}, it suffices to consider a monotone extremizer. Without loss of generality, let
\[
f(a_1)\ge f(a_2)\ge\cdots\ge f(a_n).
\]
Then there exists $k\in\{1,\dots,n-1\}$ such that $f(a_k)\ge m_f\ge f(a_{k+1})$.

Using the identity
\[
\sum_{j=1}^{k} f(a_j) - k m_f = (n-k)m_f - \sum_{j=k+1}^{n} f(a_j),
\]
and applying Jensen's inequality to the concave function $x\mapsto x^p$, we get
\begin{align*}
&\sum_{j=1}^{n}\vert f(a_j)-m_f\vert^p
=\sum_{j=1}^{k}(f(a_j)-m_f)^p+\sum_{j=k+1}^{n}(m_f-f(a_j))^p\\
&\leq 
k^{1-p}\Big(\sum_{j=1}^{k}f(a_j)-km_f\Big)^p
+(n-k)^{1-p}\Big((n-k)m_f-\sum_{j=k+1}^{n}f(a_j)\Big)^p\\
&\leq (k^{1-p}+(n-k)^{1-p})\Big((n-k)m_f - \sum_{j=k+1}^{n} f(a_j)\Big)^p \\
&\leq \frac{k(n-k)^p+k^p(n-k)}{n^p}(f(a_1)-f(a_n))^p.
\end{align*}
The equality is attained if and only if $f(a_1)=f(a_2)=\cdots=f(a_k)$, and $f(a_{k+1})=f(a_{k+2})=\cdots=f(a_n)$. 

It follows from
\begin{align*}
\sum_{j=1}^{n-1}(f(a_j)-f(a_{j+1}))^p
\geq\Big(\sum_{j=1}^{n-1}(f(a_j)-f(a_{j+1}))\Big)^p
=(f(a_1)-f(a_n))^p
\end{align*}
that 
\[
\frac{\Vert f-m_f\Vert_p^p}{{\rm Var}_{p}^{P_n}(f)^p}
\leq
\frac{k(n-k)^p+k^p(n-k)}{n^p}.
\]

For $0<p<1$, define $\psi(x):=x(n-x)^p+x^p(n-x)$ on $[1,n-1]$. Observe that $\psi(x)=\psi(n-x)$ for $x\in [1,n-1]$. For $x\in[1,n/2]$, it follows from $\psi'(x)\geq 0$ that $\psi$ attains its maximum at $x=n/2$. Consequently, the upper bound over integer $k$ is maximal when $k=\lfloor\frac{n}{2}\rfloor$ or $k=\lceil\frac{n}{2}\rceil$. A function with two distinct constant levels separated at either such index attains this bound. Conversely, strict concavity and strict subadditivity force equality to have two levels and precisely one transition along the original path, giving the equality cases in (i).

For $p=1$, set $d_j:=f(a_j)-f(a_{j+1})$ for $1\leq j\leq n-1$. A direct calculation gives
\begin{equation}\label{PATHTMATRIX}
f(a_i)-m_f=\sum_{j=1}^{n-1}\mathcal{T}_{ij}d_j,\qquad
\mathcal{T}_{ij}:=
\begin{cases}
\frac{n-j}{n} & {i\leq j,}\\[0.1cm]
-\frac{j}{n} & {i>j.}
\end{cases}
\end{equation}
Hence
\begin{align*}
\Vert f-m_f\Vert_1
&\leq \sum_{j=1}^{n-1}\left(\sum_{i=1}^n|\mathcal{T}_{ij}|\right)|d_j|
=\sum_{j=1}^{n-1}\frac{2j(n-j)}{n}|d_j|\\
&\leq\frac{2\lfloor n/2\rfloor\lceil n/2\rceil}{n}{\rm Var}_1^{P_n}(f).
\end{align*}
The two-level functions in (i) attain equality. Conversely, equality for a nonconstant function requires $d_j=0$ except at $j\in\{\lfloor n/2\rfloor,\lceil n/2\rceil\}$. For odd $n$, both central differences cannot be nonzero: equality in the first row of the triangle inequality requires them to have the same sign, whereas equality in the central row requires opposite signs. Thus precisely one central difference is nonzero, giving the stated equality cases for $p=1$.

(ii) For $p=2$, let $f$ be an extremizer for \eqref{EXTREPN}, normalized by $m_f=0$ and $\Vert f\Vert_2=1$. The gradients of the constraints $\sum_jx_j=0$ and $\sum_jx_j^2=1$ are linearly independent. Consider the Lagrangian
\[
\mathfrak{L}_2({\bf x},\lambda,\mu)
:=\sum_{j=1}^{n-1}(x_j-x_{j+1})^2-\lambda\Big(\sum_{j=1}^{n}x_j^2-1 \Big)-2\mu\sum_{j=1}^n x_j,
\]
where ${\bf x}:=(x_1,\ldots,x_n)^T$. Then we have that
\[
\frac{\partial}{\partial x_j}\mathfrak{L}_2({\bf x},\lambda,\mu)\Bigg\vert_{x_j=f(a_j)}=0
\]
holds for all $1\leq j\leq n$. Summing these equations and using $\sum_{j=1}^n f(a_j)=0$ gives $n\mu=0$. Thus $\mu=0$, and the stationary equations are equivalent to 
\begin{equation}\label{LAEQ}
\begin{cases}
x_1-x_2=\lambda x_1, & \\
2x_j-x_{j-1}-x_{j+1}=\lambda x_j & \text{for}\ {2\leq j\leq n-1,}\\
x_n-x_{n-1}=\lambda x_n. & 
\end{cases}
\end{equation}
It then follows that
\begin{align*}
&\sum_{j=1}^{n-1}\bigl(f(a_j)-f(a_{j+1})\bigr)^2 
= \sum_{j=1}^{n-1}\Bigl(f(a_j)^2-2f(a_j)f(a_{j+1})+f(a_{j+1})^2\Bigr) \\
&= f(a_1)^2+2\sum_{j=2}^{n-1}f(a_j)^2+f(a_n)^2-2\sum_{j=1}^{n-1}f(a_j)f(a_{j+1}) \\
&= f(a_1)^2+2\sum_{j=2}^{n-1}f(a_j)^2+f(a_n)^2 -\sum_{j=2}^{n-1} f(a_j)f(a_{j-1})-f(a_{n-1})f(a_n) \\
&\quad - f(a_1)f(a_2)-\sum_{j=2}^{n-1}f(a_j)f(a_{j+1}) \\
&= f(a_1)\bigl(f(a_1)-f(a_2)\bigr) 
+\sum_{j=2}^{n-1}f(a_j)\Bigl(2f(a_j)-f(a_{j-1})-f(a_{j+1})\Bigr) \\
&\quad +f(a_n)\bigl(f(a_n)-f(a_{n-1})\bigr) \\
&= \lambda\sum_{j=1}^{n}f(a_j)^2.
\end{align*}
Hence, we have that 
\[
\frac{\Vert f-m_f\Vert_2^2}{{\rm Var}_2^{P_n}(f)^2}=\frac{1}{\lambda}\leq\frac{1}{\lambda_{\min}},
\]
where $\lambda_{\min}$ is the minimal positive value of $\lambda$ such that \eqref{LAEQ} holds. 

In order to find $\lambda$, we observe that \eqref{LAEQ} is equivalent to 
\[
A_n(\lambda)\mathbf{x}:=
\begin{pmatrix}
1-\lambda & -1        & 0  & \cdots & 0      & 0 \\
-1        & 2-\lambda & -1 & \cdots & 0      & 0 \\
0         & -1        & 2-\lambda   & \cdots & 0 & 0 \\
\vdots    & \vdots    & \vdots      &        & \vdots    & \vdots \\
0         & 0         & 0           & \cdots & 2-\lambda & -1  \\
0         & 0         & 0           & \cdots & -1        & 1-\lambda
\end{pmatrix}
_{ n\times n}
\mathbf{x}=\mathbf{0},
\]
and since $\mathbf{x}$ is nonzero, it follows that
\[
{\rm det }(A_n(\lambda))=0.
\]

Denote by $P_k(\lambda)$ ($1\le k\le n-1$) the $k$-th leading principal minor of the matrix $A_n(\lambda)$. With $P_0(\lambda)=1$ and $P_1(\lambda)=1-\lambda$, the recurrence relation for $2\le k\le n-1$ reads
\[
P_k(\lambda)-(2-\lambda)P_{k-1}(\lambda)+P_{k-2}(\lambda)=0
\]
and
\[
P_{n}(\lambda)=(1-\lambda)P_{n-1}(\lambda)-P_{n-2}(\lambda).
\]

Multiplying \eqref{LAEQ} by $x_j$ and summing shows that $\lambda=\sum_{j=1}^{n-1}(x_j-x_{j+1})^2>0$. Also,
\[
\sum_{j=1}^{n-1}(x_j-x_{j+1})^2\leq 2\sum_{j=1}^{n-1}(x_j^2+x_{j+1}^2)
=4\sum_{j=1}^n x_j^2-2(x_1^2+x_n^2)\leq4\sum_{j=1}^n x_j^2.
\]
Equality between the first and last expressions would require $x_1=x_n=0$ and $x_j=-x_{j+1}$ for all $j$, forcing the entire vector to vanish. Thus the combined bound is strict. Thus $0<\lambda<4$.

The general term $P_k(\lambda)$ ($2\le k\le n-1$) is determined by the characteristic equation
\begin{equation*}
r^2-(2-\lambda)r+1=0.
\end{equation*}
Its roots are conjugate; denote them by $r_1$ and $r_2$, so that
\begin{equation}\label{CE1}
r_1+r_2=2-\lambda,\qquad r_1r_2=1.
\end{equation}
Thus we may write $r_1=e^{-i\theta}$, $r_2=e^{i\theta}$ for some $\theta\in(0,\pi)$. The closed-form expression for $P_k(\lambda)$ is therefore
\[
P_k(\lambda)=A r_1^k+ B r_2^k,
\]
where the constants $A$ and $B$ are determined by the initial conditions $P_0(\lambda)=1$ and $P_1(\lambda)=1-\lambda$. Explicitly,
\begin{equation}\label{eqn668143}
\begin{cases}
A+B=1, & \\
r_1 A+r_2 B=1-\lambda, & 
\end{cases}
\end{equation}
where $1-\lambda=r_1+r_2-1$.

By \eqref{eqn668143}, we have 
\[A=\frac{1-r_1}{r_2-r_1}\quad \text{and}\quad B=\frac{r_2-1}{r_2-r_1}.\] 
Hence, 
\begin{align*}
P_{k}(\lambda)
=\frac{r_2^{k+1}-r_1^{k+1}-(r_2^k-r_1^k)}{r_2-r_1},
\end{align*}
and thus
\begin{align*}
P_{n}(\lambda)
= \frac{\lambda(r_1^n-r_2^n)}{r_2-r_1}
= \frac{-2{\rm i}\lambda\sin(n\theta)}{r_2-r_1}.
\end{align*}

It follows from $P_{n}(\lambda)=0$ that $\theta=\frac{k\pi}{n}$, $k=1,\ldots,n-1$. Meanwhile, from \eqref{CE1}, we have that 
\[
\lambda=2-2\cos\Big(\frac{k\pi}{n}\Big),
\]
and consequently, 
\[
\lambda_{\min}
=\min\limits_{k=1,\ldots,n-1}\Big\{2-2\cos\Big(\frac{k\pi}{n}\Big)\Big\}
=4\sin^2\Big(\frac{\pi}{2n}\Big).
\]
We have thus obtained the upper bound
\[C_{P_{n},2}^2=
\frac{\Vert f-m_f\Vert_2^2}{{\rm Var}_2^{P_n}(f)^2}=\frac{1}{\lambda}
\leq
\frac{1}{4\sin^2\big(\frac{\pi}{2n}\big)}.
\]

To see its sharpness, choose the test function  $f(a_j)=\sqrt{\frac{2}{n}}\cos\Big(\frac{(2j-1)\pi}{2n}\Big)$
for $1\leq j\leq n$, then $\sum_{j=1}^{n}f(a_j)=0$, $\sum_{j=1}^{n}f(a_j)^2=1$, and we have
\[
C_{P_{n},2}^2
\geq
\frac{1}{4\sin^2\big(\frac{\pi}{2n}\big)}.
\]
This establishes the sharpness of the constant. The first equation in \eqref{LAEQ}, together with the recurrence, determines every coordinate from $x_1$. Consequently, the eigenspace corresponding to $\lambda_{\min}$ is one-dimensional and is spanned by the cosine vector above. Translation and dilation invariance give exactly the equality cases in (ii).

(iii) For $1<p<+\infty$, we consider the lower bound of $C_{P_n,p}$ in the two separate cases where $n$ is even or odd.

If $n=2k+1$ with $k\ge2$, then choosing
\[
f_1(a_j)=k+1-j,\qquad g(a_j)=\cos\Big(\frac{(2j-1)\pi}{2n}\Big),\qquad 1\le j\le n,
\]
gives
\begin{align}\label{eqn66814}
C_{P_{n},p}^p
\geq
\frac{\Vert f_1-m_{f_1}\Vert_p^p}{{\rm Var}_p^{P_n}(f_1)^p}
=\frac{\sum_{j=1}^{k}j^p}{k}
\geq
\frac{1}{k}\int_{0}^{k}x^p\ dx
=\frac{(n-1)^p}{(p+1)2^p},
\end{align}
and
\begin{align*}
C_{P_{n},p}^p
&\geq
\frac{\Vert g-m_g\Vert_p^p}{{\rm Var}_p^{P_n}(g)^p}
=
\frac{2\sum_{j=1}^{k}\cos^p\Big(\frac{(2j-1)\pi}{2n}\Big)}
{\sum_{j=1}^{n-1}\Big(\cos\Big(\frac{(2j-1)\pi}{2n}\Big)-\cos\Big(\frac{(2j+1)\pi}{2n}\Big)\Big)^p}\\
&=
\frac{1}{\Big(2\sin\frac{\pi}{2n}\Big)^p}
\geq
\frac{n^p}{\pi^p}.
\end{align*}

If $n=2k$ with $k\geq2$, then taking 
\[f_2(a_j)=2k+1-2j,\qquad g(a_j)=\cos\Big(\frac{(2j-1)\pi}{2n}\Big), \qquad 1\leq j\leq n,\] gives 
\begin{align}\label{PNLBE}
C_{P_{n},p}^p\geq \frac{\Vert f_2-m_{f_2}\Vert_p^p}{{\rm Var}_p^{P_n}(f_2)^p} 
=\frac{\sum_{j=1}^{k}(2j-1)^p}{(2k-1)2^{p-1}}
\geq\frac{(n-1)^p}{(p+1)2^p},
\end{align}
and
\begin{align*}
C_{P_{n},p}^p
&\geq
\frac{\Vert g-m_g\Vert_p^p}{{\rm Var}_p^{P_n}(g)^p}
=
\frac{2\sum_{j=1}^{k}\sin^p\Big(\frac{(2j-1)\pi}{2n}\Big)}
{\Big(2\sin\frac{\pi}{2n}\Big)^p\sum_{j=1}^{n-1}\sin^p\Big(\frac{j\pi}{n}\Big)}.
\end{align*}
By the concavity of $\sin$ and Jensen's inequality, and then using $(a+b)^p\ge a^p+b^p$ for nonnegative $a,b$ and $p>1$, we get
\[
\sin^p\left(\frac{(2j-1)\pi}{2n}\right)
\ge \frac{1}{2^p}\left(\sin^p\left(\frac{j\pi}{n}\right)+\sin^p\left(\frac{(j-1)\pi}{n}\right)\right).
\]
Hence
\begin{align*}
C_{P_{n},p}^p
&\geq
\frac{2\sum_{j=1}^{k}\Big(\sin^p\Big(\frac{j\pi}{n}\Big)+\sin^p\Big(\frac{(j-1)\pi}{n}\Big)\Big)}
{2^{p}\Big(2\sin\frac{\pi}{2n}\Big)^p\Big(2\sum_{j=1}^{k-1}\sin^p\Big(\frac{j\pi}{n}\Big)+1\Big)}\\
&=
\frac{2\Big(2\sum_{j=1}^{k-1}\sin^p\Big(\frac{j\pi}{n}\Big)+1\Big)}
{2^{p}\Big(2\sin\frac{\pi}{2n}\Big)^p\Big(2\sum_{j=1}^{k-1}\sin^p\Big(\frac{j\pi}{n}\Big)+1\Big)}\\
&\geq
\frac{n^p}{2^{p-1}\pi^p}.
\end{align*}

To estimate the upper bound of $C_{P_n,p}$ for $1<p<+\infty$, let $\mathcal{T}$ be the linear map defined in \eqref{PATHTMATRIX}. Every vector of differences $d\in\mathbb R^{n-1}$ is realized by a function on $P_n$, so
\[
C_{P_n,p}=\Vert \mathcal{T}\Vert_{\ell^p\to\ell^p}.
\]
Part (i) and part (ii) give
\[
\Vert \mathcal{T} \Vert_{\ell^1\to\ell^1}
=\frac{2\lfloor n/2\rfloor\lceil n/2\rceil}{n},\qquad
\Vert \mathcal{T}\Vert_{\ell^2\to\ell^2}
=\frac{1}{2\sin(\frac{\pi}{2n})}.
\]
For the endpoint $p=+\infty$, the absolute row sums give directly
\begin{align*}
\Vert \mathcal{T}\Vert_{\ell^\infty\to\ell^\infty}
&=\max_{1\leq i\leq n}\sum_{j=1}^{n-1}|T_{ij}|\\
&=\max_{1\leq i\leq n}
\frac{i(i-1)+(n-i)(n-i+1)}{2n}
=\frac{n-1}{2}.
\end{align*}
Applying the Riesz--Thorin interpolation theorem, it yields
\[
C_{P_n,p}\leq U_{n,p},
\]
with $U_{n,p}$ as in \eqref{PATHINTERP}. Finally, the lower bounds above and this upper bound show that $C_{P_n,p}=O(n)$ for each fixed $p>1$. The formula in (i) gives $C_{P_n,1}=O(n)$ and order $n^{1/p}$ for each fixed $0<p<1$.
This completes the proof.
\end{proof}

\begin{rem}\label{PATHSYMMETRY}
For $p=2$, the equality cases in Theorem \ref{SPCP} (ii) satisfy
\[
f(a_j)+f(a_{n-j+1})=2m_f,\qquad 1\leq j\leq n.
\]
Thus the centered extremizers are antisymmetric under reflection. For general $p>1$, $p\neq2$, this property need not hold. For example, let $n=3$ and $p=3/2$, and consider
\[
\big(f(a_1),f(a_2),f(a_3)\big)=(1,0,-1),\qquad
\big(g(a_1),g(a_2),g(a_3)\big)=(1,0,-1/2).
\]
The first function is antisymmetric and satisfies
\[
\frac{\Vert f-m_f\Vert_p}{{\rm Var}_p^{P_3}(f)}=1.
\]
The second function has $m_g=1/6$ and satisfies
\[
\frac{\Vert g-m_g\Vert_{3/2}}{{\rm Var}_{3/2}^{P_3}(g)}
=\left(\frac{5\sqrt5+1+8}{6\sqrt6(1+1/(2\sqrt2))}\right)^{2/3}
\approx1.009603>1.
\]
Every nonconstant function on $P_3$ that is antisymmetric about its arithmetic mean is a translation and dilation of $(1,0,-1)$. Hence no such function is an extremizer at $p=3/2$. We leave an explicit characterization of extremizers for general $p>1$, $p\neq2$, including the determination of when reflection antisymmetry holds, as an open problem.
\end{rem}

\begin{thm}\label{thm4.3}
	Let $K_n$ be the complete graph, $S_n$ the star graph, and $P_n$ the path graph on $n\ge3$ vertices. Then 
	\[
	\lim_{p\to+\infty}C_{K_n,p}=1-\frac{1}{n},\quad \lim_{p\to+\infty}C_{S_n,p}=2-\frac{3}{n},
	\]
	and
	\[
	\lim_{p\to+\infty}C_{P_n,p}=\frac{n-1}{2}
	\]
	when $n\geq4$.
\end{thm}

\begin{proof}
From \cite[Theorem 1]{GRM3} and Theorem \ref{KNSC} (iv), we have that
\begin{align*}
\lim_{p\to+\infty}C_{K_n,p}
&=
\lim_{p\to+\infty}
\Big(\frac{1+(n-1)^{p-1}}{n^p}\Big)^{\frac{1}{p}}\\
&=
\lim_{p\to+\infty}\exp
\Big\{\frac{(p-1)\ln(n-1)+\ln(1+(n-1)^{-(p-1)})}{p}-\ln n\Big\}\\
&=
1-\frac{1}{n}.
\end{align*}

For star graphs, the case $n=3$ is trivial.
For $2<p<+\infty$ and star graph $S_n$ with $n\geq4$ vertices $a_1$, $a_2,\ldots, a_n$, where $a_1$ is the center, from Theorem \ref{EOSN4}, we know that we can assume $f:V_{S_n}\to \mathbb{R}$ has the form 
\begin{align*}
\Big(f(a_1),f(a_2),\ldots,f(a_{n-1}),f(a_n)\Big)
=
\left(
0,
\underbrace{a,\ldots,a}_{u\text{-times}},
\underbrace{-b,\ldots,-b}_{v\text{-times}}
\right)
\end{align*}
with $1\leq v<u\leq n-2$, $u+v=n-1$, $a,b>0$, satisfying
\[
C_{S_n,p}^p
=
\frac{\Vert f-m_f\Vert_p^p}{{\rm Var}_{p}^{S_n}(f)^p}
=
\frac{m_f^p+u(a-m_f)^p+v(m_f+b)^p}{ua^p+vb^p}
>1.
\]
Then on the one hand, from
\[
m_f
=
\frac{ua-vb}{n}
\leq
\frac{u\max\{a,b\}}{n}
\leq
\frac{n-2}{n}\max\{a,b\}
<
\Big(2-\frac{3}{n}\Big)\max\{a,b\},
\]
\[
a-m_f
=
\frac{na-(ua-vb)}{n}\leq\frac{2v+1}{n}\max\{a,b\}\leq\Big(2-\frac{3}{n}\Big)\max\{a,b\},
\]
and
\[
m_f+b
=
\frac{ua-vb+nb}{n}\leq\frac{2u+1}{n}\max\{a,b\}\leq\Big(2-\frac{3}{n}\Big)\max\{a,b\},
\]
it follows that
\begin{align}\label{SNPINFI1}
C_{S_n,p}
=
\Big(\frac{m_f^p+u(a-m_f)^p+v(m_f+b)^p}{ua^p+vb^p}\Big)^{1/p}
<n^{\frac{1}{p}}\Big(2-\frac{3}{n}\Big).
\end{align}
On the other hand, by taking $\widetilde{f}(a_1)=0$, $\widetilde{f}(a_2)=\cdots=\widetilde{f}(a_{n-1})=1$, $\widetilde{f}(a_n)=-1$, we have that
\begin{align}\label{SNPINFI2}
C_{S_n,p}
&\geq
\frac{\Vert \widetilde{f}-m_{\widetilde{f}}\Vert_p}{{\rm Var}_{p}^{S_n}(\widetilde{f})}
=
\Bigg(\frac{\Big(1-\frac{3}{n}\Big)^p+(n-2)\Big(\frac{3}{n}\Big)^p+\Big(2-\frac{3}{n}\Big)^p}{n-1}
\Bigg)^{\frac{1}{p}}\\
&>
\frac{1}{(n-1)^{1/p}}\Big(2-\frac{3}{n}\Big). \nonumber
\end{align}
Hence, combining \eqref{SNPINFI1} and \eqref{SNPINFI2}, we obtain 
\[
\lim_{p\to+\infty}C_{S_n,p}=2-\frac{3}{n}.
\]

For path graphs with $n\geq4$, the upper bound in Theorem \ref{SPCP} (iii) gives
\[
\limsup_{p\to+\infty}C_{P_n,p}
\leq\lim_{p\to+\infty}U_{n,p}=\frac{n-1}{2}.
\]
On the other hand, taking $\widetilde{f}(a_j)=\frac{n+1}{2}-j$ for $1\leq j\leq n$, we have $m_{\widetilde{f}}=0$ and
\[
C_{P_n,p}\geq
\bigg(\frac{\sum_{j=1}^n|\frac{n+1}{2}-j|^p}{n-1}\bigg)^{1/p}.
\]
The right hand side tends to $(n-1)/2$ as $p\to+\infty$. Therefore,
\[
\lim_{p\to+\infty}C_{P_n,p}=\frac{n-1}{2}.
\]
This completes the proof.
\end{proof}

\section{Estimates of the best constant on general finite graphs}\label{s5}

As obtained in Theorem \ref{thm4.3}, the limits are
\[
\lim_{p\to+\infty}C_{K_n,p}=1-\frac{1}{n},\qquad
\lim_{p\to+\infty}C_{S_n,p}=2-\frac{3}{n},\qquad
\lim_{p\to+\infty}C_{P_n,p}=\frac{n-1}{2}
\]
for the complete graph $K_n$, the star graph $S_n$, and the path graph $P_n$, respectively, under the corresponding restrictions on $n$. In this section, we establish sharp upper and lower bounds for the best constant $C_{G_n,p}$ in the Poincar\'{e}-Wirtinger inequality on arbitrary finite connected graphs $G_n$ with $n\ge2$ vertices. For $1<p<+\infty$, the extremal graphs for the sharp bounds are the complete graph and the path graph, which also characterize the equality cases.

The case $p=+\infty$ is of a different nature. The constant $C_{G_n,\infty}$ is related to the vertex deviation (see \cite{Z}) and to the remoteness of the graph $G_n$ (see \cite{ACH}). Here the remoteness of $G_n$, i.e., the maximum average distance from a vertex to all vertices, with the distance from a vertex to itself taken as zero, is defined by
\[
\operatorname{rem}(G_n)
:=
\max_{a_i\in V_{G_n}}\frac{1}{n}\sum_{a_j\in V_{G_n}}d_{G_n}(a_i,a_j).
\]

This observation yields the following explicit estimate, which can be regarded as a general graph version of the limits obtained in the previous section for $K_n$, $S_n$, and $P_n$. In particular, for these graphs the upper bound below is attained and recovers the values
\[
1-\frac{1}{n},\qquad 2-\frac{3}{n},\qquad \frac{n-1}{2},
\]
respectively.

\begin{thm}
	Let $G_n=(V_{G_n},E_{G_n})$ be a simple connected graph with $n\ge2$ vertices. Then the best constant in the Poincar\'{e}-Wirtinger inequality at $p=\infty$ is given by
	\[
	C_{G_n,\infty}=\operatorname{rem}(G_n).
	\]
Equivalently, for every $f:V_{G_n}\to\mathbb{R}$,
	\[
	\max_{a_j\in V_{G_n}}|f(a_j)-m_f|
	\le
	\operatorname{rem}(G_n)
	\max_{\{a_i,a_j\}\in E_{G_n}}|f(a_i)-f(a_j)|,
	\]
	and the constant $\operatorname{rem}(G_n)$ is sharp.
\end{thm}

\begin{proof} The case $n=2$ is trivial. For $n\ge3$, let $V_{G_n}:=\{a_1,\dots,a_n\}$ denote the vertex set. Without loss of generality, we may normalize $f$ by scaling so that \[
\max_{\{a_i,a_j\}\in E_{G_n}} |f(a_i)-f(a_j)|=1.
\]
For any $a_j\in V_{G_n}$, since 
\begin{align*}
\vert f(a_j)-m_f\vert
&=\frac{1}{n}\Big\vert n f(a_j)-\sum_{a_i\in V_{G_n}}f(a_i)\Big\vert
\leq\frac{1}{n}\sum_{a_i\in V_{G_n}} d_{G_n}(a_i,a_j)\\
&\leq
\frac{1}{n}\max_{a_i\in V_{G_n}}\sum_{a_j\in V_{G_n}}d_{G_n}(a_i,a_j).
\end{align*}
Since $G_n$ is finite, there exists $a_{i_0}$ such that
\[
\sum_{a_j\in V_{G_n}}d_{G_n}(a_{i_0},a_j)=\max_{a_i\in V_{G_n}}\sum_{a_j\in V_{G_n}}d_{G_n}(a_i,a_j).
\]
By taking $\widetilde{f}(a_j)=d_{G_n}(a_j,a_{i_0})$ for any $a_j\in V_{G_n}$, we have 
\begin{align*}
C_{G_n,+\infty}
&\geq
\frac{\max_{a_j\in V_{G_n}}\vert \widetilde{f}(a_j)-m_{\widetilde{f}}\vert}{\max_{\{a_i,a_j\}\in E_{G_n}}\vert \widetilde{f}(a_i)-\widetilde{f}(a_j)\vert}\\
&\geq\vert m_{\widetilde{f}}\vert\\
&=
\frac{1}{n}\max_{a_i\in V_{G_n}}\sum_{a_j\in V_{G_n}}d_{G_n}(a_i,a_j),
\end{align*}
indicating that the inequality is sharp. This completes the proof.
\end{proof}

\begin{thm}\label{BFG}
Let $G_n$ be a simple connected graph with $n\geq3$ vertices. For $0<p<+\infty$, we have that
\begin{align}\label{kgp}
C_{K_n,p}\leq C_{G_n,p}\leq C_{P_n,p}.
\end{align}
Moreover,

(i) For any $0<p<+\infty$, $G_n=K_n$ if and only if $C_{G_n,p}=C_{K_n,p}$;

(ii) For any $1<p<+\infty$, $G_n$ is isomorphic to $P_n$ if and only if $C_{G_n,p}=C_{P_n,p}$.
\end{thm}

\begin{proof}
For $n=3$, connectedness forces $G_3$ to be isomorphic to either $K_3$ or $P_3$. From \cite[Theorem 1]{GRM3} together with Theorems \ref{KNSC} and \ref{SPCS}, using $P_3$ is isomorphic to $S_3$, we have $C_{K_3,p}<C_{P_3,p}$ for all $0<p<+\infty$. It remains to consider the case $n\ge4$.

We first prove \eqref{kgp}. Let $V_n:=\{a_1,\dots,a_n\}$. Suppose $G_n=(V_n,E_{G_n})$ and $H_n=(V_n,E_{H_n})$ are two arbitrary finite connected graphs such that $E_{G_n}\subseteq E_{H_n}$. Then for any nonconstant function $f:V_n\to\mathbb{R}$, we have that
\begin{align*}
{\rm Var}_{p}^{G_n}(f)
&=\Big(\frac{1}{2}\sum_{a_i,a_j\in V_{G_n},a_i\sim a_j}\vert f(a_i)-f(a_j)\vert^p\Big)^{1/p}\\
&=\Big(\sum_{\{a_i, a_j\}\in E_{G_n}}\vert f(a_i)-f(a_j)\vert^p\Big)^{1/p}\\
&\leq
\Big(\sum_{\{a_i, a_j\}\in E_{H_n}}\vert f(a_i)-f(a_j)\vert^p\Big)^{1/p}\\
&={\rm Var}_{p}^{H_n}(f).
\end{align*}
Since $K_n$ has the maximal edge set among all graphs on $V_n$,
\[
C_{K_n,p}\le C_{G_n,p}.
\]

Let $e=\{a_i,a_j\}\in E_{G_n}$ and $G'_n:=G_n-e$. If $G'_n$ is connected, then
\[
C_{G_n,p}\leq C_{G'_n,p}.
\]
Repeatedly deleting edges from $G'_n$ yields a sequence of graphs. This process terminates when the graph becomes a tree, denoted by $T_n=(V_{T_n},E_{T_n})$.

By \cite[Theorem 4.3]{BM}, $T_n$ has $n-1$ edges. Let $V_{T_n}:=\{a_1,\dots,a_n\}$.  Order the vertices so that $f(a_1)\ge\cdots\ge f(a_n)$. If an edge of $T_n$ does not join consecutive vertices in this ordering, choose such an edge $e=\{a_i,a_j\}$ with $j>i+1$. Removing $e$ yields two components $A$ and $B$, with $a_i\in A$, $a_j\in B$. Thus for some $k$ with $i\le k<j$, we have $a_k$ and $a_{k+1}$ in different components. Adding $e':=\{a_k,a_{k+1}\}$ to $T_n- e$ gives a new graph $T'_n$.

Since $\vert f(a_i)-f(a_j)\vert^p\geq\vert f(a_k)-f(a_{k+1})\vert^p$ for all $0<p<+\infty$, it follows that  
\[
{\rm Var}_{p}^{T_n}(f)\geq{\rm Var}_{p}^{T'_n}(f).
\]   
Every such replacement decreases the number of edges joining nonconsecutive vertices in the ordering. By performing finitely many deletions and additions of edges, we eventually arrive at a graph with
\[
a_1\sim a_2\sim\cdots\sim a_n,
\]
which is precisely a path graph. From the preceding analysis, the inequality 
\[
C_{G_n,p}\leq C_{P_n,p}
\]
is valid for every $0<p<+\infty$, which establishes \eqref{kgp}.

We now turn to the proof of (i). The necessity follows directly from \cite[Theorem 1]{GRM3} and Theorem \ref{KNSC}, while the sufficiency will be proved by contradiction.  Assume that $G_n\neq K_n$ while $C_{G_n,p}=C_{K_n,p}$. We shall derive a contradiction by showing that $G_n\neq K_n$ implies $C_{G_n,p}>C_{K_n,p}$.

Since $G_n\neq K_n$, there exist two non-adjacent vertices $a_i,a_j\in V_{G_n}$ with $i<j$, and there exists $k\in\mathbb N$ such that $i\leq k<j$. For any integers $s,i,j\in\{1,\dots,n\}$, define
\begin{equation*}
\rho_{k}(a_s):=
\begin{cases}
1 & {1\leq s\leq k,}\\
0 & {k+1\leq s\leq n,}
\end{cases}
\end{equation*}
and for $i\ne j$, define
\begin{equation*}
\sigma_{i,j}(a_s):=\delta_i(a_s)-\delta_j(a_s)=
\begin{cases}
1 & {s=i,}\\
-1 & {s=j,}\\
0 & {\text{otherwise}.}
\end{cases}
\end{equation*}
In the case $i=j$, $\sigma_{i,i}\equiv0$, which is trivial and will be excluded from the subsequent discussion. We treat the following cases according to the interval in which $p$ lies.

(a) If $p\in(0,1)\cup(p_n,+\infty)$, where $p_n$ is the number from Lemma \ref{REU}, then it follows that
\begin{align*}
C_{G_n,p}^p
&\geq
\frac{\Vert \delta_{i}-m_{\delta_{i}}\Vert_p^p}{{\rm Var}_{p}^{G_n}(\delta_{i})^p}
=
\frac{(1-\frac{1}{n})^p+(n-1)(\frac{1}{n})^p}{d_{G_n}(a_i)}\\
&\geq
\frac{(n-1)^p+n-1}{n^p(n-2)}
>
\frac{1+(n-1)^{p-1}}{n^p}\\
&=
C_{K_n,p}^p.
\end{align*}

(b) If $p\in[1,2]$ and $1\le k\le n-1$, relabel the vertices so that the nonadjacent vertices $a_i$ and $a_j$ lie on opposite sides of the partition at $k$. Then it follows that
\begin{align*}
C_{G_n,p}^p
&\geq
\frac{\Vert \rho_{k}-m_{\rho_{k}}\Vert_p^p}{{\rm Var}_{p}^{G_n}(\rho_{k})^p}
\geq
\frac{k(1-\frac{k}{n})^p+(n-k)(\frac{k}{n})^p}{\sum\limits_{\substack{1\leq s\leq k\\ k+1\leq t\leq n}}\vert \rho_k(a_s)-\rho_k(a_t)\vert^p-1}\\
&=
\frac{k(n-k)^p+(n-k)k^p}{n^p(k(n-k)-1)}
>
\frac{(n-k)^{p-1}+k^{p-1}}{n^p}.
\end{align*}
A straightforward calculation yields
\[
C_{G_n,p}^p
>
\max_{1\leq k\leq n-1}\frac{(n-k)^{p-1}+k^{p-1}}{n^p}
=
\frac{\lceil\frac{n}{2}\rceil^{p-1}+\lfloor\frac{n}{2}\rfloor^{p-1}}{n^p}
=
C_{K_n,p}^p.
\]

(c) If $p\in(2,p_n]$ and $i\neq j$, then it follows that
\begin{align*}
C_{G_n,p}^p
&\geq
\frac{\Vert \sigma_{i,j}-m_{\sigma_{i,j}}\Vert_p^p}{{\rm Var}_{p}^{G_n}(\sigma_{i,j})^p}\\
&\geq
\frac{\vert\sigma_{i,j}(a_i)\vert^p+\vert\sigma_{i,j}(a_j)\vert^p}
{\sum_{s=1}^{n}\Big(\vert \sigma_{i,j}(a_s)-\sigma_{i,j}(a_i)\vert^p+\vert \sigma_{i,j}(a_s)-\sigma_{i,j}(a_j)\vert^p\Big)-2\cdot2^p}\\
&=
\frac{1}{n-2}
>
\frac{1}{2^{p-1}+n-2}\\
&=
C_{K_n,p}^p.
\end{align*}

The assumption is therefore false, and $G_n=K_n$ follows.

We now proceed to prove (ii).  Let $P_n=(V_n, E_{P_n})$ with $E_{P_n}=\{\{a_i,a_{i+1}\}:1\leq i\leq n-1\}$.
The preceding argument shows that for any finite graph $G_n=(V_n, E_{G_n})$, and any function $f$ on $V_n$,
\begin{align}\label{VPG}
	{\rm Var}_{p}^{G_n}(f)^p\geq{\rm Var}_{p}^{P_n}(f)^p
\end{align}
holds for all $0<p<+\infty$. 
Moreover, if $f(a_i)\neq f(a_j)$ for every pair of adjacent vertices $a_i$ and $a_j$, then \eqref{VPG} is strict, and consequently,
\[
\frac{\Vert f-m_f\Vert_p^p}{{\rm Var}_{p}^{G_n}(f)^p}
<
\frac{\Vert f-m_f\Vert_p^p}{{\rm Var}_{p}^{P_n}(f)^p}
\leq
C_{P_n,p}^p.
\]   
If there exist adjacent vertices $a_i$ and $a_j$ such that $f(a_i)= f(a_j)$, then by Theorem \ref{EOPN}, we have that
\[
\frac{\Vert f-m_f\Vert_p^p}{{\rm Var}_{p}^{G_n}(f)^p}
\leq
\frac{\Vert f-m_f\Vert_p^p}{{\rm Var}_{p}^{P_n}(f)^p}
<
C_{P_n,p}^p,
\]   
for all $1<p<+\infty$. This concludes the proof.
\end{proof}

For $p \in (0,1]$, there actually exist other finite connected graphs whose sharp constant in the Poincar\'{e}-Wirtinger inequality coincides with that of the path graph. We give the following example.
\begin{rem}\label{CTE} 
For $n\geq 4$, let $D_n$ be the dumbbell graph formed by two disjoint complete graphs $K_{\lfloor n/2\rfloor}$ and $K_{\lceil n/2\rceil}$ joined by a single edge. Define the test function $f$ on $V_{D_n}$ by
\[
f(v)=
\begin{cases}
	1 & v\in V_{K_{\lfloor n/2\rfloor}},\\
	0 & v\in V_{K_{\lceil n/2\rceil}}.
\end{cases}
\]	
Then we obtain that 
\[
C_{K_{\lfloor\frac{n}{2}\rfloor}\bigoplus K_{\lceil\frac{n}{2}\rceil},p}
\geq
\Big(\frac{\lfloor\frac{n}{2}\rfloor\lceil\frac{n}{2}\rceil^p
+\lfloor\frac{n}{2}\rfloor^p\lceil\frac{n}{2}\rceil}{n^p}\Big)^{\frac{1}{p}}.
\]
Meanwhile, from \eqref{kgp}, which is from Theorem \ref{BFG}, together with Theorem \ref{SPCP} (i), it follows that
\[
C_{K_{\lfloor\frac{n}{2}\rfloor}\bigoplus K_{\lceil\frac{n}{2}\rceil},p}
\leq
\Big(\frac{\lfloor\frac{n}{2}\rfloor\lceil\frac{n}{2}\rceil^p
+\lfloor\frac{n}{2}\rfloor^p\lceil\frac{n}{2}\rceil}{n^p}\Big)^{\frac{1}{p}}
=
C_{P_n,p},
\]
implying
\[
C_{K_{\lfloor\frac{n}{2}\rfloor}\bigoplus K_{\lceil\frac{n}{2}\rceil},p}
=
C_{P_n,p}.
\]
 
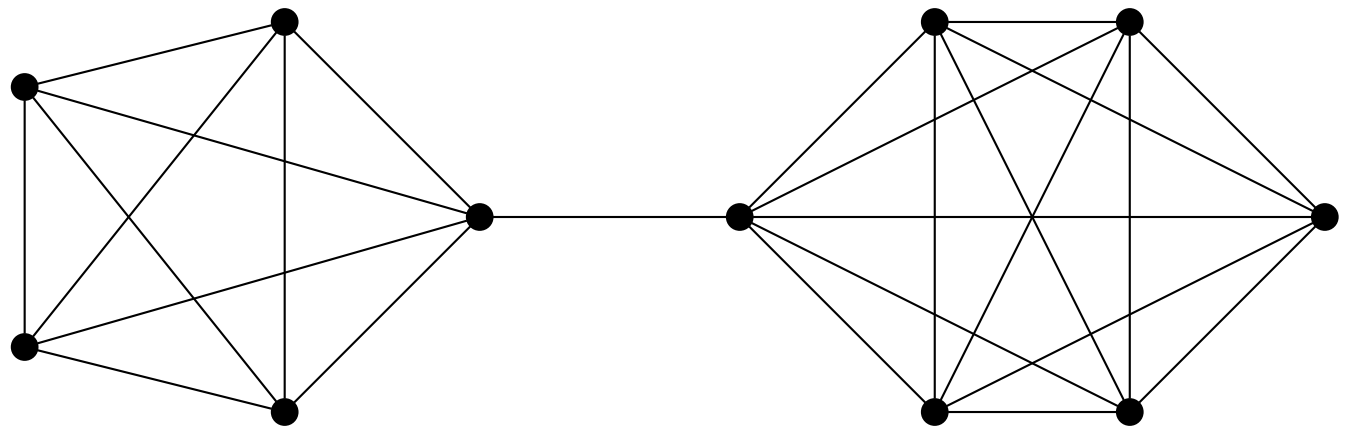
\begin{figure}[!htbp]
  \centering
\tikzstyle{every node}=[scale=0.5]
\begin{tikzpicture}[scale=0.43]
	\begin{pgfonlayer}{nodelayer}
		\node [style=Black] (218) at (-2, 0) {};
		\node [style=Black] (224) at (-5, 3) {};
		\node [style=Black] (225) at (-5, -3) {};
		\node [style=Black] (226) at (2, 0) {};
		\node [style=Black] (227) at (5, 3) {};
		\node [style=Black] (228) at (5, -3) {};
		\node [style=Black] (280) at (8, 3) {};
		\node [style=Black] (281) at (8, -3) {};
		\node [style=Black] (282) at (11, 0) {};
		\node [style=Black] (283) at (-9, 2) {};
		\node [style=Black] (285) at (-9, -2) {};
	\end{pgfonlayer}
	\begin{pgfonlayer}{edgelayer}
		\draw (218) to (226);
		\draw (226) to (227);
		\draw (228) to (226);
		\draw (218) to (224);
		\draw (218) to (225);
		\draw (282) to (280);
		\draw (282) to (281);
		\draw (281) to (227);
		\draw (280) to (228);
		\draw (281) to (228);
		\draw (280) to (227);
		\draw (281) to (280);
		\draw (228) to (227);
		\draw (282) to (226);
		\draw (280) to (226);
		\draw (281) to (226);
		\draw (228) to (282);
		\draw (282) to (227);
		\draw (225) to (285);
		\draw (224) to (283);
		\draw (283) to (285);
		\draw (283) to (225);
		\draw (224) to (225);
		\draw (224) to (285);
		\draw (218) to (283);
		\draw (218) to (285);
	\end{pgfonlayer}
\end{tikzpicture}
  \caption{An example of $K_{\lfloor\frac{n}{2}\rfloor}\bigoplus K_{\lceil\frac{n}{2}\rceil}$ when $n=11$}
\end{figure}
\end{rem}

\section{Cheeger-type estimates on finite graphs}\label{s6}
In this section, we establish an upper bound for the best constant $C_{G_n,p}$ on finite connected graphs in terms of the Cheeger constant $h_{\mu,\nu}(G_n)$. This estimate complements the existing eigenvalue-based bounds and highlights a useful connection between the Poincar\'{e}-Wirtinger constant and the isoperimetric geometry of the graph. The method presented here may also be of independent interest for problems where Cheeger-type information is more readily available.

To illustrate the approach, combining Theorem \ref{BFG} with \cite[Theorem 1.1]{WH} gives the explicit upper bound
\[
C_{G_n,p}\le C_{P_n,p}\le (n(n-1)^{p-1})^{1/p},\qquad 1<p<+\infty,\ n\ge4.
\]

Our main result applies \cite[Theorem 3]{A1} together with Theorem \ref{BFG} to derive a Cheeger-type upper bound for $C_{G_n,p}$ that holds for any finite connected graph $G_n$ with $n\ge4$ vertices. The following elementary lemma will be used in the proof.

\begin{lem}[\cite{A1}]\label{A2LE}
If $a$, $b\geq0$, $a\neq b$, and $p>1$ is conjugate to $p'$, then
\[
\Big(\frac{1}{p}\frac{b^p-a^p}{b-a}\Big)^{p'}\leq\frac{1}{2}(a^p+b^p).
\]
\end{lem}

\begin{prop}\label{p41}
Let $G_n$ be a simple connected graph with $n\geq4$ vertices. For $1<p<+\infty$, the inequality
\[
C_{G_n,p}
\leq
\frac{2^{\frac{1}{p}}p\Delta_{G_n}^{1-\frac{1}{p}}}{h_{\mu,\nu}(G_n)}
\]
holds. 
\end{prop}

\begin{proof}
Let $V_{G_n}=\{a_1,a_2,\ldots, a_n\}$. Let $f:V_{G_n}\to\mathbb{R}$ be a nonzero eigenfunction of the $p$-Laplacian associated with $\lambda_{1,p}(G_n)$. Replacing $f$ by $-f$ if necessary, we may assume that at most $n/2$ vertices have positive values. Define
\[
f_+(a_j):=\max\{f(a_j),0\},\quad f_-(a_j):=\max\{-f(a_j),0\},
\] 
and 
\[
E_t(f_+):=\{a_j\in V_{G_n}:f_+(a_j)>t\}, \quad E_t(f_-):=\{a_j\in V_{G_n}:f_-(a_j)>t\},
\]
for any $t\geq 0$. We may assume without loss of generality that $|E_t(f_+)|_\mu \le n/2$.

For $1<p<+\infty$, applying Lemma \ref{A2LE} and H\"{o}lder's inequality gives:
	\begin{align*}
	&\frac{1}{2}\sum_{a_i,a_j\in V_{G_n}, a_i\sim a_j}\vert f_+(a_i)^p-f_+(a_j)^p\vert\\
	&\leq {\rm Var}_{p}^{G_n}(f_+) \Big(\frac{1}{2}\sum_{a_i,a_j\in V_{G_n},a_i\sim a_j, f_+(a_i)\neq f_+(a_j)}\!\!\!\!
	\Big(\frac{f_+(a_i)^p-f_+(a_j)^p}{f_+(a_i)-f_+(a_j)}\Big)^{p'}\Big)^{\frac{1}{p'}}\\
	&\leq {\rm Var}_{p}^{G_n}(f_+)\Big(\frac{p^{p'}}{4}\sum_{a_i,a_j\in V_{G_n},a_i\sim a_j}\!\!\!\! (f_+(a_i)^p+f_+(a_j)^p)\Big)^{\frac{1}{p'}}\\
	&\leq {\rm Var}_{p}^{G_n}(f_+) \Big(\frac{p^{p'}\Delta_{G_n}}{2}\sum_{a_j\in V_{G_n}}f_+(a_j)^p\Big)^{\frac{1}{p'}}\\
	&=
	p\Big(\frac{\Delta_{G_n}}{2}\Big)^{1-\frac{1}{p}}\Vert f_+\Vert_p^{p-1}{\rm Var}_{p}^{G_n}(f_+).
\end{align*}

Moreover, from \eqref{DOCC}, we have
\begin{align*}
&\frac{1}{2}\sum_{a_i,a_j\in V_{G_n},a_i\sim a_j}\vert f_+(a_i)^p-f_+(a_j)^p\vert\\
&=
\sum_{\substack{
  a_i,a_j\in V_{G_n},a_i\sim a_j\\
  f_+(a_j)> f_+(a_i)
}} 
\int_{f_+(a_i)}^{f_+(a_j)}pt^{p-1}\ dt\\
&=
\int_{0}^{+\infty}pt^{p-1}\sum_{\substack{a_i,a_j\in V_{G_n},
  a_i\sim a_j\\
  f_+(a_j)> f_+(a_i)
}} 
\chi_{[f_+(a_i),f_+(a_j))}(t)\ dt\\
&=
\int_{0}^{+\infty}pt^{p-1}\vert \partial E_t(f_+)\vert_\nu\ dt\\
&\geq
h_{\mu,\nu}(G_n)\int_{0}^{+\infty}pt^{p-1}\vert E_t(f_+)\vert_\mu\ dt\\
&=
h_{\mu,\nu}(G_n)\Vert f_+\Vert_p^p.
\end{align*}
This leads to 
\begin{align}\label{p411}
\frac{\Vert f_+\Vert_p}{{\rm Var}_{p}^{G_n}(f_+)}
\leq
p\Big(\frac{\Delta_{G_n}}{2}\Big)^{1-\frac{1}{p}}\frac{1}{h_{\mu,\nu}(G_n)}.
\end{align}

For any $a_j\in E_0(f_+)$, since $f$ is an eigenfunction associated with $\lambda_{1,p}(G_n)$, it follows that 
\[
\Delta_p f(a_j)
=-\lambda_{1,p}(G_n)f(a_j)^{p-1}
\]
Since $x\mapsto\vert x\vert^{p-2}x$ is monotonically increasing on $(-\infty,0)\cup(0,+\infty)$, it follows that
\begin{align*}
\Delta_p f_+(a_j)
&=
\Bigg(
\sum_{\substack{
  a_i\in E_0(f_+)\\
  a_i\sim a_j\in E_{G_n}
}}
+
\sum_{\substack{
  a_i\notin E_0(f_+)\\
  a_i\sim a_j\in E_{G_n}
}}
\Bigg)\vert f_+(a_i)-f_+(a_j)\vert^{p-2}(f_+(a_i)-f_+(a_j))\\
&=
\sum_{\substack{
  a_i\in E_0(f_+)\\
  a_i\sim a_j\in E_{G_n}
}}\vert f(a_i)-f(a_j)\vert^{p-2}(f(a_i)-f(a_j))
-
\sum_{\substack{
  a_i\notin E_0(f_+)\\
  a_i\sim a_j\in E_{G_n}
}}
f(a_j)^{p-1}\\
&\geq
\Delta_p f(a_j)\\
&=-\lambda_{1,p}(G_n)f(a_j)^{p-1}
\end{align*}
Thus,
\begin{align}\label{p412}
\nonumber 
\lambda_{1,p}(G_n)\Vert f_+\Vert_p^p
&=
\lambda_{1,p}(G_n)\sum_{a_i\in V_{G_n}}f_+(a_i)^{p-1}f_+(a_i)
\\ \nonumber
&=
\lambda_{1,p}(G_n)\Big(\sum_{f_+(a_i)>0}+\sum_{f_+(a_i)=0}\Big)f_+(a_i)^{p-1}f_+(a_i)
\\ \nonumber
&\geq
-\sum_{a_i\in E_0(f_+)}f_+(a_i)\Delta_p f_+(a_i)\\ \nonumber
&=
\frac{1}{2}\sum_{a_i,a_j\in V_{G_n},a_i\sim a_j}\vert f_+(a_i)-f_+(a_j)\vert^p\\
&=
{\rm Var}_{p}^{G_n}(f_+)^p.
\end{align}
Therefore, combining \eqref{REVPC}, \eqref{p411} and \eqref{p412}, we conclude that
\[
C_{G_n,p}
\leq
\frac{2^{\frac{1}{p}}p\Delta_{G_n}^{1-\frac{1}{p}}}{h_{\mu,\nu}(G_n)}.
\]
This completes the proof.
\end{proof}
For $n\geq4$, $h_{\mu,\nu}(P_n)=\frac{1}{\lfloor n/2\rfloor}$ (see \cite[Chapter 2]{M}), and $\Delta_{P_n}=2$. Hence for $1<p<+\infty$, $p\neq2$, from Theorem \ref{BFG} and Proposition \ref{p41}, we have that  
\begin{align*}
C_{G_n,p}
\leq
C_{P_n,p}
\leq
2p\left\lfloor\frac{n}{2}\right\rfloor.
\end{align*} 
However, from Theorem \ref{SPCP} (iii), we have that
\[
C_{G_n,p}\leq C_{P_n,p}
\leq
\left(\frac{2\lfloor n/2\rfloor\lceil n/2\rceil}{n}\right)^{\frac2p-1}
\left(\frac{1}{2\sin(\frac{\pi}{2n})}\right)^{2-\frac2p}
\]
if $1<p<2$, and
\[
C_{G_n,p}\leq C_{P_n,p}
\leq
\left(\frac{1}{2\sin(\frac{\pi}{2n})}\right)^{\frac2p}
\left(\frac{n-1}{2}\right)^{1-\frac2p}
\]
if $2<p<+\infty$. Then for sufficiently large $n\in\mathbb{Z}^+$, we observe that
\[
(n(n-1)^{p-1})^{\frac1p}\approx n,\qquad
2p\left\lfloor\frac n 2\right\rfloor\approx pn,
\]
and
\[
U_{n,p}\approx\frac{n}{2^{\frac2p-1}\pi^{2-\frac2p}},\qquad 1<p<2,
\]
\[
U_{n,p}\approx\frac{n}{2^{1-\frac2p}\pi^{\frac2p}},\qquad 2<p<+\infty.
\]
One can check that
\[
\frac{1}{\pi}<\frac{1}{2^{\frac2p-1}\pi^{2-\frac2p}}<\frac12<\min\{1,p\},\qquad 1<p<2,
\]
and
\[
\frac{1}{\pi}<\frac{1}{2^{1-\frac2p}\pi^{\frac2p}}<\frac12<\min\{1,p\},\qquad 2<p<+\infty,
\]
which further indicates that Theorems \ref{SPCP} and \ref{BFG} provide a better upper bound for the best constant in the Poincar\'{e}-Wirtinger inequality on any finite connected graph than the path-based Cheeger estimate above.


\begin{thebibliography}{99}

\bibitem{AM} N. Alon, V.D. Milman, $\lambda_1$, isoperimetric inequalities for graphs, and superconcentrators, \emph{J. Combin. Theory Ser. B} \textbf{38} (1985), no. 1, 73--88.

\bibitem{A1} S. Amghibech, Eigenvalues of the discrete $p$-Laplacian for graphs, \emph{Ars Combin.} \textbf{67} (2003), 283--302.
    
\bibitem{ACH} M. Aouchiche, G. Caporossi, P. Hansen, Variable neighborhood search for extremal graphs. XX. Automated comparison of graph invariants, \emph{MATCH Commun. Math. Comput. Chem.} \textbf{58} (2007), no. 2, 365--384.
    
\bibitem{BAR} M.T. Barlow, Random walks and heat kernels on graphs, \emph{London Math. Soc. Lecture Note Ser.} vol. 438, Cambridge: Cambridge University Press, 2017.
    
\bibitem{BG} S.G. Bobkov, F. G\"{o}tze, Discrete isoperimetric and Poincar\'{e}-type inequalities, \emph{Probab. Theory Related Fields}  \textbf{114} (1999), no. 2, 245--277.
    
\bibitem{BH} T. B\"uhler, M. Hein, An inverse power method for nonlinear eigenproblems with applications in $1$-spectral clustering and sparse PCA, \emph{NIPS} (2010), 847--855.
    
\bibitem{BM} J.A. Bondy, U.S.R. Murty, Graph theory, \emph{Grad. Texts in Math.} vol. 244, New York: Springer,  2008.
    
\bibitem{BRE} S. Brenner, Discrete Sobolev and Poincar\'{e} inequalities for piecewise polynomial functions, \emph{Electron. Trans. Numer. Anal.} \textbf{18} (2004), 42--48.
    
\bibitem{BPV} T. Bruno, M.M. Peloso, M. Vallarino, Local and non-local Poincar\'{e} inequalities on Lie groups, \emph{Bull. Lond. Math. Soc.} \text{54} (2022), no. 6, 2162--2173.
    
\bibitem{CH} C. Cartensen, F, Hellwig, Constants in discrete Poincar\'{e} and Friedrichs inequalities and discrete quasi-interpolation, \emph{Comput. Methods Appl. Math.} \textbf{18} (2018), no. 3, 433--450.
    
\bibitem{CHA} K.C. Chang, Spectrum of the $1$-Laplacian and Cheeger's constant on graphs, \emph{J. Graph Theory} \textbf{81} (2016), no. 2, 167--207.    

\bibitem{C} J. Cheeger, A lower bound  for the smallest value of the Laplacian, in Problems in Analysis (Paper dedicated to Salomon Bochner, 1969), \emph{Princeton University Press} Princeton, N.J.: 1970, 195--199. 
    
\bibitem{CGY} F.R.K. Chung, A. Grigor'yan, S.-T. Yau, Upper bounds for eigenvalues of the discrete and continuous Laplace operators, \emph{Adv. Math.} \textbf{117} (1996), no. 2, 165--178.

\bibitem{DE} T. Delmotte, Parabolic Harnack inequality and estimates of Markov chains on graphs, \emph{Rev. Mat. Ieroam.} \textbf{15} (1999), no. 1, 181--232.    

    
\bibitem{DO} J. Dodziuk, Difference equations, isoperimetric inequality and transience of certain random walks, \emph{Trans. Amer. Math. Soc.} \textbf{284} (1984), no. 2, 787--794.
    
    
\bibitem{ELP} L.B. Efraim, F. Lust-Piquard, Poincar\'{e} type inequalities on the discrete cube and in the CAR algebra, \emph{Probab. Theory Related Fields} \textbf{65} (2008), no. 3--4, 569--602. 
    
\bibitem{E} L.C. Evans, Partial differential equations, \emph{Grad. Stud. Math.} vol. 19, Providence, RI: American Mathematical Society, 1998.
    
\bibitem{FS} K. Funano, Y. Sakurai, Upper bounds for higher-order Poincar\'{e} constants, \emph{Trans. Amer. Math. Soc.} \textbf{373} (2020), no. 6, 4415--4436.
    
\bibitem{GHJ} H. Ge, B. Hua, W. Jiang, The limit of first eigenfunctions of the $p$-Laplacian on graphs, \emph{Pacific J. Math.} \textbf{312} (2021), no. 1, 103--112.
    
\bibitem{GRM3} C. Gonz\'{a}lez-Riquelme, J. Madrid, Sharp Poincar\'{e}-Wirtinger inequalities on complete graphs, \emph{Ann. Fenn. Math.} \textbf{51} (2026), no. 1, 163--175.
    
\bibitem{HLP} G.H. Hardy, J.E. Littlewood, G. P\'olya, Inequalities, 2nd ed., \emph{Cambridge Mathematical Library}, Cambridge: Cambridge University Press, 1952. 

\bibitem{HLLY} P. Horn, Y. Lin, S. Liu, S.-T. Yau, Volume doubling, Poincar\'{e} inequality and Gaussian heat kernel estimate for non-negatively curved graphs, \emph{J. Reine. Angew. Math.} \textbf{51} (2019), 89--130.
    
\bibitem{HW} B. Hua, L. Wang, Dirichlet $p$-Laplacian eigenvalues and Cheeger constants on symmetric graphs, \emph{Adv. Math.} \textbf{364} (2020), 106997, 34 pp.
    
\bibitem{HMT} D. Hume, J.M. Mackay, R. Tessera, Poincar\'{e} profiles of groups and spaces, \emph{Rev. Mat. Iberoam.} \textbf{36} (2020), no. 6, 1835--1886.
    
\bibitem{ILVHV} P. Ivanisvili, D. Li, R. van Handel, A. Volberg, Improving constant in end-point Poincar\'{e} inequality on Hamming cube, arXiv:1811.05584, 2019.
    
\bibitem{IV} P. Ivanisvili, A. Volberg, Poincar\'{e} inequality $3/2$ on the Hamming cube, \emph{Rev. Mat. Iberoam.} \textbf{36} (2020), no. 1, 79--97.


\bibitem{Ka}
J. Karamata, Sur une in\'{e}galit\'{e} relative aux fonctions convexes, \emph{Publ. Math. Univ. Belgrade} \textbf{1} (1932), 145--148.

    
\bibitem{KZ} S. Keith, X. Zhong, The Poincar\'{e} inequality is an open ended condition, \emph{Ann. of Math. (2)} \textbf{167} (2008), no. 2, 575--599.
    
\bibitem{KM} M. Keller, D. Mugnolo, General Cheeger inequalities for $p$-Laplacians on graphs, \emph{Nonlinear Anal.} \textbf{147} (2016), 80--95.
    
\bibitem{LO} A.H. Le, P. Omnes, Discrete Poincar\'{e} inequalities for arbitrary meshes in the discrete duality finite volume context, \emph{Electron. Trans. Numer. Anal.} \textbf{40} (2013), 94--119.
    
\bibitem{LSTV} M. Levi, F. Santagati, A. Tabacco, M. Vallarino, Poincar\'{e} inequalities on graphs, \emph{Anal. Math.} \textbf{49} (2023), no. 2, 529--544.
    
\bibitem{LW} P. Li, J. Wang, Weighted Poincar\'{e} inequality and rigidity of complete manifolds, \emph{Ann. Sci. \'{E}cole Norm. Sup. (4)} \textbf{39} (2006), no. 6, 921--982.
    
\bibitem{M} B. Mohar, Isoperimetric numbers of graphs, \emph{J. Combin. Theory Ser. B} \textbf{47} (1989), no. 3, 274--291. 
    
\bibitem{P} H. Poincar\'{e}, Sur les equations aux D\'{e}riv\'{e}es partielles de la physique math\'{e}matique, \emph{Amer. J. Math.} \textbf{12} (1890), no. 3, 211--294. 

\bibitem{ST} J. Soria, P. Tradacete, Best constants for the Hardy-Littlewood maximal operator on finite graphs, \emph{J. Math. Anal. Appl.} \textbf{436} (2016), no. 2, 661--682.
    
\bibitem{WH} Y. Wang, H. Huang, Eigenvalue estimates of the $p$-Laplacian on finite graphs, \emph{Differential Geom. Appl.} \textbf{74} (2021), Paper No. 101697, 11pp.
    
\bibitem{XZ} C. Xiang, F. Zhao, Sharp constants for fractional Hardy-Littlewood maximal operators on finite graphs, \emph{Forum Math.} \textbf{38} (2026), no. 2, 411--434.
    
\bibitem{Y} M. Yamasaki, Discrete potentials on an infinite network, \emph{Mem. Fac. Sci. Shimane Univ.} \textbf{13} (1979), 31--44.
    
\bibitem{Z} B. Zelinka, Medians and peripherians of trees, \emph{Arch. Math. (Brno)}  \textbf{4} (1968), 87--95.


\end{thebibliography}
\end{document}